\documentclass[11pt, a4paper]{amsart}

\usepackage[T1]{fontenc}
\usepackage[utf8]{inputenc}
\usepackage[english]{babel}
\usepackage[leqno]{amsmath}
\usepackage{amssymb}
\usepackage{tabularx}
\usepackage{a4wide}
\usepackage{url}
\usepackage{tikz-cd}
\usepackage{csquotes}
\usepackage{enumitem}
\usepackage{mathtools}
\usepackage{manfnt}
\usepackage{mathrsfs}
\usepackage{aurical}
\usepackage{bbding}
\usepackage{scalerel}

\usepackage[disable]{todonotes} 

\definecolor{cadmiumgreen}{rgb}{0.0, 0.42, 0.24}
\definecolor{darkred}{rgb}{.85,0,0}
\usepackage[
  colorlinks,
  linkcolor=darkred,
  citecolor=cadmiumgreen,
  pdfstartview ={FitV},
]{hyperref}
\hypersetup{
  pdfauthor={Omid Amini, Matthieu Piquerez},
  pdftitle={Hodge theory for tropical fans},
  bookmarksdepth = 2}

\usepackage[
  alphabetic,
  msc-links,
  nobysame,
  lite,
]{amsrefs}

\usetikzlibrary{calc}
\tikzset{vcenter/.style={baseline={([yshift=-.8ex]current bounding box.center)}}}
\tikzset{dot/.style={insert path={node {\tikz[baseline=.6pt]\filldraw[black] (0,0) circle (1.2pt);}}}}

\usetikzlibrary{cd, backgrounds, calc, arrows.meta, decorations.pathmorphing, intersections, 3d, fit}

\setlist[itemize,1]{itemsep=\smallskipamount}
\setlist[enumerate,1]{itemsep=\smallskipamount, label=\textnormal{(\arabic*)}}

\newtheorem{thm}{Theorem}[section]
\newtheorem{lemma}[thm]{Lemma}
\newtheorem{prop}[thm]{Proposition}
\newtheorem{claim}[thm]{Claim}
\newtheorem{cor}[thm]{Corollary}

\theoremstyle{definition}
\newtheorem{defii}[thm]{Definition}
\newenvironment{defi}
  {\pushQED{\qed}\defii}
  {\popQED\enddefii}
\newtheorem{remm}[thm]{Remark}
\newenvironment{remark}
  {\pushQED{\qed}\remm}
  {\popQED\endremm}
\newtheorem{exx}[thm]{Example}
\newenvironment{example}
  {\pushQED{\qed}\exx}
  {\popQED\endexx}
\newtheorem{question}[thm]{Question}

\numberwithin{equation}{section}

\DeclareFontFamily{U}{mathx}{}
\DeclareFontShape{U}{mathx}{m}{n}{<-> mathx10}{}
\DeclareSymbolFont{mathx}{U}{mathx}{m}{n}
\DeclareMathAccent{\widehat}{0}{mathx}{"70}
\DeclareMathAccent{\widecheck}{0}{mathx}{"71}

\DeclareRobustCommand\longtwoheadrightarrow
     {\relbar\joinrel\twoheadrightarrow}

\newcommand{\cf}{cf.~}
\newcommand{\ie}{i.e.}

\newcommand{\resp}{resp.\ }
\newcommand{\myand}{\ \textrm{and}\ }
\renewcommand{\~}{\widetilde}
\renewcommand{\hat}{\widehat}
\renewcommand{\check}{\widecheck}
\newcommand{\Q}{\mathbb{Q}}
\newcommand{\Z}{\mathbb{Z}}
\newcommand{\N}{\mathbb{N}}
\newcommand{\R}{\mathbb{R}}

\newcommand{\longhookrightarrow}{\lhook\joinrel\longrightarrow}

\let\oldchi\chi
\newcommand{\raisechi}[2]{\raisebox{.4ex}{$#1#2$}}
\renewcommand{\chi}{{\mathpalette\raisechi\oldchi}}

\let\oldforall\forall
\renewcommand{\forall}{\oldforall\:}
\let\oldbigwedge\bigwedge
\renewcommand{\bigwedge}{{\textstyle\oldbigwedge\!}}
\renewcommand{\emptyset}{\varnothing}
\renewcommand{\epsilon}{\varepsilon}
\renewcommand{\geq}{\geqslant}
\renewcommand{\leq}{\leqslant}
\renewcommand{\setminus}{\smallsetminus}

\makeatletter
\let\oldsum\sum
\renewcommand{\sum}{\@ifnextchar_\@mysum\oldsum}
\def\@mysum_#1{\oldsum_{\substack{#1}}}
\let\oldbigoplus\bigoplus
\renewcommand{\bigoplus}{\@ifnextchar_\@mybigoplus\oldbigoplus}
\def\@mybigoplus_#1{\oldbigoplus_{\substack{#1}}}
\let\oldprod\prod
\renewcommand{\prod}{\@ifnextchar_\@myprod\oldprod}
\def\@myprod_#1{\oldprod_{\substack{#1}}}
\makeatother

\newcommand{\rquot}[2]{#1\big/#2}
\newcommand{\rest}[1]{\raisebox{-1pt}{$\vert$}_{#1}}
\newcommand{\card}[1]{\vert#1\vert} 
\newcommand{\Card}[1]{\mathrm{card}(#1)} 
\newcommand{\simto}{\xrightarrow{\raisebox{-3pt}[0pt][0pt]{\small$\hspace{-1pt}\sim$}}}
\newcommand{\id}{\mathrm{id}} 
\newcommand{\dual}{\star}
\newcommand{\st}{\bigm|} 
\newcommand{\Bigst}{\Bigm|} 
\newcommand{\bul}{\bullet}

\DeclareMathOperator{\rk}{rk} 
\newcommand{\rank}{\mathrm{rank}}
\let\hom\relax
\DeclareMathOperator{\hom}{Hom} 

\newcommand{\Vect}{\mathrm{Vect}}
\DeclareMathOperator{\ord}{ord} 

\renewcommand{\i}{\mathrm i} 
\renewcommand{\k}{\Bbbk} 

\DeclareMathOperator{\PD}{PD} 

\renewcommand{\P}{\mathbb P} 

\newcommand{\conezero}{{\underline0}}
\newcommand{\e}{\mathfrak e} 
\newcommand{\nvect}{\mathfrak n} 

\newcommand{\shiftcomp}[2][0]{{}\mkern#1mu\overline{\mkern-#1mu#2}}
\newcommand{\comp}[1]{\if#1X \shiftcomp[3]{#1}\else\if#1Z \shiftcomp[3]{#1} \else \shiftcomp{#1}\fi\fi} 

\newcommand{\suppaux}[2]{\scalebox{1}[1.4]{$#1\lvert$}#2\scalebox{1}[1.4]{$#1\rvert$}}
  \newcommand{\supp}[1]{\mathpalette\suppaux{#1}}
\newcommand{\dimsaux}[2]{\raisebox{.2ex}{\scalebox{1}[.8]{$#1\lvert$}}#2\raisebox{.2ex}{\scalebox{1}[.8]{$#1\rvert$}}}
  \newcommand{\dims}[1]{\mathpalette\dimsaux{#1}}
\newcommand{\subfaceeq}{\preceq}
\newcommand{\supfaceeq}{\succeq}
\newcommand{\subface}{\prec}
\newcommand{\ssubface}{\mathbin{\mathchoice
  {\subface\!\!\!\cdot}%
  {\subface\!\!\!\cdot}%
  {\subface\!\cdot}%
  {\subface\!\cdot}%
}} 
\newcommand{\supface}{\succ}
\newcommand{\ssupface}{\mathbin{\mathchoice
  {\cdot\!\!\!\supface}%
  {\cdot\!\!\!\supface}%
  {\cdot\!\supface}%
  {\cdot\!\supface}%
}}

\newcommand{\LL}{\mathscr L} 
\newcommand{\eLL}{\hat\LL} 
\newcommand{\1}{{\mathbf 1}} 
\newcommand{\0}{{\mathbf 0}} 

\renewcommand{\div}{\mathrm{div}} 
\newcommand{\Div}{\mathrm{Div}} 
\newcommand{\Prin}{\mathrm{Prin}} 

\newcommand{\MW}{\mathrm{MW}} 
\newcommand{\W}{\mathrm{W}} 
\newcommand{\x}{\textnormal{x}}
\newcommand{\y}{\textnormal{y}}

\newcommand{\mT}{\mathcal T\!} 
\newcommand{\Csh}{\mathscr C} 
\newcommand{\Bsh}{\mathscr B} 
\newcommand{\Bsho}{\ssub{\mathcal B}!} 
\newcommand{\Ssh}{\mathscr S} 
\newcommand{\Sh}{\Ssh} 

\newcommand{\gst}[2]{\ssub{\langle#1\rangle}!_{#2}} 

\newcommand{\tropmod}[2]{{\mathcal{T\!M}}_{\!#1}(#2)} 
\newcommand{\basetm}[1]{#1_o} 
\newcommand{\symbuptm}{\tikz[scale=.2, baseline=-.1]{\draw(.1,-.05)to[out=50,in=-130](.9,.05) (.5,0)--++(0,.5);}}
\newcommand{\uptm}[1]{{#1}_{\!\symbuptm}} 

\newcommand{\etm}{\e_{\!\symbuptm}} 
\newcommand{\prtm}{\mathfrak p} 

\newcommand{\spec}{\mathrm{Spec}} 
\DeclareMathOperator{\gys}{Gys} 
\newcommand{\mer}{\mathcal M} 
\renewcommand{\O}{\mathcal O}

\newcommand{\filt}{{\mathscr G}}
\newcommand{\F}{\mathscr F}
\renewcommand{\H}{\mathscr H}

\newcommand{\Ma}{{\mathpalette\doscaleMa\relax}}
\newcommand{\doscaleMa}[2]{\scalebox{1.24}{$#1\mathfrak m$}}
\newcommand{\Fl}{\mathscr{F}} 
\newcommand{\Cl}{\mathcal C\!\ell} 
\newcommand{\prCl}{\mathcal C\!\ell_0} 
\newcommand{\crct}{\mathfrak C} 
\newcommand{\distel}{\ast} 
\newcommand{\bases}{\mathfrak B} 
\newcommand{\ind}{\mathfrak I} 
\newcommand{\rkm}{\mathrm{rk}} 
\newcommand{\contr}[1]{/#1} 
\newcommand{\del}{\setminus} 
\DeclareMathOperator{\cl}{cl} 

\let\cech\v
\renewcommand{\v}{\mathbf v} 
\newcommand{\Bl}[2]{\mathcal B\ell_{#2}(#1)} 
\newcommand{\cycl}{\mathrm{cl}} 

\newcommand{\Ploc}[1]{#1^{\scaleto{\bigstar}{1ex}}} 

\newcommand{\corps}{\mathbb K} 
\newcommand{\HL}{\mathrm{HL}} 
\newcommand{\HR}{\mathrm{HR}} 

\newcommand{\J}{\mathfrak{J}} 

\NewDocumentCommand{\ssub}{O{0pt} O{.8} m t! e{_^}}{
  #3%
  \IfValueT{#5}{
    \IfBooleanTF{#4}{\sb{\hspace{#1}\scaleobj{#2}{#5}}}{\sb{#5}}
  }
  \IfValueT{#6}{\sp{#6}}
}

\ExplSyntaxOn
\NewDocumentCommand{\tossub}{o o m}{
  \expandafter\let\csname old\cs_to_str:N #3\endcsname#3
  \renewcommand#3%
  {\ssub[#1][#2]{\csname old\cs_to_str:N #3\endcsname}}
}
\ExplSyntaxOff

\tossub\omega
\tossub\varpi
\tossub\deg
\tossub\LL
\tossub\eLL
\tossub\e
\tossub\distel
\tossub\rk
\tossub\P
\tossub[-2pt][1]\Gamma
\tossub\Div
\tossub\Prin
\tossub\cycl
\newcommand{\sssigma}{\ssub[-2pt]{\sigma}!}
\newcommand{\ssM}{\ssub{M}!}
\newcommand{\ssN}{\ssub{N}!}
\newcommand{\ssI}{\ssub{I}!}
\newcommand{\ssJ}{\ssub{J}!}
\newcommand{\ssSigma}{\ssub{\Sigma}!}
\newcommand{\sshatSigma}{\ssub{\hat\Sigma}!}
\newcommand{\sstildeSigma}{\ssub{\~\Sigma}!}

\def\xx{1}
\def\xy{0}
\def\yx{0}
\def\yy{1}
\def\zx{.3}
\def\zy{.6}

\def\ratio{1.3}
\def\ratioo{.5}
\def\ratioe{.6}

\newcommand{\BFcoordinates}{
  \coordinate (I) at (0,0);
  \coordinate (O) at (0,0);

  \coordinate (e1) at (\ratioe*\xx, \ratioe*\xy);
  \coordinate (e2) at (\ratioe*\yx, \ratioe*\yy);
  \coordinate (e3) at (\ratioe*\zx, \ratioe*\zy);
  \coordinate (e0) at ($(0,0)-(e1)-(e2)-(e3)$);

  \foreach \i/\j in {1/2, 1/3, 1/0, 2/3, 2/0, 3/0} {
    \coordinate (e\i\j) at ($(e\i)+(e\j)$);
  }

  \foreach \i/\j/\ax/\ay/\bx/\by in {1/2/\xx/\xy/\yx/\yy, 1/3/\xx/\xy/\zx/\zy, 2/3/\yx/\yy/\zx/\zy} {
    \pgftransformcm{\ax}{\ay}{\bx}{\by}{\pgfpoint{0}{0}}
    \coordinate (I\i\i\j) at (2,1);
    \coordinate (I\j\i\j) at (1,2);
    \coordinate (I\i\j) at ($(0,0)!(I\i\i\j)!(1,1)$);
    \coordinate (I\i) at ($(0,0)!(I\i\i\j)!(1,0)$);
    \coordinate (I\j) at ($(0,0)!(I\j\i\j)!(0,1)$);
  }

  \foreach \i/\ax/\ay in {1/\xx/\xy, 2/\yx/\yy, 3/\zx/\zy} {
    \pgftransformcm{\ax}{\ay}{-\xx-\yx-\zx+\ax}{-\xy-\yy-\zy+\ay}{\pgfpoint{0}{0}}
    \coordinate (I\i0) at (0,\ratio);
    \coordinate (I\i\i0) at (2,\ratio);
    \coordinate (I0\i0) at (-\ratioo,\ratio);
  }

  \coordinate (I0) at ($.333*(I010)+.333*(I020)+.333*(I030)$);
}

\def\roundface#1#2[#3]{\fill[#3] (I) -- (I#1) .. controls ($.3*(I#1)+.7*(I#1#1#2)$) and ($.3*(I#1#2)+.7*(I#1#1#2)$) .. (I#1#2) .. controls ($.3*(I#1#2)+.7*(I#2#1#2)$) and ($.3*(I#2)+.7*(I#2#1#2)$) .. (I#2) -- cycle}

\begin{document}

\allowdisplaybreaks

\title{Hodge theory for tropical fans}

\author{Omid Amini}
\address{CNRS - CMLS, École polytechnique, Institut polytechnique de Paris.}
\email{\href{omid.amini@polytechnique.edu}{omid.amini@polytechnique.edu}}

\author{Matthieu Piquerez}
\address{LS2N, Inria, Nantes Université}
\email{\href{matthieu.piquerez@univ-nantes.fr}{matthieu.piquerez@univ-nantes.fr}}

\date{\today}

\begin{abstract}
This paper is the first in a series devoted to the development of a Hodge theory for tropical varieties. We introduce a notion of $\mT$-stability for tropical fans and prove that various geometric properties of tropical fans are $\mT$-stable. As a consequence, we establish Kähler properties for the Chow ring in a large class of tropical fans, going beyond the case of matroids and their Bergman fans. As a by-product, we obtain a new proof of the Kähler package for combinatorial geometries. The approach makes it possible to deal with tropical fans with general weights.
\end{abstract}

\maketitle

\setcounter{tocdepth}{1}

\tableofcontents

\section{Introduction} \label{sec:intro}

The work in this paper is motivated by the recent development in combinatorial Hodge theory, which has undergone an expansion following the work of Adiprasito-Huh-Katz~\cite{AHK} and further follow-ups, leading to the resolution of several open problems in the theory of matroids and their applications. We refer to the expository papers by Huh~\cites{Huh18, Huh22}, Okounkov~\cite{Okou22}, Ardila~\cites{Ard18,Ard22} and Baker~\cite{Baker18} for an overview of these exciting developments in combinatorics.

\smallskip
Our starting point was the following question.
\begin{question}
How Chow rings of matroids and their Hodge theory is linked to homological properties of tropical varieties? To what extent the setting of combinatorial Hodge theory can be extended beyond the setting of matroids and their Bergman fans? Is there a Kähler geometry in the tropical setting?
\end{question}

Our aim in this work and its companions is to propose answers to these questions by developing Hodge theoretic aspects of tropical geometry.

\smallskip
The specific work presented in the current paper is of local nature and concerns geometric properties of tropical fans captured in their Chow rings. Tropical fans and their support named tropical fanfolds are building blocks for the construction of more general tropical varieties. We establish the Kähler package for the Chow ring of a large class of tropical fans which can be recursively constructed by using three basic operations that preserve the balancing condition (orientability in tropical geometry). As a by-product, we obtain a new proof of the Kähler package for combinatorial geometries that circumvents some of the main difficulties encountered in~\cite{AHK}. At the same time, our approach makes it possible to deal with tropical fans with more general weights.

\smallskip
The paper is written to be accessible with no previous knowledge of algebraic and tropical geometry. In particular, for algebraic geometric statements that play a crucial role in combinatorial and tropical Hodge theory we provide new combinatorial proofs.

\smallskip
In the rest of this introduction, we provide an overview of our results.

\subsection{Fans and their Chow rings}

Consider a lattice of finite rank $N \simeq \Z^n$ and let $N_\R$ be the vector space generated by $N$. A \emph{rational fan} in $N_\R$ is a non-empty collection $\Sigma$ of strongly convex rational polyhedral cones that verifies the following two properties:
\begin{enumerate}
\item if $\sigma$ is a cone in $\Sigma$, then any face $\tau$ of $\sigma$ belongs to $\Sigma$.
\item for a pair of cones $\sigma$ and $\eta$ in $\Sigma$, their intersection $\sigma \cap \eta$ is a common face of both $\sigma$ and $\eta$.
\end{enumerate}

We denote the support of $\Sigma$ by $\supp\Sigma$ and call it a \emph{fanfold}. For each integer $k$, $\Sigma_k$ denotes the set of $k$-dimensional cones in $\Sigma$.

Let $\Sigma$ be a rational fan that we assume to be simplicial, meaning that each cone in $\Sigma$ is generated by as many rays as its dimension. The \emph{Chow ring} of $\Sigma$ denoted by $A^\bul(\Sigma)$ is defined by generators and relations. Consider the polynomial ring $\Z[\x_\zeta]_{\zeta\in \Sigma_1}$ with indeterminate variables $\x_\zeta$ associated to rays $\zeta$ in $\Sigma_1$. Then, $A^\bul(\Sigma)$ is the quotient ring
\[ A^\bul(\Sigma) \coloneqq \rquot{\Z[\x_\zeta]_{\zeta\in \Sigma_1}}{\bigl(I + J\bigr)} \]
where
\begin{itemize}
\item $I$ is the ideal generated by the products $\x_{\rho_1}\!\cdots \x_{\rho_k}$, for $k\in \N$, such that $\rho_1, \dots, \rho_k$ are non-comparable rays in $\Sigma$, that is, they do not form a cone in $\Sigma$, and
\item $J$ is the ideal generated by the elements of the form
\[ \sum_{\zeta\in \Sigma_1} m(\e_\zeta)\x_\zeta, \qquad  m \in M \coloneqq N^\dual,\]
with $\e_\zeta$ the primitive vector of the ray $\zeta$.
\end{itemize}

The ideal $I+J$ is homogeneous and the Chow ring inherits a graded ring structure. Moreover, for degree larger than the dimension of $\Sigma$, the corresponding graded piece vanishes. Denoting by $d$ the dimension of $\Sigma$, we can thus write
\[ A^\bul(\Sigma) = \bigoplus_{k= 0}^d A^k(\Sigma) \]
with the $k$-th degree piece $A^k(\Sigma)$, $0\leq k \leq d$, generated as $\Z$-module by degree $k$ monomials. When the fan is unimodular, the Chow ring $A^\bullet(\Sigma)$ coincides with the Chow ring of the toric variety $\P!_\Sigma$ associated to $\Sigma$, see~\cite{Bri96} and~\cites{Dan78, BDP90, FS97}.

Similarly as above, we define the Chow ring with rational and real coefficients that we denote by $A^\bul(\Sigma, \Q)$ and $A^\bul(\Sigma, \R)$, respectively.

Chow rings of fans have a rich combinatorics. A discussion of these properties is provided in Section~\ref{sec:chow_ring}.

\medskip

\begin{center}
*\, *\, *
\end{center}

We fix some terminology before proceeding. Each cone $\sigma$ in a rational fan $\Sigma$ defines a sublattice of $N$ denoted by $N_\sigma$ that has rank equal to the dimension of $\sigma$. The sublattice $N_\sigma$ is given by the integral points of the vector space generated by $\sigma$. We denote by $\dims{\sigma}$ the dimension of a face $\sigma$. A cone $\sigma$ is called \emph{simplicial} if it is generated by $\dims{\sigma}$ rays, equivalently, if $\sigma$ is of the form $\sum_i \R_{\geq 0}\e_i$ for a family of independent vectors $\e_1, \dots, \e_{\dims{\sigma}}$ in $N_\R$. A rational cone $\sigma$ is called \emph{unimodular} if it is generated by $\dims{\sigma}$ vectors which form a basis of $N_\sigma$. A fan $\Sigma$ is called simplicial, \resp unimodular, if all its cones are simplicial, \resp unimodular. A \emph{facet} is a maximal cone. A fan is \emph{pure dimensional} if all its facets have the same dimension. We denote the dimension of $\Sigma$ by $d$.

In this paper, the star fan $\Sigma^\sigma$ of $\Sigma$ at a cone $\sigma \in \Sigma$ refers to the fan in $\rquot {N_\R}{N_{\sigma,\R}}$ induced by the cones $\eta$ in $\Sigma$ which contain $\sigma$ as a face.

\subsection{Tropical fans}

An \emph{orientation} of a rational fan $\Sigma$ of pure dimension $d$ is an integer valued map
\[ \omega\colon \Sigma_d \to \Z\setminus\{0\} \]
which verifies the so-called \emph{balancing condition}: for any cone $\tau$ in $\Sigma$ of codimension one, we have the vanishing of the following sum in the quotient lattice $\rquot{N}{N_\tau}$
\[ \sum_{\sigma \supset \tau} \omega(\sigma)\e_{\sigma}^{\tau} =0 \]
where the sum is over facets $\sigma$ of $\Sigma$ which contain $\tau$, and $\e_{\sigma}^{\tau}$ is the generator of the quotient $\rquot{(\sigma \cap N)}{(\tau \cap N)} \simeq \Z_{\geq 0}$. The balancing condition is the analog in polyhedral geometry of the orientability property for manifolds, and leads to the definition of a fundamental class that plays a central role in the treatment of Poincaré duality and other refined geometric properties in polyhedral geometry.

A \emph{tropical fan} in $N_\R$ is a pair $(\Sigma, \omega!_\Sigma)$ consisting of a pure dimensional rational fan $\Sigma$ and an orientation $\omega!_\Sigma$ as above. We will call \emph{tropical fanfold} the support of any tropical fan. We say that the tropical fan $(\Sigma, \omega!_\Sigma)$ is \emph{unitary} if $\omega!_\Sigma$ only takes values $\pm1$.

The class of tropical fans is closed under products (with the product orientation). Moreover, a star fan of a tropical fan is again tropical (with the induced orientation).

\subsection{Bergman fans of matroids}

Associated with any matroid $\Ma$ on a basis set $E$ is a fan $\ssSigma_{\Ma}$ called the \emph{Bergman fan of $\Ma$} that lives in the real vector space $\rquot{\R^E}{\R(1, \dots, 1)}$ and is quasi-projective and unimodular with respect to the lattice $\rquot{\Z^E}{\Z(1, \dots, 1)}$. The structure of the Bergman fan reflects the combinatorics of the matroid~\cite{AK06}. The \emph{augmented Bergman fan} of $\Ma$ is a unimodular fan $\ssub{\widehat\Sigma}!_\Ma$ in $\R^E$ defined in~\cite{BHMPW}. Both these fans are tropical with respect to the orientation which takes value one on each facet.

We will call \emph{Bergman fanfold} the support of a Bergman fan. We call \emph{generalized Bergman fan} any fan structure on a Bergman fanfold. Any augmented Bergman fan is a generalized Bergman fan. When the matroid is defined by an arrangement of hyperplanes, the Bergman fanfold can be identified with the tropicalization of the complement of the hyperplane arrangement, for coordinates given by the linear functions that define the arrangement~\cite{AK06}.

\subsection{Chow-Kähler tropical fans} \label{sec:local-intro}

A tropical fan $(\Sigma,\omega!_\Sigma)$ of dimension $d$ comes with a degree map defined by the orientation $\omega!_\Sigma$
\[\deg!_\Sigma \colon A^d(\Sigma) \to \Z\]
that leads to the pairing
\[ \begin{array}{ccc}
  A^k(\Sigma) \times A^{d-k}(\Sigma) &\to& \Z\\
  (\alpha, \beta)\quad\, &\mapsto& \deg!_\Sigma(\alpha \cdot \beta)
\end{array} \]
for any $k=0, \dots, d$. We say that $(\Sigma, \omega!_\Sigma)$ verifies \emph{Poincaré duality for the Chow ring with integer coefficients} denoted $\PD_\Z$ if the above pairing is perfect. We say $(\Sigma, \omega!_\Sigma)$ verifies \emph{Poincaré duality for the Chow ring with rational coefficients} denoted $\PD_\Q$ if the above pairing becomes perfect after tensoring with $\Q$.

For an element $\ell \in A^1(\Sigma,\Q)$, we say that the pair $(\Sigma, \ell)$ verifies the \emph{Hard Lefschetz property} denoted $\HL(\Sigma, \ell)$ if $(\Sigma,\omega!_\Sigma)$ verifies $\PD_\Q$ and the following holds:

\smallskip
\noindent (Hard Lefschetz) for any non-negative integer $k \leq \frac d2$, the multiplication map by $\ell^{d-2k}$ induces an isomorphism
\[ A^k(\Sigma, \Q) \simto A^{d-k}(\Sigma,\Q)\ , \qquad a \mapsto \ell^{d-2k}\cdot a. \]

\medskip

We say the pair $(\Sigma, \ell)$ verifies \emph{Hodge-Riemann bilinear relations} denoted $\HR(\Sigma, \ell)$ if $(\Sigma,\omega!_\Sigma)$ verifies $\PD_\Q$ and the following holds:

\smallskip
\noindent (Hodge-Riemann bilinear relations) for any non-negative integer $k \leq \frac d2$, the symmetric bilinear form
\[ Q_{\ell}^k\colon A^k(\Sigma,\Q) \times A^{k}(\Sigma,\Q) \to \Q\ ,\qquad (a,b) \mapsto (-1)^k\deg_\Sigma(\ell^{d-2k}\cdot a\cdot b) \]
is non-degenerate and its signature is given by the sum
\[ \sum_{i=0}^k (-1)^i\Bigl(\rank(A^i(\Sigma))-\rank(A^{i-1}(\Sigma))\Bigr). \]

We note that $\HR(\Sigma, \ell)$ implies $\HL(\Sigma, \ell)$. Moreover, denoting by $P_\ell^k(\Sigma)$ the \emph{primitive part} of $A^k(\Sigma,\Q)$ defined as the kernel of the multiplication map by $\ell^{d-2k+1}$ from $A^k(\Sigma,\Q)$ to $A^{d-k+1}(\Sigma,\Q)$, the property $\HR(\Sigma, \ell)$ becomes equivalent to requiring that $(-1)^kQ_\ell^k$ be positive definite on $P_\ell^k(\Sigma) \subseteq A^k(\Sigma,\Q)$.

\smallskip
Let now $(\Sigma,\omega!_\Sigma)$ be a unimodular tropical fan with positive weight function $\omega!_\Sigma$. An element~$\ell$ in $A^1(\Sigma,\Q)$ is called \emph{ample} if $\ell$ can be represented in the form $\sum_{\zeta \in \Sigma_1} f(\e_\zeta) \x_\zeta$ for a strictly convex conewise linear function $f$ on the fan $\Sigma$ which takes rational values on lattice points of~$\Sigma$ (see Section~\ref{sec:function_theory} for the definition of strict convexity). A fan which admits such a function is called \emph{quasi-projective}.

\smallskip
A tropical fan $(\Sigma, \omega!_\Sigma)$ is called \emph{Chow-Kähler} if it is quasi-projective and moreover, $(\Sigma, \omega!_{\Sigma})$ and, more generally, any star fan $(\Sigma^\sigma, \omega!_{\Sigma^\sigma})$, $\sigma\in\Sigma$, verifies the Hodge-Riemann bilinear relations $\HR(\Sigma^\sigma, \ell)$ for any ample element $\ell \in A^1(\Sigma^\sigma)$. We refer to Section~\ref{sec:forthcoming} for a discussion of the terminology.

\subsection{Divisor theory on tropical fans}

Let $(\Sigma, \omega!_{\Sigma})$ be a tropical fan of dimension $d$. A \emph{meromorphic function} on $\Sigma$ is by definition a continuous conewise integral linear function on $\supp\Sigma$. We denote by $\mer(\Sigma, \omega!_\Sigma)$ the set of all meromorphic functions on $\Sigma$.

To any face $\tau$ of codimension one in $\Sigma$, we associate the \emph{order of vanishing function} $\ord_\tau \colon \mer(\Sigma, \omega!_{\Sigma}) \to \Z$ defined at any $f\in \mer(\Sigma, \omega!_{\Sigma})$ by the \emph{sum of slopes of $f$ along the adjacent facets}, see~\cite{AR10} or Section~\ref{sec:divisors} for the precise definition. The \emph{divisor of a meromorphic function $f$} is by definition the pair $\div(f) = (\Delta, \omega!_{\Delta})$ consisting of the subfan $\Delta$ of~$\Sigma$ defined by the set of all codimension one faces $\tau$ of $\Sigma$ with $\ord_\tau(f) \neq 0$, and the map $\omega!_{\Delta} \colon \Delta_{d-1}\to \Z \setminus \{0\}$ which takes value $\ord_\tau(f)$ at any face $\tau \in \Delta_{d-1}$. We say that $\div(f)$ is \emph{trivial} if $\ord_\tau(f) = 0 $ for all $\tau \in \Sigma_{d-1}$. The divisor $\div(f)$ is a tropical fan of dimension $d-1$ provided that it is non-trivial.

\smallskip
We call $f\in \mer(\Sigma, \omega!_{\Sigma})$ \emph{holomorphic on $\Sigma$} provided that $\ord_\tau(f) \geq 0$ for all $\tau \in \Sigma_{d-1}$.

\smallskip
A \emph{divisor} $D$ on $(\Sigma, \omega!_{\Sigma})$ is a tropical fan $(\Delta, \omega!_{\Delta})$ of dimension $d-1$ with $\Delta$ a subfan of $\Sigma$. Divisors of the form $\div(f)$ for $f\in \mer(\Sigma, \omega!_{\Sigma})$ are called \emph{principal}. We say $D$ is \emph{$\Q$-principal} if an integer multiple $aD$ with $a$ non-zero is principal.

\subsection{$\mT$-stability} \label{subsec:intro_T_stability}

We introduce a notion of \emph{$\mT$-stability} for tropical fans and their geometric properties. In practice, this allows to proceed by induction and reduce to the simplest possible tropical fans. We give an idea here and refer to Section~\ref{sec:operations} for more details.

The definition is based on three types of operations on tropical fans: \emph{products}, \emph{stellar subdivisions} and their inverse \emph{stellar assemblies}, and \emph{tropical modifications}. The first two operations are classical in the theory of fans. Stellar subdivision in particular corresponds to the fundamental notion of blow-up in algebraic geometry. The third operation, revealed in the pioneering work by Mikhalkin~\cites{Mik06, Mik07}, is specific to the tropical setting and its importance lies in the possibility of producing richer tropicalizations out of the existing ones by introducing \emph{new coordinates}. Given a meromorphic function $f$ on a tropical fan $(\Sigma, \omega!_\Sigma)$, the tropical modification of $(\Sigma, \omega!_\Sigma)$ along (the divisor of) $f$ is the result of modifying the graph of $f$ into a tropical fan $(\~\Sigma, \omega!_{\~\Sigma})$ by introducing new cones lying above $\div(f)$. We refer to Section~\ref{subsec:tropical_modification} for the definition.

Let $\Csh$ be a class of tropical fans and $\Sh\subseteq \Csh$ a subclass. We say that \emph{$\Sh$ is $\mT$-stable in $\Csh$}, or simply \emph{$\mT$-stable} in the case $\Csh$ is the class of all tropical fans, if the following properties hold:
\begin{itemize}
\item (Stability under products) For a pair of tropical fans $(\Sigma, \omega!_\Sigma)$ and $(\Sigma', \omega!_{\Sigma'})$ in $\Sh$, if the product $(\Sigma \times \Sigma', \omega!_{\Sigma\times \Sigma'})$ belongs to $\Csh$, then it is in $\Sh$.

\item (Stability under tropical modifications along a divisor in the subclass) Given a tropical fan $(\Sigma, \omega!_\Sigma)$ in $\Sh$ and a meromorphic function $f$ on $\Sigma$ such that the divisor $\div(f)$ of $f$ is in $\Sh$, the tropical modification $(\~\Sigma, \omega!_{\~\Sigma})$ of $\Sigma$ along (the divisor of) $f$ is in $\Sh$ provided that it belongs to the class $\Csh$.

\item (Stability under stellar subdivisions and stellar assemblies with center in the subclass) For an element $(\Sigma, \omega!_\Sigma)$ in $\Csh$, and for a cone $\sigma \in \Sigma$ which has the property that the stellar subdivision $\Sigma'$ of $\Sigma$ at $\sigma$ belongs to $\Csh$ and the star fan $\Sigma^\sigma$ belongs to $\Sh$, we have $(\Sigma, \omega!_{\Sigma}) \in \Sh$ if and only if $(\Sigma', \omega!_{\Sigma'}) \in \Sh$.
\end{itemize}

Examples of classes $\Csh$ of tropical fans which are of interest to us are \emph{all}, \resp \emph{simplicial}, \resp \emph{unimodular}, \resp \emph{quasi-projective}, \resp \emph{principal}, \resp \emph{$\Q$-principal}, \resp \emph{div-faithful}, \resp \emph{locally irreducible} and \resp \emph{$\Q$-locally irreducible} tropical fans. The last five classes and their properties are introduced and studied later in the paper, see Section~\ref{sec:principal_divfaithful_locallyirreducible-intro} for a discussion.

The subclass $\Csh\subseteq \Csh$ is obviously $\mT$-stable in $\Csh$, and it is easy to see that intersection of two $\mT$-stable classes is again $\mT$-stable. This leads to the following definition.

Let $\Bsh$ be a subset of $\Csh$. We refer to $\Bsh$ as the \emph{base set}. The \emph{$\mT$-stable subclass of\/ $\Csh$ generated by $\Bsh$} denoted by $\gst{\Bsh}{\Csh}$ is by definition the smallest subclass of $\Csh$ which contains $\Bsh$ and which is $\mT$-stable; it is obtained by taking the intersection of all $\mT$-stable subclasses $\Sh$ of $\Csh$ that contain $\Bsh$. If $\Csh$ is the class of all tropical fans, we just write $\gst{\Bsh}{}$.

An important example of the base set is the set $\Bsho$ defined as follows. Denote by $\conezero$ the cone~$\{0\}$. By an abuse of the notation, we denote by $\conezero$ the fan consisting of unique cone $\conezero$. Denote by $(\conezero, n)$, $n\in\Z\setminus\{0\}$, the tropical fan of dimension $0$ with weight function taking value $n$ on $\conezero$. Let $\Lambda$ be the complete tropical fan in $\R$ that consists of three cones $\conezero, \R_{\geq 0}$, and $\R_{\leq 0}$ endowed with constant weight equal to $1$ on $\R_{\geq 0}$ and $\R_{\leq 0}$. We set
\[\Bsho\coloneqq \left\{ (\conezero,n), n\in\Z\setminus\{0\} \right\} \cup \left\{ \Lambda \right\}.\]
As we show later, already in this very simple case, $\gst{\Bsho}{}$ contains many interesting fans. For instance, generalized Bergman fans are all in this class, but $\gst{\Bsho}{}$ is strictly larger.

A tropical fan is called \emph{quasilinear} if it belongs to $\gst{\Bsho}{}$. By definition, this means that the tropical fan can be obtained from the collection $\Bsho$ by performing a sequence of three above operations on tropical fans. More precisely, a tropical fan is quasilinear if it either belongs to~$\Bsho$, or, is a product of two quasilinear fans, or, is the stellar subdivision or stellar assembly of a quasilinear fan, or, is obtained as a result of tropical modification of a quasilinear fan along a tropical divisor which is itself quasilinear.

\smallskip
$\mT$-stability can be defined for properties of tropical fans. If $P$ is a predicate on tropical fans and $\Csh$ is a class of tropical fans, then $P$ is called \emph{$\mT$-stable in $\Csh$} if the subclass of tropical fans in $\Csh$ which verify $P$ is $\mT$-stable in $\Csh$.

\begin{remark}
In earlier version of our work, we were using \emph{tropical shellability} instead of \emph{$\mT$-stability}. In Appendix~\ref{sec:multimagmoid}, we introduce an algebraic framework that gives a conceptual formulation of the notion of $\mT$-stability, and justifies the naming. The terminology \emph{quasilinear} instead of our former \emph{tropically shellable} is borrowed from the work of Nolan Schock~\cite{Sch21} which applies our results in the study of tropical compactifications of moduli spaces.
\end{remark}

\subsection{$\mT$-stability results} Various geometric properties of tropical fans are shown to be $\mT$-stable within an appropriate class.

\subsubsection{Normality, local irreducibility, and div-faithfulness} \label{sec:principal_divfaithful_locallyirreducible-intro}

A tropical fan $(\Sigma, \omega!_{\Sigma})$ of dimension $d$ is called \emph{normal}, \resp \emph{$\Q$-normal}, if for any cone $\tau \in \Sigma_{d-1}$, any orientation of $\Sigma^\tau$ is an integer, \resp rational, multiple of the orientation $\omega!_{\Sigma^\tau}$ induced from $\omega!_{\Sigma}$. The following theorem is proved in Section~\ref{sec:normal_stable}.

\begin{thm} \label{thm:normal_irreducible_stable-intro}
Being $\Q$-normal is $\mT$-stable. Being normal is $\mT$-stable in the class of unitary tropical fans.
\end{thm}

A tropical fan $(\Sigma, \omega!_{\Sigma})$ is called \emph{irreducible, \resp $\Q$-irreducible, at a face $\eta \in \Sigma$} provided that any orientation of $\Sigma^\eta$ is an integer, \resp rational, multiple of the orientation $\omega!_{\Sigma^\eta}$ induced from $\omega!_\Sigma$. We call $(\Sigma, \omega!_\Sigma)$ \emph{locally irreducible}, \resp \emph{$\Q$-locally irreducible}, if it is irreducible, \resp $\Q$-irreducible, at any face $\eta\in\Sigma$. The following theorem is proved in Section~\ref{sec:local_irreducibility_stable}.

\begin{thm}
$\Q$-local irreducibility is $\mT$-stable. Local irreducibility is $\mT$-stable in the class of unitary tropical fans.
\end{thm}

We say that a fan $\Sigma$ is \emph{divisorially faithful} or simply, \emph{div-faithful}, \emph{at a face $\eta$} if the following holds: For any meromorphic function $f$ on $\Sigma^\eta$, if $\div(f)$ is trivial, then $f$ is a linear function on $\Sigma^\eta$. We call the tropical fan $\Sigma$ \emph{div-faithful} if $\Sigma$ is div-faithful at any cone $\eta\in \Sigma$.

\smallskip
The next theorem is proved in Section~\ref{sec:stability_div-faithful}.

\begin{thm} \label{thm:stability_div-faithful-intro}
The property of being div-faithful is $\mT$-stable.
\end{thm}

Div-faithful property plays an important role in our treatment of Kähler geometry for tropical fans as Chow rings behave nicely under tropical modifications for div-faithful tropical fans, see Section~\ref{sec:tropmod-intro}.

\subsubsection{Principality} We say that a tropical fan $(\Sigma, \omega!_\Sigma)$ is \emph{principal at $\eta$} if any divisor on $\Sigma^\eta$ is the divisor of a meromorphic function on $(\Sigma^\eta, \omega!_{\Sigma^\eta})$. We call the tropical fan $(\Sigma,\omega!_\Sigma)$ \emph{principal} if $\Sigma$ is principal at any cone $\eta \in \Sigma$. We say $(\Sigma,\omega!_\Sigma)$ is \emph{$\Q$-principal} if for any divisor $D$, an integer multiple $aD$ with $a \in \Z \setminus \{0\}$ is principal.

 In Section~\ref{sec:stability_principality}, we prove the following theorem.

\begin{thm} \label{thm:stability_principal-intro}
Being $\Q$-principal is $\mT$-stable within the class of\/ $\Q$-locally irreducible tropical fans. Being principal is $\mT$-stable within the class of locally irreducible and unitary tropical fans.
\end{thm}

\subsubsection{Poincaré duality for the Chow ring} We have the following theorem whose proof is given in Section~\ref{sec:PD_stable}.

\begin{thm} \label{thm:PD_stable-intro}
Properties $\PD_\Z$ and $\PD_\Q$ are both $\mT$-stable in the class of div-faithful unimodular tropical fans.
\end{thm}

In particular, we deduce the following result.

\begin{thm}
Any unimodular quasilinear fan verifies $\PD_\Q$. Any unitary unimodular quasilinear fan verifies $\PD_\Z$.
\end{thm}

\subsubsection{Chow-Kähler property}

In Section~\ref{sec:kahler}, we establish the following theorem.
We say a tropical fan $(\Sigma, \omega!_\Sigma)$ is \emph{effective} provided that the orientation $\omega!_\Sigma$ takes positive values.

\begin{thm} \label{thm:chow-KP_stable-intro}
The property of being Chow-Kähler is $\mT$-stable in the class of effective quasi-projective unimodular tropical fans.
\end{thm}

An immediate corollary of this theorem is the following result.

\begin{thm}\label{thm:chow-KP_quasilinear-intro}
Any effective unimodular quasi-projective quasilinear fan is Chow-Kähler.
\end{thm}

\subsection{Three key ingredients in the proofs}

In establishing the above $\mT$-stability results, we need to describe the behavior of the Chow rings under tropical modifications and stellar subdivisions, and give a Chow-theoretic description of the introduced geometric notions. We use three main ingredients in doing so. We briefly discuss them here.

\subsubsection{Chow ring of a tropical modification}\label{sec:tropmod-intro} Section~\ref{sec:chow_tropmod} is devoted to the study of the behavior of Chow rings under tropical modifications. We prove the following important result.

\begin{thm}[Stability of the Chow ring under tropical modifications] \label{thm:invariant_chow_tropical_modification-intro}
Let $(\~\Sigma, \omega!_{\~\Sigma})$ be the tropical modification of the tropical fan $(\Sigma, \omega!_{\Sigma})$ along (the divisor of) a meromorphic function~$f$. Assume that $\Sigma$ is div-faithful. We have an isomorphism
\[ A^\bul(\~\Sigma,\Q) \simeq A^\bul(\Sigma,\Q) \]
and an isomorphism between Minkowski weights
\[ \MW_{\bul}(\~\Sigma, \Z) \simeq \MW_{\bul}(\Sigma, \Z). \]
If in addition, either $\Sigma$ is saturated or $A^\bul(\Sigma)$ torsion-free, then we have an isomorphism between the Chow rings with integral coefficients
\[ A^\bul(\~\Sigma) \simeq A^\bul(\Sigma). \]
The isomorphisms are all compatible with the degree maps.
\end{thm}

For the definition of Minkowski weights, see Section \ref{sec:MW}. The div-faithfulness is needed in the theorem, see Example~\ref{ex:necessity_invariant_chow_tropical_modification}. Without this assumption, we can prove that there is always a surjection from the Chow ring of $\Sigma$ to that of $\~\Sigma$. Saturation property is discussed in Section~\ref{sec:saturation}.

\subsubsection{Keel's lemma} Consider a unimodular fan $\Sigma$ of dimension $d$ and let $\sigma$ be a cone in~$\Sigma$. Let $\Sigma'$ be the fan obtained by the unimodular stellar subdivision of $\sigma$ in $\Sigma$. Denote by $A^\bul(\Sigma)[T]$ the polynomial ring in one variable $T$ over the coefficient ring $A^\bul(\Sigma)$. For each ray $\zeta \subface \sigma$, let $x_\zeta$ be the element in $A^1(\Sigma)$ given by the generator $\x_\zeta$.

\begin{thm}[Keel's lemma] \label{thm:keel-intro}
We have an isomorphism
\[ A^\bul(\Sigma')\simeq \rquot{A^\bul(\Sigma)[T]}{(\J T+P(T))} \]
where $\J$ is the kernel of the restriction map from $A^\bul(\Sigma)$ to $A^\bul(\Sigma^\sigma)$ and $P(T)$ is the product $\prod_{\zeta}(x_\zeta+T)$ over rays $\zeta$ of $\sigma$.
\end{thm}

This is a special case of a more general statement proved in~\cite{Kee92}*{Theorem 1 in the appendix} on Chow rings of blow-ups of algebraic varieties. We give a combinatorial proof of this key result in Section~\ref{sec:keel_proof}.

\subsubsection{Localization lemma} Let $\Sigma$ be a simplicial rational fan $\Sigma$ in $N_\R$ and consider the Chow ring $A^\bul(\Sigma)$. Localization lemma provides an alternative presentation of each homogeneous piece $A^k(\Sigma)$ of the Chow ring. This turns out to be quite powerful in combinatorial treatment of Chow rings.

For each cone $\sigma$ of dimension $k$ in $\Sigma$, let $\x_\sigma$ be an indeterminate variable and let $Z^k(\Sigma)$ be the free abelian group generated by $\x_\sigma$, $\sigma\in \Sigma_k$. Sending $\x_\sigma$ to the product $\prod_{\rho}\x_\rho$, $\rho$ a ray of $\sigma$, gives an embedding $Z^k(\Sigma) \hookrightarrow \Z[\x_\rho \mid \rho\in \Sigma_1]$. Passing to the quotient by the ideal $I+J$ gives a map $Z^k(\Sigma) \to A^k(\Sigma)$.

\begin{thm}[Localization lemma] \label{thm:kernel_Z^k-intro}
The map $Z^k(\Sigma) \to A^k(\Sigma)$ is surjective and its kernel is generated by elements of the form
\[ \sum_{\sigma\supseteq\tau \\ \sigma\in\Sigma_k} m(\e_{\sigma}^{\tau})\x_\sigma \]
for $\tau$ a face in $\Sigma_{k-1}$ and $m$ an element of $M$ which vanishes on $N_\tau$.
\end{thm}

Note that since $m$ is vanishing on $N_\tau$, it defines a linear map $m \colon \rquot{N}{N_\tau} \to \Z$, and so $m(\e_{\sigma}^{\tau})$ is well-defined.

This important result is a special case of a more general statement proved in~\cite{FMSS}*{Theorem 1}. We will provide a combinatorial proof in Section \ref{sec:kernel_Z^k}.

\subsection{Sketch of the proof of the Kähler package for combinatorial geometries}

A special case of our Theorem~\ref{thm:chow-KP_stable-intro} is a new proof of the Kähler property for Chow rings of matroids, established in~\cite{AHK}. In order to describe the geometric content of our $\mT$-stability results, we provide a sketch of our proof of the Kähler property for matroids if we were to rewrite it avoiding the use of $\mT$-stability. Note that this proof applies equally to any effective unimodular quasilinear fan, and leads to Theorem~\ref{thm:chow-KP_quasilinear-intro}.

Consider a Bergman fan $\ssSigma_\Ma$ of dimension $d$ and let $X = \supp{\ssSigma_\Ma}$ be the corresponding Bergman fanfold. By a star fanfold of $X$ we mean the fanfold of any star fan $\ssSigma_\Ma^\sigma$ for any cone $\sigma \in \ssSigma_\Ma$. All the star fanfolds of $X$ are Bergman. We prove that any quasi-projective unimodular generalized Bergman fan with support $X$ is Chow-Kähler.

The proof goes as follows.

\smallskip
\noindent $(1)$ A star fan of a quasi-projective unimodular generalized Bergman fan is itself a quasi-projective unimodular generalized Bergman fan.

\smallskip
\noindent $(2)$ Proceeding by induction on the rank of $\Ma$, we can assume that all star fans $\Sigma^\sigma$, $\sigma \neq \conezero$, of any quasi-projective unimodular generalized Bergman fan $\Sigma$ with support $X$ are Chow-Kähler. Using Keel's lemma and weak factorization theorem~\cites{Wlo97, Mor96}, we deduce that it will be enough to show the existence of one Chow-Kähler tropical fan with support $X$.

\smallskip
\noindent $(3)$ If the matroid $\Ma$ is free, the Bergman fanfold $X$ has support the entire space. In this case, the product of $d$ copies of the projective line $\Lambda$ (introduced in Section \ref{subsec:intro_T_stability}) gives a Chow-Kähler fan with support $X$. Otherwise, there exists an element $e$ of $\Ma$ such that the deletion $\Ma\setminus e$ has the same rank as $\Ma$. In this case, $\ssSigma_{\Ma/e}$ is a tropical divisor in $\ssSigma_{\Ma \setminus e}$. Proceeding by a second induction on the size of the ground set, we can suppose that $\ssSigma_{\Ma \setminus e}$ is Chow-Kähler. By the Chow-ring characterization of principal and div-faithfulness given in Section~\ref{subsec:characterization_principal_div-faithful}, it follows $\ssSigma_{\Ma \setminus e}$ is both div-faithful and principal. We infer that $\ssSigma_{\Ma/e}$ is the divisor of a holomorphic function on $\ssSigma_{\Ma \setminus e}$. This function is moreover unique up to addition by a linear function. The tropical modification $\sstildeSigma_{\Ma \setminus e}$ of $\ssSigma_{\Ma \setminus e}$ along $\ssSigma_{\Ma/e}$ is well-defined and turns out to have support equal to the Bergman fanfold $X$, a result proved by Shaw in their PhD thesis~\cite{Sha13a}, see Lemma~\ref{lem:shaw}. By our Theorem~\ref{thm:invariant_chow_tropical_modification-intro}, the Chow ring of $\sstildeSigma_{\Ma \setminus e}$ is isomorphic to the Chow ring of $\ssSigma_{\Ma \setminus e}$. We infer the existence of a Chow-Kähler fan with support $X$. By the previous point, we conclude.

\subsection{Examples}

We include in the last section of this paper a collection of examples dealing with various aspects of the geometry of tropical fans, to which we refer in the other sections. These examples are supplemented by further questions and remarks that, we hope, will clarify the concepts introduced in the paper.

\subsection{Related work}

Kähler package for the Chow ring of generalized Bergman fans was established by Ardila, Denham and Huh in~\cite{ADH} using weak factorization theorem and the work~\cite{AHK}. An alternative proof of the Chow-Kähler property for matroids using semismall decompositions is provided in the work by Braden, Huh, Matherne, Proudfoot, and Wang~\cites{BHMPW, BHMPW20b}. Parallel to these two works, a proof of the Kähler package for the Chow ring of generalized Bergman fans along the above sketched lines appeared in the first version of our work~\cite{AP-tht}.

Poincaré duality for the Chow ring of generalized Bergman fans is also proved using weak factorization in the work by Gross and Shokrieh~\cite{GS19}. Chow-Kähler property for the degree one part of the Chow ring of matroids (used in applications to log-concavity statements) is established in the work by Backman, Eur and Simpson~\cite{BES} by using a simplicial representation of the Chow ring of matroids.

We refer to the survey papers~\cites{Ard18, Baker18, Huh18, Ard22, Huh22, Okou22} for an overview of the current developments of combinatorial Hodge theory.

\subsection{Forthcoming work}\label{sec:forthcoming} In our companion work, we will provide several applications of the results of this paper.

\smallskip
\noindent -- We prove in~\cite{AP-him} a precise link between Chow rings of unimodular fans and the tropical cohomology rings of their canonical compactifications. In the case of matroids, when the matroid $\Ma$ is realizable over a field, by the work of Feichtner and Yuzvinsky~\cite{FY04}, the Chow ring $A^\bul(\ssSigma_{\Ma})$ is the Chow ring of a smooth projective variety over the same field. Moreover, as it is noted in~\cite{AHK}*{Theorem 5.12}, a Chow equivalence between a smooth projective variety over a field and the toric variety $\P!_{\ssSigma_{\Ma}}$ implies that the matroid $\Ma$ is realizable over that field. In this regard, it came somehow as a surprise that in the non-realizable case, the Chow ring $A^\bul(\ssSigma_{\Ma})$, which is the Chow ring of a non-complete smooth toric variety, verifies all the nice properties enjoyed by the cohomology rings of complex projective manifolds. That said, the fact that any matroid is realizable over the tropical hyperfield (see the recent work of Baker and Bowler~\cite{BB19}), suggested that a similar geometric picture to the realizable case arise in the general situation. The Hodge isomorphism theorem proved in~\cite{AP-him} confirms this by showing that the ring $A^\bul(\ssSigma_{\Ma})$ is still the cohomology ring of a smooth projective tropical variety. The result gives as well an alternative representation of the Chow rings of matroids (and more general quasilinear tropical fans), adding a tropical viewpoint to the work of Feichtner and Yuzvinsky~\cite{FY04}, Brion~\cite{Bri96}, and Billera~\cite{Bil89}. We use these results to establish a tropical analog of Kleiman's criterion of ampleness.

\smallskip
\noindent -- Homological properties of tropical fans are further studied in our paper~\cite{AP-homology}. We establish an analogue of the Deligne weight spectral sequence for tropical fans, relating the cohomology of the fan to the Chow rings of its star fans.

\smallskip
\noindent -- Tropical Hodge theory in the global setting is developed in our work~\cites{AP-tht}. Using Chow rings of tropical fans, we introduce Kähler tropical varieties and establish a Hodge theory for them. Note that we use the terminology Chow-Kähler in this paper because in the development of Kähler geometry for tropical varieties, a tropical fan is called Kähler if in addition to being Chow-Kähler, its corresponding tropical fanfold is a tropical homology manifold (\ie, the tropical cohomology of any open subset verifies Poincaré duality).

\smallskip
\noindent -- In~\cite{AP20-HC}, we prove a tropical analogue of the Hodge conjecture for Kähler tropical varieties which admit a rational triangulation. We moreover prove the Grothendieck standard conjecture that numerical and homological equivalences coincide for tropical varieties.

\smallskip
\noindent -- In~\cite{AP23-MA}, we use the set-up of this paper to properly formulate and study the tropical Monge-Ampère equation on tropical varieties. This is motivated by the work of Yang Li~\cite{Li20} which reduces the SYZ conjecture in maximally degenerate families of complex algebraic varieties to the existence of solutions to a tropical Monge-Ampère equation (once this has been properly formulated).

\smallskip
\noindent -- Although at some occasions in the paper we assume that tropical fans are unimodular, with a little extra effort, our approach can be generalized to deal with rational simplicial fans when working with Chow rings with rational coefficients, and more general simplicial fans when working with Chow rings with real coefficients. Part of this generalization is discussed in the first chapter of the second named author's doctoral thesis~\cite{Piq-thesis}. A full treatment of the topic in this paper would have resulted in an increase in the length and technicality of the article. For this reason, and to simplify presentation, we defer discussion of these results to a future publication.

\subsection*{Basic notations} \label{sec:intro-basic-notations}

The set of natural numbers is denoted by $\N = \{1, 2, 3, \dots\}$. For any natural number $n$, we denote by $[n]$ the set $\{1,\dots, n\}$.

\smallskip
The set of non-negative real numbers is denoted by $\R_{\geq 0}$.

\smallskip
For a lattice $N$, we view its dual $M=N^\dual$ as linear forms on $N$.

\smallskip
Given a poset $(P, \subfaceeq)$ and a functor $\phi$ from $P$ to a category $\mathcal C$, if $\phi$ is covariant (\resp contravariant), then for a pair of elements $\tau \subfaceeq \sigma$ in $P$, we denote by $\phi_{\tau \subfaceeq \sigma}$ (\resp $\phi_{\sigma \supfaceeq \tau}$), the corresponding map $\phi(\tau) \to \phi(\sigma)$ (\resp $\phi(\sigma) \to \phi(\tau)$) in $\mathcal C$, the idea being that in the subscript of the map $\phi_{\bul}$ representing the arrow in $\mathcal C$, the first item refers to the source and the second to the target. This convention will be in particular applied to the poset of faces in a fan.

For subsets $A$ and $B$ of a real vector space $V$, we write $A+B$ for the subset of $V$ consisting of all the sums $a+b$ with $a\in A$ and $b\in B$.

\subsection*{Acknowledgments}

We thank the organizers and participants of the Banff workshop on Algebraic Aspects of Matroid Theory (23w5149) who suggested a change in the terminology regarding our earlier use of tropical shellability instead of $\mT$-stability.

Content of this paper was part of a course taught by one of us at the Berlin Mathematical School during the academic year 2022-2023 and a minicourse at EPFL Bernoulli center. We thank the organizers and participants of these events for their constructive questions and remarks. We warmly thank Edvard Aksnes and Kris Shaw for discussions and collaboration related to the subject of this paper.

O.A. is part of the ANR project ANR-18-CE40-0009, and thanks Math+, the Berlin Mathematics Research Center, for support. M.P. has received funding from the European Research Council (ERC) under the European Union’s Horizon 2020 research and innovation program (grant agreement No. 101001995).

\section{Preliminaries} \label{sec:prel}

The aim of this section is to introduce basic notations and definitions which will be used all through the paper.

\smallskip
Throughout, $N$ will be a free $\Z$-module of finite rank and $M=N^\dual = \hom(N, \Z)$ will be the dual of $N$. We denote by $N_\Q$, $N_\R$, $M_\Q$, $M_\R$ the corresponding rational and real vector spaces. We thus have $M_\Q = N_\Q^\dual$ and $M_\R = N_\R^\dual$. For a polyhedral cone $\sigma$ in $N_\R$, we use the notation $N_{\sigma, \R}$ to denote the real vector subspace of $N_\R$ generated by elements of $\sigma$ and set $N^\sigma_\R \coloneqq \rquot{N_\R}{N_{\sigma, \R}}$. If the cone $\sigma$ is rational, we get natural lattices $N^\sigma$ and $N_\sigma$ in $N^\sigma_\R$ and $N_{\sigma, \R}$, respectively, which are both of full rank. The duals of $N^\sigma$ and $N_\sigma$ are denoted by $M^\sigma$ and $M_\sigma$, respectively.

For the ease of reading, we adopt the following convention. We use $\sigma$ (or any other face of $\Sigma$) as a superscript where referring to the quotient of some space by $N_{\sigma, \R}$ or to the elements related to this quotient. In contrast, we use $\sigma$ as a subscript for subspaces of $N_{\sigma,\R}$ or for elements associated to these subspaces.

\subsection{Fans} \label{subsec:fans}

Let $\Sigma$ be a fan of dimension $d$ in $N_\R$. The \emph{dimension} of a cone $\sigma$ in $\Sigma$ is denoted by $\dims\sigma$. The set of $k$-dimensional cones of $\Sigma$ is denoted by $\Sigma_k$, and elements of $\Sigma_1$ are called \emph{rays}. We denote by $\conezero$ the cone $\{0\}$. Any $k$-dimensional cone $\sigma$ in $\Sigma$ is determined by its set of rays in $\Sigma_1$. The \emph{support} of $\Sigma$ denoted $\supp \Sigma$ is the closed subset of $N_\R$ obtained by taking the union of all the cones in $\Sigma$. A fan $\Sigma$ with $\supp\Sigma = N_\R$ is called \emph{complete}. A \emph{facet} of $\Sigma$ is a cone which is maximal for the inclusion. $\Sigma$ is \emph{pure dimensional} if all its facets have the same dimension. The \emph{$k$-skeleton of\/ $\Sigma$} is by definition the subfan of $\Sigma$ consisting of all the cones of dimension at most $k$.

\subsection{Face poset}

For a fan $\Sigma$, we denote by $\LL!_\Sigma$ the face poset of $\Sigma$ in which the partial order $\subfaceeq$ is given by the inclusion of faces: we write $\tau \subfaceeq \sigma$ if $\tau \subseteq \sigma$ and $\tau \subface \sigma$ if $\tau \subset \sigma$. Set $\0 \coloneqq \conezero$. The \emph{extended poset $\eLL!_\Sigma$} is defined as $\eLL!_\Sigma \coloneqq \LL!_\Sigma \sqcup\{\1\}$ obtained by adding an element $\1$ and extending the partial order to $\eLL!_\Sigma$ by declaring $\sigma \subface \1$ for all $\sigma \in \LL!_\Sigma$.

The \emph{join} and \emph{meet} operations $\vee$ and $\wedge$ on $\LL!_\Sigma$ are defined as follows. For two cones $\sigma$ and $\tau$ of $\Sigma$, we set $\sigma \wedge \tau \coloneqq \sigma \cap \tau$. To define the operation $\vee$, note that the set of cones in $\Sigma$ which contain both $\sigma$ and $\tau$ is either empty or has a minimal element $\eta \in \Sigma$. In the former case, we set $\sigma \vee \tau \coloneqq \1$, and in the latter case, $\sigma \vee \tau \coloneqq \eta$. The two operations are extended to the augmented poset $\eLL!_\Sigma$ by $\sigma \wedge \1 = \sigma$ and $\sigma \vee \1 = \1$ for any cone $\sigma$ of $\Sigma$.

\subsection*{Notations}

The above discussion leads to the following notations. Let $\tau$ and $\sigma$ be a pair of faces in $\Sigma$. We say $\sigma$ \emph{covers} $\tau$ and write $\tau \ssubface \sigma$ if $\tau\subface\sigma$ and $\dims{\tau} = \dims\sigma-1$. A family of faces $\sigma_1,\dots,\sigma_k$ are called \emph{comparable} if $\sigma_1 \vee \cdots \vee \sigma_k \neq \1$. Moreover, we use the notation $\sigma\sim\sigma'$ for two faces $\sigma$ and $\sigma'$ if $\sigma\wedge\sigma' = \0$ and $\sigma\vee\sigma' \neq \1$.

\subsection{Star fan}

The \emph{star fan} $\Sigma^\sigma$ refers to the fan in $N^\sigma_\R=\rquot {N_\R}{N_{\sigma,\R}}$ induced by the cones $\eta$ in $\Sigma$ which contain $\sigma$ as a face. This is consistent with the terminology used in~\cite{AHK} and differs from the one in~\cites{Kar04, BBFK02} where this is called \emph{transversal fan}.

\subsection*{Notations}

For a cone $\eta$ in $\Sigma$ which is comparable to $\sigma$, we denote by $\eta^\sigma$ the corresponding cone in $\Sigma^\sigma$. That is, $\eta^\sigma \coloneqq \rquot{(\eta+N_{\sigma,\R})}{N_{\sigma,\R}}$.

In the opposite direction, if $\xi$ is a face of $\Sigma^\sigma$, then the set of faces of $\Sigma$ corresponding to $\xi$ is by definition the set of all cones $\eta$ in $\Sigma$ which are comparable to $\sigma$, and verify $\eta^\sigma =\xi$. They form a sublattice of the lattice of faces of $\Sigma$. We use the notation $\hat\xi$ for the smallest, and $\check{\xi}$ for the largest of these elements. (The notation is justified by the operations of meet $\wedge$ and join $\vee$ in lattice theory.) Note that $\dims{\hat\xi}=\dims\xi$ and that $\check\xi = \hat\xi\vee\sigma$.

In order to simplify the notation, we often drop the above notations for rays, and denote in the same way a ray in the fan and in its star fans. That is, for a ray $\zeta$ in $\Sigma^\sigma$, $\hat\zeta$ is denoted $\zeta$.

\subsection{Unit normal vectors} \label{subsec:orientation}

Let now $\Sigma$ be a rational fan of pure dimension $d$. Let $\sigma$ be a cone of $\Sigma$ and let $\tau$ be a face of codimension one in $\sigma$. Then, $N_{\tau,\R}$ cuts $N_{\sigma,\R}$ into two closed half-spaces only one of which contains $\sigma$. Denote this half-space by $H_\sigma$. By a \emph{unit normal vector to $\tau$ in $\sigma$} we mean any vector $v$ of $N_\sigma \cap H_\sigma$ such that $N_\tau + \Z v = N_\sigma$. We denote such an element by $\nvect_{\sigma/\tau}$. Two different choices of $\nvect_{\sigma/\tau}$ have the same projection in the quotient lattice $N^\tau_\sigma = \rquot{N_\sigma}{N_\tau}$. We thus get a well-defined vector in $N^\tau_\sigma$ that we denote by $\e^\tau_\sigma$. In the case $\sigma=\rho$ is a ray and $\tau =\conezero$, we simply use the notation $\e_\rho$ instead of $\nvect_{\rho/\conezero} = \e^{\conezero}_\rho$.

\subsection{Stellar subdivision}

\begin{figure}
  \begin{tikzpicture}[scale=2.5, x={(1.2cm,0cm)}, line join=round, baseline=1]
    \coordinate (O) at (0,0);
    \coordinate (I) at (1,0,-.5);
    \coordinate (J) at (1,0,.5);
    \coordinate (X) at (1,0,0);
    \coordinate (A) at (1,1,0);
    \coordinate (B) at (1,-1,0);
    \draw[thick, fill=gray!70] (J) -- (O) -- (I);
    \draw[thick, fill=red!21] (J) -- (O) -- (B);
    \fill[thick, fill=red!21] (J) -- (I) -- (B);
    \draw[-latex, dashed] (O) -- ($1.2*(X)$) node[right] {$\rho$};
    \fill[thick, opacity=.7, fill=blue!30] (J) -- (I) -- (A);
    \draw[thick, opacity=.7, fill=blue!30] (J) -- (O) -- (A);
    \draw[blue!70] (I) -- (A) -- (J);
    \draw[red!70] (I) -- (B) -- (J);
    \draw[gray] (I) -- (J);
  \end{tikzpicture}
  $\longrightarrow$
  \begin{tikzpicture}[scale=2.5, x={(1.2cm,0cm)}, line join=round, baseline=1]
    \coordinate (O) at (0,0);
    \coordinate (I) at (1,0,-.5);
    \coordinate (J) at (1,0,.5);
    \coordinate (X) at (1,0,0);
    \coordinate (A) at (1,1,0);
    \coordinate (B) at (1,-1,0);
    \draw[thick, fill=gray!70] (J) -- (O) -- (I);
    \draw[thick, opacity=.7, fill=gray!70] (A) -- (O) -- (B);
    \draw[thick, fill=red!21] (J) -- (O) -- (B);
    \fill[thick, fill=red!21] (J) -- (X) -- (B);
    \fill[thick, fill=green!21] (I) -- (X) -- (B);
    \fill[thick, opacity=.7, fill=brown!30] (I) -- (X) -- (A);
    \draw[green!70] (I) -- (B);
    \draw[thick, -latex] (O) -- ($1.2*(X)$) node[right] {$\rho$};
    \fill[thick, opacity=.7, fill=blue!30] (J) -- (X) -- (A);
    \draw[thick, opacity=.7, fill=blue!30] (J) -- (O) -- (A);
    \draw[blue!70] (J) -- (A);
    \draw[brown!70] (I) -- (A);
    \draw[red!70] (J) -- (B);
    \draw[gray] (I) -- (J);
    \draw[gray] (A) -- (B);
  \end{tikzpicture}

  \caption{A blow-up along the ray $\rho$. \label{fig:blowup}}
\end{figure}
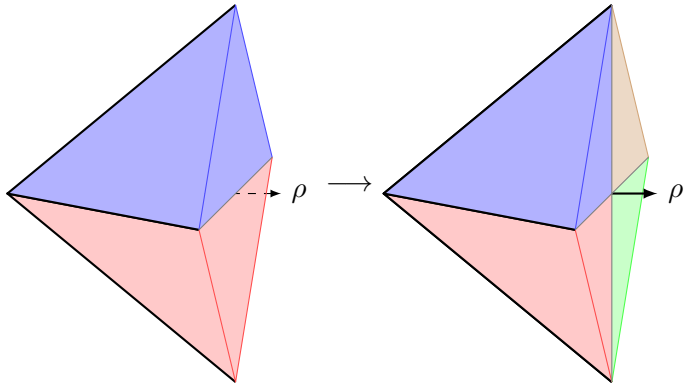

Let $\Sigma$ be a rational fan. Let $\sigma \in \Sigma$ be a cone of $\Sigma$. Let $\rho$ be a rational ray generated by a vector in the relative interior of $\sigma$. The \emph{blow-up of\/ $\Sigma$ along $\rho$}, also called the \emph{stellar subdivision of\/ $\Sigma$ along $\rho$}, is the rational fan $\Sigma_{(\rho)}$ defined as follows (see Figure \ref{fig:blowup}). For any cone $\eta \supfaceeq \sigma$, we remove $\eta$ from $\Sigma$, and replace it by the cones of the form $\tau+\rho$ with $\tau$ any proper face of $\eta$ which does not contain $\sigma$ and such that $\tau+\rho$ intersects the interior of $\eta$. We obtain a new fan with the same support which we denote by $\Sigma_{(\rho)}$. If $\Sigma$ is the blow-up of some fan $\Sigma'$ along a ray $\rho$, then $\Sigma'$ is called the \emph{stellar assembly of\/ $\Sigma$} or \emph{blow-down of\/ $\Sigma$ along $\rho$}. By an abuse of the terminology, we say a fan $\~\Sigma$ is obtained \emph{by a blow-up of\/ $\Sigma$ along $\sigma$} if there exists a ray $\rho$ in the relative interior of $\sigma$ as above so that $\~\Sigma$ coincides with $\Sigma_{(\rho)}$.

\smallskip
If $\Sigma$ is unimodular, a blow-up of $\Sigma$ along $\sigma$ is called \emph{unimodular} if $\Sigma_{(\rho)}$ is still unimodular. Such a blow-up is in fact unique. Indeed, for any face $\sigma$ with rays $\zeta_1, \dots, \zeta_k$, there is exactly one ray $\rho$ such that the blow-up along $\rho$ is unimodular. This ray is generated by $\e_{\zeta_1} + \dots + \e_{\zeta_k}$. Via the link to toric geometry, the unimodular blow-up $\~\Sigma$ of $\Sigma$ along $\sigma$ corresponds to the blow-up of the toric variety $\P_\Sigma$ along the closure $D^\sigma$ of the torus orbit $T^\sigma$ associated to $\sigma$. For this reason, we denote this blow-up by $\Bl{\Sigma}{\sigma}$.

\subsection{Local and fan irrelevant properties}

A property $P$ of rational fans is called \emph{local} if for any fan $\Sigma$ verifying $P$, all the star fans $\Sigma^\sigma$ for $\sigma \in \Sigma$ also verify $P$.

Let $N$ be a lattice and $N_\R$ be the corresponding real vector space. Consider another lattice $N'$ in $N_\R$. Let $\Sigma$ be a rational fan in $N_\R$ and $\Sigma'$ a rational fan in $N'_\R$. We say that $\Sigma$ and $\Sigma'$ \emph{have the same support} if we have $\supp\Sigma = \supp{\Sigma'}$ and $\supp\Sigma \cap N = \supp{\Sigma'}\cap N'$.

We say the property $P$ is \emph{fan irrelevant}, or a \emph{property of the support}, or we say \emph{$P$ only depends on the support}, if the following property holds: For any two rational fans with the same support $\Sigma$ and $\Sigma'$ in $N_\R$ and $N'_\R$, respectively, with $N_\R=N'_\R$, we have $\Sigma$ verifies $P$ if and only if $\Sigma'$ verifies $P$.

\subsection{Function theory} \label{sec:function_theory}

Let $\Sigma$ be a fan in $N_\R$ and let $\supp{\Sigma}$ be its support. The set of \emph{linear functions on $\Sigma$} is defined as the restriction to $\supp\Sigma$ of linear functions on $N_\R$; such a linear function is defined by an element of $M_\R$. In the case $\Sigma$ is rational, a linear function on $\Sigma$ is called \emph{integral}, \resp \emph{rational}, if it is defined by an element of $M$, \resp $M_\Q$.

Let $f\colon \supp\Sigma \to \R$ be a continuous function. We say that $f$ is \emph{conewise linear on $\Sigma$} if on each face $\sigma$ of $\Sigma$, the restriction $f\rest\sigma$ of $f$ to $\sigma$ is linear. In such a case, we simply write $f\colon \Sigma \to \R$, and denote by $f_\sigma$ the linear form on $N_{\sigma,\R}$ which coincides with $f\rest\sigma$ on $\sigma$. If the linear forms $f_\sigma$ are all integral, then we say $f$ is \emph{meromorphic on $\Sigma$}. We denote by $\mer(\Sigma)$ the set of all meromorphic functions on $\Sigma$. Pointwise addition of functions turns $\mer(\Sigma)$ into a group of finite rank.

\smallskip
A conewise linear function $f \colon \Sigma \to \R$ is called \emph{convex}, \resp \emph{strictly convex}, if for each face $\sigma$ of $\Sigma$, there exists a linear function $\ell$ on $\Sigma$ such that $f-\ell$ vanishes on $\sigma$ and is non-negative, \resp strictly positive, on $\eta \setminus \sigma$ for any cone $\eta \supfaceeq \sigma$ in $\Sigma$.

A fan $\Sigma$ is called \emph{quasi-projective} if it admits a strictly convex conewise linear function. A \emph{projective} fan is a fan which is both quasi-projective and complete. In the case $\Sigma$ is rational, it is quasi-projective, \resp projective, if and only if the toric variety $\P_\Sigma$ is quasi-projective, \resp projective.

\subsection{Saturation}\label{sec:saturation}

When working with integral coefficients, we need to require a saturation property for rational fans to get rid of torsion in Chow groups. A rational fan $\Sigma$ is called \emph{saturated at $\sigma$}, for a face $\sigma\in\Sigma$, if the set of integral linear functions on $\Sigma^\sigma$ (\ie, those induced by $M^\sigma$) coincides with the set of rational linear functions (\ie, elements of $M^\sigma_\Q$) whose restriction to any cone of $\Sigma^\sigma$ is integral. This is equivalent to requiring the lattice generated by the integral points $\supp{\Sigma^\sigma}\cap N^\sigma$ be saturated in $N^\sigma$.

We call a fan $\Sigma$ \emph{saturated} if it is saturated at all its faces.

For a discussion of the saturation with examples, we refer to Section~\ref{subsec:examples_non_saturated_non_unimodular}.

\subsection{Set-theoretical convention}

In this paper we work with fans modulo isomorphisms. For a rational fan $\Sigma$, we denote by $\ssN_\Sigma$ the restriction of the ambient lattice to the vector subspace of $N_\R$ spanned by $\Sigma$. Two rational fans $\Sigma$ and $\Sigma'$ are called isomorphic if there exists an integral linear isomorphism $\phi\colon \ssN_\Sigma\otimes\R \simto \ssN_{\Sigma'}\otimes\R$ inducing an isomorphism between $\ssN_\Sigma$ and $\ssN_{\Sigma'}$ such that for each face $\sigma \in \Sigma$, $\phi(\sigma)$ is a face of $\Sigma'$ and for each face $\sigma' \in \Sigma'$, $\phi^{-1}(\sigma')$ is a face of $\Sigma$.

Any rational fan is isomorphic to a fan in the space $\R^n$ endowed with the lattice $\Z^n$ for a sufficiently big $n$. Hence we can talk about the set of isomorphism classes of rational fans. In practice, by an abuse of the language, we will make no difference between an isomorphism class of rational fans and one of its representative.

\subsection{Minkowski weights on rational fans} \label{sec:MW}

Let $\Sigma$ be a rational fan in $N_\R$. Let $\corps = \Z$ or $\Q$. A \emph{Minkowski weight} of dimension $p$ on $\Sigma$ with coefficient in $\corps$ is a map $w \colon \Sigma_{p} \to \corps$ which verifies the following \emph{balancing condition}:

\[ \forall\:\tau\in \Sigma_{p-1},\qquad \sum_{\sigma\ssupface\tau} w(\sigma) \e^\tau_\sigma = 0 \in N^{\tau}. \]

We denote by $\MW_{p}(\Sigma)$ the set of all Minkowski weights of dimension $p$ on $\Sigma$ with integral coefficients. Addition of weights cell by cell turns $\MW_{p}(\Sigma)$ into a group. The set of Minkowski weights with rational coefficients is the vector space generated by $\MW_{p}(\Sigma)$, and is denoted by $\MW_{p}(\Sigma, \Q)$.

The \emph{support} of a Minkowski weight $w$ of dimension $p$ is the set of all cones $\sigma \in \Sigma_p$ with $w(\sigma) \neq 0$. A Minkowski weight is called \emph{effective} if $w$ takes only non-negative values. It is called \emph{reduced} if all the non-zero weights are equal to one, and called \emph{unitary} if the non-zero weights are equal to $+1$ or $-1$.

\subsection{Orientability}

We call a rational fan $\Sigma$ of pure dimension $d$ \emph{orientable} if there exists an element $\omega \in \MW_d(\Sigma)$ that has full support. We call any such $\omega$ an \emph{orientation} of $\Sigma$.

\begin{prop} \label{prop:orientable_stable}
Orientability is both local and a property of the support. Moreover, the product of two orientable fans is orientable.
\end{prop}

\begin{proof}
Let $\Sigma$ be a fan, and let $\omega \in \MW_d(\Sigma)$ be an orientation. For any cone $\sigma \in \Sigma$, the star fan $\Sigma^\sigma$ gets the induced orientation $\omega^\sigma$ defined by $\omega^\sigma(\eta^\sigma) = \omega(\eta)$, for any facet $\eta \supfaceeq \sigma$.

If $\Sigma'$ is another fan with the same support as $\Sigma$, then we get an induced orientation $\omega'$ on $\Sigma'$ by setting for any facet $\eta'\in\Sigma'$, $\omega'(\eta') \coloneqq \omega(\eta)$ where $\eta$ is a facet of $\Sigma$ which intersects the interior of $\eta'$. One checks that $\omega'(\eta')$ does not depend on the choice of $\eta$, and $\omega'$ is an orientation.

Finally, if $\Sigma'$ is any fan with an orientation $\omega'$, the fan $\Sigma\times\Sigma'$ can be endowed with the orientation $\omega\times\omega'$ defined by $\omega\times\omega'(\eta\times\eta') = \omega(\eta)\cdot\omega'(\eta')$ for any facets $\eta\in\Sigma$ and $\eta'\in\Sigma'$. Obviously, $\omega \times \omega'$ has full support.

A codimension one cone in $\Sigma \times \Sigma'$ is of the form $\sigma \times \tau'$ or $\tau \times \sigma'$ for facets $\sigma$ and $\sigma'$ of $\Sigma$ and $\Sigma'$, and codimension one faces $\tau$ and $\tau'$ of $\Sigma$ and $\Sigma'$, respectively. In either cases, the balancing condition in $\Sigma\times \Sigma'$ is a consequence of the balancing condition in $\Sigma$ or $\Sigma'$.
\end{proof}

\subsection{Tropical fans}

A \emph{tropical fan} is a pair $(\Sigma, \omega!_\Sigma)$ consisting of an orientable rational fan $\Sigma$ of pure dimension $d$ endowed with the choice of an orientation $\omega!_\Sigma$ for $\Sigma$. We call $\omega!_\Sigma$ the \emph{underlying orientation} of $\Sigma$. For any $\tau$ of dimension $d-1$, we have the balancing condition
\begin{align} \label{eq:codimension_one_balancing}
\sum_{\sigma \ssupface \tau} \omega!_\Sigma(\sigma)\e^\tau_{\sigma} =0 \quad \textrm{in $N^\tau$}.
\end{align}
In the following, sometimes we omit the mention of $\omega!_\Sigma$ and simply write $\Sigma$ when referring to a tropical fan $(\Sigma, \omega!_\Sigma)$. Also, if there is no risk of confusion, we simplify $\omega!_\Sigma$ to $\omega$.

A tropical fan $\Sigma$ is called \emph{effective} if its \emph{underlying orientation} $\omega!_\Sigma$ takes only positive values. It is called \emph{unitary} if $\omega!_\Sigma$ takes values 1 or $-1$ on any facet, and \emph{reduced} if it is both effective and unitary, \ie, $\omega!_\Sigma$ is constant equal to 1. As a direct consequence of the proof of Proposition~\ref{prop:orientable_stable}, we get the following.

\begin{prop}
Being effective, \resp unitary, \resp reduced, are local properties, properties of the support, and closed by products.
\end{prop}

Tropical fans arise naturally in connection with tropicalizations of subvarieties of algebraic tori. A tropical fan $(\Sigma, \omega!_\Sigma)$ which arises as the tropicalization of a subvariety $\mathbf{X}$ of an algebraic torus $\mathbf{T} = \spec(\k[M])$ over a trivially valued field $\k$ is called \emph{realizable over $\k$}. A tropical fan is called \emph{realizable} if it can be realized over some field $\k$. Results related to algebraic and complex geometry can be found in~\cites{AR10, GKM09, KM09, Katz12, Bab14, BH17, Gro18} and~\cites{MS15, MR18, BIMS}.

\subsection{Generalized tropical line}

Let $(\Sigma, \omega!_\Sigma)$ be a tropical fan of dimension one in $N_\R$. We say $\Sigma$ is a \emph{generalized tropical line} if any proper subset of the set of vectors $\e_\rho$, $\rho \in \Sigma_1$, form an independent set of vectors in $N_\R$ and moreover the integers $\omega!_\Sigma(\rho)$, $\rho\in \Sigma_1$, are coprime. This is equivalent to requiring $\MW_1(\Sigma) = \Z\,\omega!_\Sigma$.

\subsection{Connectedness through codimension one}

For any fan $\Sigma$ of pure dimension $d$, we define the (\emph{top dimensional}) \emph{dual graph} of $\Sigma$ as follows. This is the graph $G=(\Sigma_d, E)$ whose vertex set is equal to the set of facets $\Sigma_{d}$ and has an edge connecting any pair of facets $\sigma$ and $\eta$ which share a codimension one face in $\Sigma$. We say $\Sigma$ is \emph{connected through codimension one} if its dual graph is connected.

\medskip

\begin{center}
*\, *\, *
\end{center}

\medskip

We now discuss two basic algebrogeometric terminologies for tropical fans.

\subsection{Normal tropical fans}

Let $\Sigma$ be a tropical fan with underlying orientation $\omega!_\Sigma$. We call $\Sigma$ \emph{normal} if for any cone $\tau \in \Sigma_{d-1}$, we have $\MW_1(\Sigma^\tau) = \Z \, \omega!_{\Sigma^\tau}$. That is, we require the one-dimensional star fan $\Sigma^{\tau}$ endowed with the weights $\omega!_{\Sigma^\tau}$ be a generalized tropical line. This translates into the following property: for any cone $\tau \in \Sigma_{d-1}$ and for any collection of integer numbers $a_\sigma$ for $\sigma \ssupface \tau$, the relation
\[ \sum_{\sigma \ssupface \tau} a_\sigma\e^\tau_\sigma = 0 \in N^\tau \]
implies that the coefficients $a_\sigma$ are all a multiple by a common integer $\lambda$ of the weights $\omega!_\Sigma(\sigma)$.

We say that $\Sigma$ is \emph{$\Q$-normal} if $\MW_1(\Sigma^\tau, \Q) = \Q \, \omega!_{\Sigma^\tau}$ for any $\tau\in \Sigma_{d-1}$.

We have the following proposition.

\begin{prop} \label{prop:trop_normal_support_local_product}
Being normal, \resp $\Q$-normal, for tropical fans is a local property. The product of two normal, \resp $\Q$-normal, tropical fans is normal, \resp $\Q$-normal. Being $\Q$-normal is a property of the support. Moreover, $\Q$-normality and normality are equivalent for unitary tropical fans.
\end{prop}

\begin{proof}
The first claim is obvious. The second follows from the definition of the product orientation and the observation that the star fan of a cone of codimension one in the product $\Sigma\times\Sigma'$ coincides either with the star fan of a cone in $\Sigma$ or with the star fan of a cone in $\Sigma'$. The two last claims are straightforward.
\end{proof}

\subsection{Irreducible tropical fans} \label{sec:irreducible}

We now present a tropical notion of irreducibility. We say that a tropical fan $\Sigma$ is \emph{irreducible at a face $\sigma \in \Sigma$} provided that we have $\MW_{d-\dims\sigma}(\Sigma^\sigma) \simeq \Z \omega!_{\Sigma^\sigma}$. A tropical fan irreducible at $\conezero$ is simply called \emph{irreducible}. We call $\Sigma$ \emph{locally irreducible} if it is irreducible at any face $\sigma\in\Sigma$. We define \emph{$\Q$-irreducibility} and \emph{$\Q$-local irreducibility} as above working with rational coefficients, that is, $\Sigma$ is $\Q$-irreducible at a face $\sigma$ provided that $\MW_{d-\dims\sigma}(\Sigma^\sigma, \Q) \simeq \Q \omega!_{\Sigma^\sigma}$, and $\Sigma$ is $\Q$-locally irreducible if it is $\Q$-irreducible at any face.

\begin{prop}
Being locally irreducible and $\Q$-locally irreducible are local properties. $\Q$-local irreducibility is a property of the support. Moreover, for unitary tropical fans, $\Q$-local irreducibility and local irreducibility are equivalent.
\end{prop}

\begin{proof}
The first claim is tautological from the definition. The second and third can be obtained by a direct verification.
\end{proof}

The following theorem gives a link between normality and local irreducibility.

\begin{thm}[Characterization of locally and $\Q$-locally irreducible tropical fans] \label{thm:characterization_irreducible}
The following assertions are equivalent for a tropical fan $\Sigma$.
\begin{enumerate}
\item $\Sigma$ is locally irreducible.
\item $\Sigma$ is normal and each star fan $\Sigma^\sigma$, $\sigma \in \Sigma$, is connected through codimension one.
\end{enumerate}
Similarly, $\Sigma$ is $\Q$-locally irreducible if and only if it is $\Q$-normal and each star fan $\Sigma^\sigma$, $\sigma \in \Sigma$, is connected through codimension one.
\end{thm}

\begin{proof}
Denote by $d$ the dimension of $\Sigma$.

\smallskip
(1) $\Rightarrow$ (2). Suppose $\Sigma$ is locally irreducible. For a cone $\tau$ of codimension one in $\Sigma$, we have $\MW_{1}(\Sigma^\tau) = \Z \omega^\tau$. This implies that $\Sigma$ is normal. It remains to prove that for any $\sigma\in \Sigma$, the star fan $\Sigma^\sigma$ is connected through codimension one. Let $G^\sigma$ be the dual graph of $\Sigma^\sigma$, and suppose for the sake of a contradiction that $G^\sigma$ is not connected. This means we can find a partition of the facets of $\Sigma^\sigma$ into a disjoint union $S_1 \sqcup \dots \sqcup S_l$, for an integer $l\geq 2$, so that $S_j$ form the connected components of $G^\sigma$. For each $j$, consider the restriction $\omega!_{\Sigma^\sigma}\rest{S_j}$ and extend it by zero to all the facets of $\Sigma^\sigma$ to obtain the weight function $\alpha_j\colon \Sigma^\sigma_{d-\dims\sigma} \to \Z$. The elements $\alpha_j$ all belong to $\MW_{d-\dims\sigma}(\Sigma^\sigma)$, and they are not scalar multiples of each other. This contradicts the irreducibility of $\Sigma^\sigma$.

\smallskip
(2) $\Rightarrow$ (1). We now prove the reverse implication. Suppose $\Sigma$ is normal and moreover, for each $\sigma \in \Sigma$, $\Sigma^\sigma$ is connected through codimension one. Let $\sigma$ be a face of $\Sigma$ and let $\alpha \in \MW_{d-\dims\sigma}(\Sigma^\sigma)$ be a Minkowski weight of top dimension. We need to show that $\alpha = \lambda \omega!_{\Sigma^\sigma}$ for an integer $\lambda$. A codimension one face $\tau$ of $\Sigma^\sigma$ is of the form $\eta^\sigma$ for a codimension one face $\eta\supfaceeq \sigma$ of $\Sigma$, and we have the equality of star fans $\Sigma^\eta = (\Sigma^\sigma)^\tau$. Applying the normality condition to the codimension one face $\tau=\eta^\sigma$ of $\Sigma^\sigma$, we infer the existence of an integer $\lambda_\tau$ such that $\alpha^\tau =\lambda_\tau \omega!_{\Sigma^\eta}$. Using now the connectivity of $\Sigma^\sigma$ through codimension one, we conclude that the scalars $\lambda_\tau$ are all equal. This shows $\alpha$ is a multiple of $\omega!_{\Sigma^\sigma}$, as required.

\smallskip
The same reasoning gives the $\Q$-statement.
\end{proof}

\begin{prop} \label{prop:irreducibility_product}
The class of locally irreducible, \resp $\Q$-locally irreducible, tropical fans is closed under products.
\end{prop}

\begin{proof}
This follows from Theorem~\ref{thm:characterization_irreducible}, Proposition~\ref{prop:trop_normal_support_local_product}, and the observation that product of two fans that are each connected through codimension one, is again connected through codimension one.
\end{proof}

Theorem \ref{thm:characterization_irreducible} justifies the following definition via Proposition \ref{prop:irreducible_components_irreducible}.

\begin{defi}[Irreducible components of a normal tropical fan]
Let $\Sigma$ be a ($\Q$-)normal, tropical fan of dimension $d$ with underlying orientation $\omega!_\Sigma$. Consider the dual graph $G$ of $\Sigma$ with the set of vertices $\Sigma_d$. Any connected component of $G$ with vertex set $V \subseteq \Sigma_{d}$ defines a subfan $\ssSigma_V$ of $\Sigma$ defined by
\[ \ssSigma_V \coloneqq \left\{ \eta \in \Sigma \st \eta \subfaceeq \sigma \textrm{ for some }\sigma \in V \right\}. \]
The subfan $\ssSigma_V$ comes with underlying orientation $\omega!_\Sigma\rest{V}$ which is tropical, ($\Q$-)normal and connected through codimension one. We refer to the subfans $\ssSigma_V$ as \emph{\textnormal{($\Q$-)}irreducible components} of $\Sigma$.
\end{defi}

We refer to Section~\ref{sec:alternate_irreducible} for an example, and an alternative definition of irreducible components. The following proposition justify the terminology.

\begin{prop} \label{prop:irreducible_components_irreducible}
Let $\Sigma$ be a \textnormal{$(\Q$-)}normal tropical fan of dimension $d$. The \textnormal{($\Q$-)}irreducible components of\/ $\Sigma$ are \textnormal{($\Q$-)}irreducible. Moreover, they induce a partition of the facets $\Sigma_d$.
\end{prop}

\begin{proof}
This follows from the proof of Theorem \ref{thm:characterization_irreducible}.
\end{proof}

Note that the statement in the proposition only claims irreducibility at $\conezero$. This is weaker than local irreducibility. We refer to Example \ref{ex:PD_chow_not_smooth} for an irreducible tropical fan at $\conezero$ that is not locally irreducible.

\section{Chow rings of fans} \label{sec:chow_ring}

In this section, we review the definition of the Chow rings associated to rational simplicial fans and discuss their basic properties. The new result here is a combinatorial proof of the Localization lemma, Theorem~\ref{thm:kernel_Z^k}, given in Section~\ref{sec:kernel_Z^k}.

We work with integer coefficients unless it is explicitly stated that the coefficients are rational numbers. We will require $\Sigma$ be unimodular for the cycle class map and for Poincaré duality with integer coefficients.

\subsection{Definition of the Chow ring}

Let $\Sigma$ be a rational simplicial fan of dimension $d$ in $N_\R$. Consider the polynomial ring $\Z[\x_\zeta]_{\zeta\in \Sigma_1}$ with indeterminate variables $\x_\zeta$ associated to rays $\zeta$ in $\Sigma_1$. The Chow ring $A^\bul(\Sigma)$ of $\Sigma$ with integer coefficients is by definition the quotient ring
\[ A^\bul(\Sigma) \coloneqq \rquot{\Z[\x_\zeta]_{\zeta\in \Sigma_1}}{\bigl(I + J\bigr)} \]
where
\begin{itemize}
\item $I$ is the ideal generated by the products $\x_{\rho_1}\!\cdots \x_{\rho_k}$, $k\in \N$, such that $\rho_1, \dots, \rho_k$ are not comparable in $\Sigma$, and
\item $J$ is the ideal generated by the elements of the form
\[ \sum_{\zeta\in \Sigma_1} m(\e_\zeta)\x_\zeta \]
for $m \in M = N^\dual$.
\end{itemize}

The ideal $I+J$ is homogeneous and the Chow ring inherits a graded ring structure. We can thus write
\[ A^\bul(\Sigma) = \bigoplus_{k\geq 0} A^k(\Sigma) \]
with the $k$-th degree piece $A^k(\Sigma)$, $k\in \Z_{\geq 0}$, generated by degree $k$ monomials as $\Z$-module.

For each ray $\zeta$ of $\Sigma$, we denote by $x_\zeta$ the image of $\x_\zeta$ in $A^1(\Sigma)$. More generally, for each cone $\sigma$ with rays $\rho_1, \dots, \rho_{\dims\sigma}$, we denote by $x_\sigma$ the product $x_{\rho_1}\dots x_{\rho_{\dims\sigma}}$.

A meromorphic function $f \in \mer(\Sigma)$ gives an element of $A^1(\Sigma)$ denoted $\ell(f)$ and defined as
\[ \ell(f)\coloneqq \sum_{\zeta \in \Sigma_1} f(\e_\zeta) x_\zeta. \]
If $\Sigma$ is unimodular, all the elements of $A^1(\Sigma)$ are of this form, and $A^1(\Sigma)$ can be identified with the quotient space $A^1(\Sigma) \simeq \rquot{\mer(\Sigma)}{M}$, the space of meromorphic functions on $\Sigma$ modulo integral linear functions.

When the fan $\Sigma$ is unimodular, we have the following characterization of the Chow ring, \cf\cites{Dan78, BDP90, Bri96, FS97}.

\begin{thm}
Let $\Sigma$ be a unimodular fan, and denote by $\P_\Sigma$ the corresponding toric variety. The Chow ring $A^\bullet(\Sigma)$ is isomorphic to the Chow ring of\/ $\P_\Sigma$.
\end{thm}

In order to distinguish the fan we are referring to, we sometimes denote by $\ssN_\Sigma$ and $\ssM_\Sigma$ the lattices underlying the definition of $\Sigma$, $\ssM_\Sigma = \ssN_\Sigma^\dual$, and denote by $\ssI_{\Sigma}$ and $\ssJ_{\Sigma}$ the ideals used in the definition of the Chow ring $A^\bul(\Sigma)$.

\smallskip
The Chow ring with rational, \resp real, coefficients, are denoted by $A^\bullet(\Sigma, \Q)$, \resp $A^\bullet(\Sigma, \R)$, and defined as
\[A^\bullet(\Sigma, \Q) = A^\bullet(\Sigma)\otimes_\Z \,\Q \qquad \textrm{and} \qquad A^\bullet(\Sigma, \R) = A^\bullet(\Sigma)\otimes_\Z \,\R.\]

\subsection{Localization lemma}

Consider a rational simplicial fan $\Sigma$ in $N_\R$. We do not assume in this subsection that $\Sigma$ is tropical.

Consider the Chow ring $A^\bul(\Sigma)$. For each cone $\sigma \in \Sigma$, let $\x_\sigma$ be an indeterminate variable and define $Z^k(\Sigma) \coloneqq \bigoplus_{\sigma\in\Sigma_k} \Z\x_\sigma.$ We have a natural embedding
\begin{align*}
Z^k(\Sigma)= \bigoplus_{\sigma\in\Sigma_k} \Z\x_\sigma &\hookrightarrow \Z[\x_\rho \mid \rho\in \Sigma_1] \\
\x_\sigma &\mapsto \prod_{\rho\subfaceeq \sigma\\\dims\rho=1}\x_\rho.
\end{align*}
This gives an additive map
\[ Z^k(\Sigma) \to A^k(\Sigma) \]
which sends $\x_\sigma$ to $x_\sigma$. We denote the kernel of this map by $Z^k_0$.

\begin{thm}[Localization lemma] \label{thm:kernel_Z^k}
The map $Z^k(\Sigma) \to A^k(\Sigma)$ is surjective and its kernel $Z^k_0$ is generated by elements of the form
\[ \sum_{\sigma \ssupface\tau} m(\nvect_{\sigma/\tau})\x_\sigma \]
for $\tau$ in $\Sigma_{k-1}$ and $m$ an element in $M$ that vanishes on $N_\tau$ (equivalently, $m\in M^\tau$).
\end{thm}
Note in particular that $A^k(\Sigma) = 0$ for $k > \dim(\Sigma)$.

\begin{proof}
This is a special case of a more general result stated in~\cite{FMSS}*{Theorem 1}, which establishes an isomorphism between the Chow groups of a variety endowed with an action of a solvable linear algebraic group on one side, and equivariant Chow groups associated to the variety, defined by cycles and relations which are invariant under the action of the group. A generalization of the result can be found in~\cite{Tot14}; see also~\cites{Fra06, Josh01, Pay06}.

We will provide a combinatorial proof of this result in Section \ref{sec:kernel_Z^k}.
\end{proof}

\subsection{Duality between Chow groups and Minkowski weights} \label{rem:chow_to_mw}

In the case the fan $\Sigma$ is unimodular, the localization lemma provides an isomorphism $\MW_k(\Sigma) \simeq A^k(\Sigma)^\dual$. Namely, for any non-negative integer $k$, consider the pairing
\[ \arraycolsep=1pt \begin{array}{rlll}
\langle \cdot\,, \cdot \rangle \colon Z^k(\Sigma) &\times& \MW_k(\Sigma) &\ \longrightarrow\  \Z \\
\x_\sigma&,& w &\ \longmapsto\  w(\sigma), \qquad \sigma\in\Sigma_k,\, w\in \MW_k(\Sigma),
\end{array} \]
extended linearly to all $Z^k$. We have the following duality theorem.

\begin{thm}[Duality Theorem] \label{thm:duality_A_MW}
Suppose $\Sigma$ is unimodular. The bilinear pairing above vanishes on the kernel $Z_0^k(\Sigma)$ of the map $Z^k(\Sigma) \to A^k(\Sigma)$. The induced map $\MW_k(\Sigma) \to A^k(\Sigma)^\dual$ is an isomorphism.
\end{thm}

\begin{proof}
See~\cite{AHK} for more details. The first assertion follows from the balancing condition for the Minkowski weight $w$ and the localization lemma. The second assertion is again a direct consequence of the localization lemma, and of the following fact (see also \cite{AHK}*{Proposition~5.6}). An integer valued function $w \colon \Sigma_k \to \Z$ is a Minkowski weight if and only if it is orthogonal to $Z_0^k$ via the pairing $\langle \cdot\,, \cdot \rangle$, that is, $\MW_k(\Sigma)$ can be identified with $Z_0^k(\Sigma)^\perp$.
\end{proof}

\begin{remark}
There is one subtle point here worth emphasizing. In general, the pairing between $A^k(\Sigma)$ and $\MW_k(\Sigma)$ is not perfect. Indeed, as we will show in Examples \ref{ex:F1_vs_N} and~\ref{ex:non_saturated_smooth_fan}, the Chow ring $A^k(\Sigma)$ can have torsion, in which case, $A^k(\Sigma)$ cannot be isomorphic to $\MW_k(\Sigma)^\dual$. We discuss torsion-freeness of Chow rings of tropical fans further in Section~\ref{sec:divisors}.
\end{remark}

\subsection{Künneth formula}

The Chow ring and Minkowski weights of a product of two simplicial fans is described as follows. The tensor products in this section are all over $\Z$.

\begin{prop}[Künneth formula]
Let $\Sigma$ and $\Sigma'$ be two rational simplicial fans. Then,
\[ A^\bul(\Sigma \times \Sigma') \simeq A^\bul(\Sigma) \otimes A^\bul(\Sigma') \quad\text{and}\quad \MW_\bul(\Sigma \times \Sigma') \simeq \MW_\bul(\Sigma) \otimes \MW_\bul(\Sigma'). \]
The first isomorphism is moreover a ring isomorphism.
\end{prop}

Künneth formula for Chow rings of toric varieties is known, see~\cite{FMSS}*{Theorem 2}.

\begin{proof}
We have $\ssM_{\Sigma\times\Sigma'}\simeq \ssM_{\Sigma} \times \ssM_{\Sigma'}$. It is easy to check that
\begin{gather*}
\ssI_{\Sigma}\otimes \Z[\x_\zeta \mid \zeta \in \Sigma'_1]\,+\,\Z[\x_\zeta \mid \zeta \in \Sigma_1]\otimes \ssI_{\Sigma'} = \ssI_{\Sigma\times\Sigma'},\quad\text{and} \\
\ssJ_{\Sigma}\otimes \Z[\x_\zeta \mid \zeta \in \Sigma'_1]\,+\,\Z[\x_\zeta \mid \zeta \in \Sigma_1]\otimes \ssJ_{\Sigma'} = \ssJ_{\Sigma\times\Sigma'}.
\end{gather*}
The statement about the Chow rings then follows because the tensor product is right-exact.

\smallskip
For Minkowski weights, there is a perfect pairing $Z^k(\Sigma) \times \W_k(\Sigma) \to \Z$ where $\W_k(\Sigma)\simeq \Z^{\Sigma_k}$ denotes the additive group of all the integer valued maps $w\colon \Sigma_k \to \Z$. Taking the orthogonal sum over different degrees leads to a perfect pairing between $Z^\bul(\Sigma) \coloneqq \bigoplus_{k=0}^d Z^k(\Sigma)$ and $\W_\bul(\Sigma)\coloneqq\bigoplus_{k=0}^d\W_k(\Sigma)$.

Let as before $Z_0^k(\Sigma)$ be the kernel of the surjective map $Z^k(\Sigma) \to A^k(\Sigma)$. By Localization lemma, we have $\MW_k(\Sigma) \simeq Z_0^k(\Sigma)^\perp$ (see the proof of Theorem~\ref{thm:duality_A_MW}). By the description of $Z^k_0$ given by Localization lemma, we deduce that
\[ Z_0^\bul(\Sigma\times\Sigma') = Z_0^\bul(\Sigma) \otimes Z^\bul(\Sigma') \,+\, Z^\bul(\Sigma) \otimes Z_0^\bul(\Sigma'). \]
Taking the orthogonal of each side in $\W_\bul(\Sigma \times \Sigma') \simeq \W_\bul(\Sigma)\otimes \W_\bul(\Sigma')$, we get the desired isomorphism $\MW_\bul(\Sigma\times\Sigma') \simeq \MW_\bul(\Sigma)\otimes\MW_\bul(\Sigma')$.
\end{proof}

\subsection{Degree maps}

Let $\Sigma$ be unimodular. Using the isomorphism $\MW_d(\Sigma) \simeq A^d(\Sigma)^\dual$, we associate to any element $\omega \in \MW_d(\Sigma)$, the corresponding \emph{degree map} $\deg_\omega$ defined by
\begin{align*}
\deg_\omega\colon A^d(\Sigma)&\to\Z \\
x_\sigma &\mapsto \omega(\sigma).
\end{align*}

If $\Sigma$ is a unimodular tropical fan with underlying orientation $\omega!_\Sigma$, we obtain a well defined degree map $\deg!_\Sigma \colon A^d(\Sigma)\to\Z$ by setting $\deg!_\Sigma \coloneqq \deg_{\omega!_\Sigma}$, that we abbreviate to $\deg$ if there is no risk of confusion. Moreover, for each cone $\sigma\in\Sigma$, we denote by $\deg_\sigma\colon A^{d-\dims\sigma}(\Sigma^\sigma)\to\Z$ the corresponding degree map, that is, $\deg_\sigma = \deg!_{\Sigma^\sigma}$.

\subsection{Cycle class map} \label{subsec:cycle_class_map}

Let $\Sigma$ be a unimodular tropical fan with underlying orientation $\omega!_\Sigma$ and the corresponding degree map $\deg!_{\Sigma}$.

The composition of the product map in the Chow ring with the degree map gives a bilinear pairing
\[ \begin{tikzcd}[column sep=small, row sep=0pt]
A^k(\Sigma) \times A^{d-k}(\Sigma) \rar& A^d(\Sigma) \rar& \Z,\qquad\\
(\alpha, \beta) \ar[rr, mapsto]&& \deg!_\Sigma(\alpha \cdot \beta),
\end{tikzcd} \]
for each non-negative integer $k\leq d$. This leads to a map $A^k(\Sigma) \to A^{d-k}(\Sigma)^\dual$ which sends an element $\alpha \in A^k(\Sigma)$ to the element $\bigl(\beta \mapsto \deg!_\Sigma(\alpha\cdot\beta)\bigr)$ of $A^{d-k}(\Sigma)^\dual$.

\begin{defi}[Cycle class map] \label{defi:map-chow-homology}
Let $\Sigma$ be a tropical fan. For any integer $k$, the \emph{cycle class map} $\cl$ is the map
\[ \cl \colon A^k(\Sigma) \to \MW_{d-k}(\Sigma) \]
given by the composition of the following maps
\[ A^k(\Sigma) \to A^{d-k}(\Sigma)^\dual \simto \MW_{d-k}(\Sigma). \qedhere \]
\end{defi}

Concretely, for each element $\alpha \in A^k(\Sigma)$, $\cl(\alpha)$ is the Minkowski weight in $\MW_{d-k}(\Sigma)$ defined by
\begin{align*}
\cl(\alpha) \colon \Sigma_{d-k} &\to \Z\\
\sigma &\mapsto \deg!_\Sigma(\alpha\cdot x_\sigma).
\end{align*}

\subsection{Restriction and Gysin maps} \label{subsec:restr_gys}

Chow rings of a rational simplicial fan $\Sigma$ and their star fans are related by two types of maps called \emph{restriction} and \emph{Gysin}, that we discuss now.

For a pair of cones $\tau \subfaceeq \sigma$ in $\Sigma$, we define the restriction and Gysin maps $\i^*_{\tau \subfaceeq \sigma}$ and $\gys_{\sigma \supfaceeq \tau}$ between the Chow rings of $\Sigma^\tau$ and $\Sigma^\sigma$.

The restriction map
\[ \i^*_{\tau \subfaceeq \sigma}\colon A^\bul(\Sigma^\tau) \to A^\bul(\Sigma^\sigma) \]
is a graded $\Z$-algebra homomorphism
\[ \i^*_{\tau \subfaceeq \sigma}\colon A^1(\Sigma^\tau) \to A^1(\Sigma^\sigma) \]
defined on generating sets by
\[ \forall \, \rho \in \Sigma_1^\tau, \qquad \i^*_{\tau \subfaceeq \sigma} (x_\rho) =
\begin{cases}
  x_\rho & \qquad  \textrm{if $\sigma \sim \rho$}, \\
  -\sum_{\zeta \in \Sigma^\sigma_1} m(\e_{\hat\zeta}) x_\zeta & \qquad \textrm{if $\rho \in \sigma$}, \\
  0 & \qquad \textrm{otherwise},
\end{cases} \]
where
\begin{itemize}
\item in the first equality, we view the ray $\rho$ of $\Sigma^\tau$ as a ray $\hat \rho$ in $\Sigma$, and identify it with a ray in $\Sigma^\sigma$.
\item in the second equality, $m$ is any element in $M^\tau = (N^\tau)^\dual$ that takes value $1$ on $\e_\rho$ and value zero on other rays of $\sigma^\tau$. The ray $\hat\zeta$ is the one in $\Sigma^\tau$ associated to $\zeta \in \Sigma_1^\sigma$.
\end{itemize}
Note that any two such choices of $m$ and $m'$ differ by an element which vanishes on $N_\sigma$, that is, by an element of $M^\sigma$. This means the element $\sum_{\zeta \in \Sigma^\sigma} m(\e_\zeta) x_\zeta$ in $A^1(\Sigma^\sigma)$ does not depend on the choice of $m$.

\smallskip
The Gysin map is the $\Z$-module morphism
\[ \gys_{\sigma \supfaceeq \tau} \colon A^{\bul}(\Sigma^\sigma) \longrightarrow A^{\bul+ \dims \sigma -\dims \tau}(\Sigma^\tau) \]
defined as follows. Let $r = \dims \sigma -\dims \tau$, and denote by $\rho_1, \dots, \rho_r$ the rays of $\sigma$ which are not in $\tau$.

Consider the $\Z$-module map
\[ \Z[\x_\zeta]_{\substack{\zeta \in \Sigma_1^\sigma}} \longrightarrow \Z[\x_\zeta]_{\substack{\zeta \in \Sigma^\tau_1}} \]
defined by multiplication by $\x_{\rho_1}\x_{\rho_2} \dots \x_{\rho_r}$. Obviously, it sends an element of the ideal $\ssI_{\Sigma^\sigma}$ in the source to an element of the ideal $\ssI_{\Sigma^\tau}$ in the target. Moreover, the projection
\[ N^\tau = \rquot{N}{N_\tau} \longtwoheadrightarrow N^\sigma = \rquot{N}{N_{\sigma}} \]
gives an injection
\[ M^\sigma = (N^\sigma)^\dual \longhookrightarrow M^\tau = (N^\tau)^\dual. \]

This shows that the elements of $\ssJ_{\Sigma^\sigma}$ in the source are sent to elements of $\ssJ_{\Sigma^\tau}$ in the target as well. Passing to the quotient, we get a $\Z$-module map
\[ \gys_{\sigma \supfaceeq \tau} \colon A^{k}(\Sigma^\sigma) \longrightarrow A^{k+ \dims \sigma -\dims \tau}(\Sigma^\tau). \]

The following proposition gathers some basic properties of the restriction and Gysin maps.

\begin{prop} \label{lem:i_gys_basic_properties-local}
Let $\tau\subfaceeq\sigma$ be a pair of faces, and let $x\in A^\bul(\Sigma^\tau)$ and $y\in A^\bul(\Sigma^\sigma)$. Let $\sigma^\tau$ be the face associated to $\sigma$ in $\Sigma^\tau$, and let $x_{\sigma/\tau}$ be the associated element of $A^{\dims\sigma-\dims\tau}(\Sigma^\tau)$. Then, we have the following compatibility properties between the restriction and Gysin maps.
\begin{gather}
\text{The restriction map $\i^*_{\tau\subfaceeq\sigma}$ is a surjective ring homomorphism.} \label{eqn:i_surjective_homeo-local} \\
\text{We have }\gys_{\sigma\supfaceeq\tau}\circ\,\i^*_{\tau\subfaceeq \sigma}(x)=x_{\sigma/\tau}\cdot x. \label{eqn:gys_circ_i-local} \\
\text{We have }\gys_{\sigma\supfaceeq\tau}(\i^*_{\tau\subfaceeq\sigma}(x)\cdot y)=x\cdot\gys_{\sigma \supfaceeq\tau}(y)
\,\,\, \text{called \textnormal{(Projection Formula)}}. \label{eqn:gys_i_simplification-local}
\end{gather}
\end{prop}

\begin{proof}
In order to simplify the presentation, we drop the indices of $\gys$ and $\i^*$. Properties~\eqref{eqn:i_surjective_homeo-local} and~\eqref{eqn:gys_circ_i-local} follow directly from the definitions. From Equation~\eqref{eqn:gys_circ_i-local}, we can deduce Equation~\eqref{eqn:gys_i_simplification-local} by the following calculation. Let $\~y$ be a preimage of $y$ by $\i^*$. Then,
\[ \gys(\i^*(x)\cdot y)=\gys(\i^*(x\cdot \~y))=x_{\delta/\gamma}\cdot x\cdot\~y=x\cdot\gys\circ\,\i^*(\~y)=x\cdot\gys(y). \qedhere \]
\end{proof}

Assume now that $\Sigma$ is a tropical fan, and for each cone $\sigma \in \Sigma$, let $\deg_\sigma\colon A^{d-\dims\sigma}(\Sigma^\sigma)\to\Z$ be the corresponding degree map.

\begin{prop} \label{lem:i_gys_basic_properties-local-tropical}
Notations as in \textnormal{Proposition~\ref{lem:i_gys_basic_properties-local}}, assume $\Sigma$ is a tropical fan. Then, we have
\begin{equation} \label{eqn:deg_circ_gys-local}
\deg_\sigma=\deg_\tau\circ\gys_{\sigma\supfaceeq\tau}.
\end{equation}
Moreover, $\gys_{\sigma\supfaceeq\tau}$ and $\i^*_{\tau\subfaceeq\sigma}$ are dual in the sense that
\begin{equation} \label{eqn:i_gys_dual-local}
\deg_\tau(x\cdot\gys_{\sigma\supfaceeq\tau}(y))=\deg_\sigma(\i^*_{\tau\subfaceeq\sigma}(x)\cdot y).
\end{equation}
\end{prop}

\begin{proof}
We drop the indices of $\gys$ and $\i^*$. For Equation~\eqref{eqn:deg_circ_gys-local}, let $\eta$ be a facet of $\Sigma^\sigma$. Let $\check\eta$ be the corresponding cone containing $\sigma^\tau$ in $\Sigma^\tau$, that is maximal among all the cones in $\Sigma^\tau$ containing $\sigma^\tau$. We have the respective elements $x_{\check\eta}\in A^{d-\dims\tau}(\Sigma^\tau)$ and $x_\eta\in A^{d-\dims\sigma}(\Sigma^\sigma)$. By definition of the degree maps, $\deg_\tau(x_{\check\eta})=\deg_\sigma(x_\eta)$. Using the definition of $\gys$, we get that $x_{\check\eta}=\gys(x_\eta)$, from which we conclude $\deg_\sigma = \deg_\tau\circ\gys$. Finally, we obtain Equation~\eqref{eqn:i_gys_dual-local} by the chain of equalities
\[ \deg_\tau(x\cdot\gys(y))=\deg_\tau(\gys(\i^*(x)\cdot y))=\deg_\sigma(\i^*(x)\cdot y). \qedhere \]
\end{proof}

\subsection{Combinatorial proof of Localization lemma~\ref{thm:kernel_Z^k}} \label{sec:kernel_Z^k}

In this section, we give a combinatorial proof of Theorem~\ref{thm:kernel_Z^k}. The result is used by Adiprasito-Huh-Katz~\cite{AHK} to deduce the duality between Chow groups and Minkowski weights of unimodular fans. While the other proofs in \cite{AHK} are written to be accessible to a combinatorial audience, this one refers to~\cite{FMSS}, which studies equivariant Chow rings of algebraic varieties admitting a solvable group action. The proof we give is elementary and arguably more transparent. Moreover, combined with a similar in spirit result proved in~\cite{Ami20} for the combinatorial Chow rings of products of graphs, it suggests a more general theory of combinatorial Chow rings associated to semistable degenerations of algebraic varieties, that merits a further study.

\subsubsection{Admissible expansions} Let $\Sigma$ be any simplicial rational fan. It will be convenient in the course of the proof to introduce the following notation. For an element $\ell \in M =N^\dual$, let
\[ \x_\ell \coloneqq \sum_{\rho \in \Sigma_1} \ell(\e_\rho) \x_\rho. \]

Let $B^k$ be the group of homogeneous polynomials of degree $k$ in $\Z[\x_\rho \mid \rho\in\Sigma_1]$ and $Z^k = \bigoplus_{\sigma \in \Sigma_k} \Z\x_\sigma$, which can be identified with a subgroup of $B^k$ via the identification between $\x_\sigma$ and the monomial $\prod_{\rho\subfaceeq \sigma \\ \rho \in\Sigma_1}\x_{\rho}$. Recall that $I$ is the ideal in $\Z[\x_\rho \mid \rho\in\Sigma_1]$ generated by the products $\x_{\rho_1}\cdots \x_{\rho_s}$, $s\in \N$, such that $\rho_1, \dots, \rho_s$ are not comparable $\Sigma$, and that $J$ is the ideal generated by the elements $\x_\ell$ for $\ell \in M$.

We have to show that $Z^k \to \rquot{B^k}{B^k \cap(I + J)} $ is surjective and determine its kernel. We will prove the latter, the former becomes clear in the course of the proof; one can also find a proof of the surjectivity in \cite{AHK}.

\smallskip
Let $K^k$ be the subgroup of $B^k$ generated by elements of the form
\[ \sum_{\rho \in \Sigma_1 \\ \rho\sim\tau} \ell(\e_\rho) \x_\rho \x_\tau, \]
for $\tau \in \Sigma_{k-1}$ and $\ell$ an element in $M$ which is orthogonal to $\tau$, that is, $\ell \in M^\tau$. Consider an element $a$ of $B^k \cap \bigl(I+J\bigr)$ which belongs to $Z^k$. We shall prove that $a \in K^k$.

\smallskip
First, notice that we can write $a$ as an element of $I$ plus a sum of monomials of the form $\x_\ell \x_{\zeta_2} \cdots \x_{\zeta_k}$ for some $\ell \in M=N^\dual$ and rays $\zeta_2, \dots, \zeta_k$. The first step is to prove the following result.

\begin{claim} \label{claim:admissible_expansion}
Each element $a \in B^k\cap (I + J) $ can be written as a sum consisting of an element in $I$ and a sum consisting of elements of the form $\x_{\ell_1} \x_{\ell_2} \cdots \x_{\ell_s} \x_{\tau}$ for $s\in \N$, $\tau \in \Sigma_{k-s}$, and $\ell_1, \dots, \ell_s \in M$.
\end{claim}

An expansion of the form described in the claim for an element $a \in B^k\cap (I + J)$ will be called \emph{admissible} in the sequel.

\begin{proof}
It will be enough to prove the statement for elements $a$ of the form $\x_{\ell_1} \cdots \x_{\ell_s} \x_{\zeta_{s+1}}^{\kappa_{s+1}} \cdots \x_{\zeta_r}^{\kappa_r}$ with $\kappa_{s+1} \geq \dots \geq \kappa_{r}$ and distinct and comparable rays $\zeta_{s+1}, \dots, \zeta_r$ that form the rays of a cone $\tau$ of dimension $r-s$ in $\Sigma$. Actually, the case $s=1$ will already give the result, but considering arbitrary values for $s$ allows to proceed by induction on the lexicographical order of the $k$-tuples of non-negative integer numbers.

To a term of the form $\x_{\ell_1} \cdots \x_{\ell_s} \x_{\zeta_{s+1}}^{\kappa_{s+1}} \cdots \x_{\zeta_r}^{\kappa_r}$ with $\kappa_{s+1} \geq \dots \geq \kappa_r$ we associate the $k$-tuple of integers $(\kappa_{s+1}, \dots, \kappa_r, 0, \dots, 0)$ with $k-r+s$ zero terms. We now show that an element as above with $\kappa_{s+1} > 1$ can be rewritten as the sum of an element of $I$ and a sum of terms of the above form having a $k$-tuple with strictly smaller lexicographical order.

So suppose $\kappa_{s+1} > 1$. We take a linear form $l$ which takes the value one on $\e_{\zeta_{s+1}}$ and $0$ on the other rays $\zeta_{s+2}, \dots, \zeta_r$. This is doable since $\tau$ is simplicial. Using
\[ \x_l = \x_{\zeta_{s+1}} + \sum_{\rho \in \Sigma_1 \\ \rho\not\subfaceeq \tau} l(\e_\rho)\x_\rho, \]
we get
\[ \begin{split} \x_{\ell_1}\! \cdots \x_{\ell_s} \x_{\zeta_{s+1}}^{\kappa_{s+1}} \cdots \x_{\zeta_r}^{\kappa_r}
  =&\: \x_{\ell_1}\! \cdots \x_{\ell_s} \x_{l} \x_{\zeta_{s+1}}^{\kappa_{s+1}-1} \x_{\zeta_{s+2}}^{\kappa_{s+2}} \cdots \x_{\zeta_r}^{\kappa_r} \\
  &\quad- \sum_{\rho \in \Sigma_1 \\ \rho \not\subfaceeq \tau} l(\e_\rho) \x_{\ell_1} \cdots \x_{\ell_s} \x_{\rho} \x_{\zeta_{s+1}}^{\kappa_{s+1}-1} \x_{\zeta_{s+2}}^{\kappa_{s+2}} \cdots \x_{\zeta_r}^{\kappa_r}.
\end{split} \]
Each term in the right hand side is either in $I$ or is of the form described above with a lower lexicographic order. Proceeding by induction, we get the claim.
\end{proof}

\subsubsection{An auxilary filtration} Using admissible expansions, we now introduce an increasing filtration $\F_\bul$ on $Z^k \cap (I +J)$ as follows.

First, for each $s>0$, denote by $\filt_s$ the group generated by the elements of $B^k$ of the form $\x_{\ell_1} \cdots \x_{\ell_s} \x_\tau$ for some linear forms $\ell_1, \dots, \ell_s \in M$ and for some $\tau \in \Sigma_{k-s}$. Moreover, for $s=0$, define the subgroup $\filt_0 \subseteq \filt_1$ as the one generated by the elements of the form $\x_\ell \x_\tau$ for some $\tau \in \Sigma_{k-1}$ and some linear form $\ell$ orthogonal to $\tau$.

For any $t\geq 0$, let $\H_t \coloneqq \filt_0 + \dots + \filt_t$ and define $\F_t \coloneqq Z^k \cap (\H_t + I)$. Note that $\H_t$ is the subgroup generated by the elements which admit an admissible expansion having only terms $\x_{\ell_1} \x_{\ell_2} \cdots \x_{\ell_s} \x_{\tau}$ with $s\leq t$. In this way, we get a filtration
\[ \F_0 \subseteq \F_1 \subseteq \cdots \subseteq \F_{k-1} \subseteq \F_k \subseteq Z^k \cap (I + J). \]

We now prove that all these inclusions are equalities.

\begin{claim} \label{claim:equality}
We have $\F_0 = \F_1 = \cdots = \F_{k-1} = \F_k= Z^k \cap (I + J)$.
\end{claim}

\begin{proof}
By Claim~\ref{claim:admissible_expansion}, we have $\F_k = Z^k \cap (I + J)$, which proves the last equality. We prove all the other equalities.

So fix $s > 0$ and let $a$ be an element of $\F_s = Z^k \cap (\H_s + I)$ admitting an admissible expansion consisting of an element of $I$ plus a sum of terms each in $\filt_1, \dots, \filt_{s-1}$ or $\filt_s$. We need to show that $a \in \F_{s-1}$.

We can assume there is a term of $a$ in this admissible expansion which lies in $\filt_s$, \ie, of the form $\x_{\ell_1} \cdots \x_{\ell_s} \x_\tau$, with $\tau \in \Sigma_{k-s}$. Otherwise, the statement $a \in \F_{s-1}$ holds trivially. Consider all the terms in the admissible expansion of $a$ which are of the form $\x_{\ell'_1} \cdots \x_{\ell'_s} \x_\tau$ with the same cone $\tau$, but with possibly different linear forms $\ell_1', \dots, \ell_s'$. We will prove that the sum of those terms that we denote by $a_\tau$ belongs to $\filt_{s-1} + I$. Applying this to any $\tau \in \Sigma_{k-s}$, we obtain $a \in \F_{s-1}$ and the claim follows.

\smallskip
For each ray $\rho$ of $\tau$, choose a linear form $\ell_{\rho, \tau}$ which takes value one on $\e_\rho$ and zero on other rays of $\tau$. This is again possible since $\Sigma$ is simplicial. Then, $\ell_1$, for instance, can be decomposed as the sum
\[ \ell_1 = l_1 + \sum_{\rho \subfaceeq \tau \\
\rho \in \Sigma_1} \ell_1(\e_\rho) \ell_{\rho, \tau} \]
with $l_1$ vanishing on $\tau$.

The observation now is that the term $\x_{l_1} \x_{\ell_2} \cdots \x_{\ell_s} \x_\tau$ belongs to $\filt_{s-1} + I$: this is by definition if $s = 1$, and for $s>1$, it is obtained by expanding the product
\[ \x_{l_1} \x_{\ell_2} \cdots \x_{\ell_s} \x_\tau = \sum_{\rho \in \Sigma_1\\
\rho \not \subfaceeq \tau} l_1(\e_\rho) \x_\rho \x_{\ell_2} \cdots \x_{\ell_s} \x_\tau \]
and by observing that each term in the right hand side is either in $I$ or is in $\filt_{s-1}$. (Note that the sum is on rays $\rho \not \subfaceeq \tau$ because $l_1$ vanishes on $\tau$.) Hence, in proving that $a_\tau$ belongs to $\filt_{s-1}+I$, we can ignore this term.

Decomposing in the same way each $\x_{\ell_i}$ which appears in a term in the initial admissible expansion of $a_\tau$, we can rewrite $a_\tau$ as a sum of terms which already belong to $\H_{s-1}+I$ plus a sum of terms of the form
\begin{equation} \label{eq:atau}
\prod_{\rho \subfaceeq \tau\\\rho\in \Sigma_1}\x_{\ell_{\rho, \tau}}^{\kappa_\rho} \x_\tau, \qquad \kappa_\rho \in \Z_{\geq 0}.
\end{equation}
Altogether, these give a new admissible expansion of $a$ in which $a_\tau$ (defined as before in the new admissible expansion) is a sum of the terms of the form in~\eqref{eq:atau}. We will work from now on with this admissible expansion of $a$.

Fix non-negative integers $\kappa_\rho$ for rays $\rho$ of $\tau$ whose sum is $s$. Denote by $\bigl[\prod_{\rho \in \tau} \x_\rho^{\kappa_\rho + 1} \bigr] a$ the coefficient of this monomial in $a$ written as sum of monomials in $B^k$. Since $a$ belongs to $Z^k$, all the monomials of $a$ are square-free. But the product $\prod_{\rho \in \tau} \x_\rho^{\kappa_\rho + 1}$ is not square-free since $s>0$. Hence, the corresponding coefficient is zero.

Consider now the monomial $\prod_{\rho \in \tau} \x_\rho^{\kappa_\rho + 1}$. We will look in which terms in the admissible expansion of $a$ it can occur. Such a monomial cannot appear in a term of the admissible expansion of $a$ which belongs to $I$. It cannot neither be in a term of $a$ of the form $\x_{\lambda_1} \cdots \x_{\lambda_{s'}} \x_\sigma$ with $\sigma \in \Sigma_{k-s'}$ with $s'<s$ nor with $\sigma \in \Sigma_{k-s}$ and $\sigma \neq \tau$.

It follows that the monomial $\prod_{\rho \in \tau} \x_\rho^{\kappa_\rho + 1}$ can only appear in the terms of $a_\tau$. More precisely, by the definition of $\ell_{\rho,\tau}$, it can only appear in each term of the form $\prod_{\rho \subfaceeq \tau}\x_{\ell_{\rho, \tau}}^{\kappa_\rho} \x_\tau$ for the chosen $\kappa_\rho$, and with coefficient one in each term. Hence, the sum of these terms have to cancel out. This proves that $a_\tau$ is zero in the new admissible expansion, which shows that $a \in \F_{s-1}$ and the claim follows.
\end{proof}

\subsubsection{End of the proof} We can now finish the proof of the lemma.

\begin{proof}[Proof of Theorem~\ref{thm:kernel_Z^k}]
Recall that $K^k$ is the set generated by elements of the form
\[ \sum_{\rho \in \Sigma_1 \\ \rho\sim\tau} \ell(\e_\rho) \x_\rho \x_\tau, \]
where $\tau \in \Sigma_{k-1}$ and $\ell$ is orthogonal to $\tau$. The following facts are then clear:
\[ \F_0 \subseteq K^k + I, \quad K^k \subseteq Z^k, \quad \text{and } I \cap Z^k = \{0\}. \]
Since $\F_0 \subseteq Z^k$, we deduce that $\F_0 \subseteq K^k$. Applying Claim~\ref{claim:equality}, we get $Z^k \cap (I + J) = \F_0$ which implies $Z^k \cap (I + J) \subseteq K^k$. The inclusion $K^k \subseteq Z^k \cap (I + J)$ is obvious. So we get $Z^k \cap (I + J) = K^k$ which concludes the proof.
\end{proof}

\section{Tropical divisors} \label{sec:divisors}

In this section, we consider a tropical fan $\Sigma$ of dimension $d$ in $N_\R$ and study divisors associated to meromorphic functions on $\Sigma$. Based on this, we define three classes of tropical fans: the class of \emph{principal}, \resp \emph{$\Q$-principal}, \resp \emph{div-faithful}, tropical fans.

\subsection{Divisors on fans} \label{sec:divisors_tropical}

Let $\Sigma$ be a fan. A \emph{divisor} of $\Sigma$ is an element $D$ of $\MW_{d-1}(\Sigma)$, that is, a Minkowski weight $D \colon \Sigma_{d-1} \to \Z$. We denote by $\Div(\Sigma)$ the group of divisors on $\Sigma$ and note that we have $\Div(\Sigma) = \MW_{d-1}(\Sigma)$. Using the terminology of Section~\ref{sec:MW}, we call a divisor $D$ \emph{effective} if the coefficients $D(\tau)$, $\tau \in \Sigma_{d-1}$, are all non-negative, and we say it is \emph{reduced} if all its weights are equal to zero or one. In the case of a non-zero reduced divisor $D$, we can identify $D$ with its support $\Delta$ that we view as a tropical subfan of $\Sigma$. If $D$ is trivial, that is, the weights in $D$ are all equal to zero, then, the corresponding support is $\Delta=\emptyset$, although, this is not a subfan of $\Sigma$.

Working with rational coefficients, we set $\Div!_{\Q}(\Sigma) = \MW_{d-1}(\Sigma, \Q)$.

\subsection{Principal divisor associated to meromorphic functions on tropical fans} \label{sec:principal_divisors_tropical}

Let $\Sigma$ be a tropical fan with underlying orientation $\omega=\omega!_\Sigma \colon \Sigma_d \to \Z$. Let $f\in \mer(\Sigma)$ be a meromorphic function on $\Sigma$. We denote by $f_\eta$ the linear form induced by $f$ on $N_{\eta,\R}$, for each $\eta \in \Sigma$.

Let $\tau$ be a face of codimension one in $\Sigma$. The \emph{order of vanishing of $f$ along $\tau$} denoted by $\ord_\tau(f)$ is defined as
\[ \ord_\tau(f) \coloneqq -\sum_{\sigma \ssupface \tau} \omega(\sigma)f_\sigma(\nvect_{\sigma/\tau}) + f_\tau\Bigl(\sum_{\sigma \ssupface \tau} \omega(\sigma) \nvect_{\sigma/\tau}\Bigr) \]
with the sums running over all the cones $\sigma \in \Sigma_d$ that contain $\tau$. Note that the last term is well-defined since the sum belongs to $N_\tau$ by the balancing condition. Moreover, if $f$ is linear, then we have $\div(f)=0$.

\begin{prop}
Notations as above, the order of vanishing $\ord_\tau(f)$ is well-defined. That is, $\ord_\tau(f)$ is independent of the choice of normal vectors $\nvect_{\sigma/\tau} \in N_\sigma$.
\end{prop}

\begin{proof}
For each pair $\sigma \ssupface \tau$, two different choices $\nvect_{\sigma/\tau}$ and $\nvect'_{\sigma/\tau}$ of normal vectors differ by a vector in $N_\tau$. It follows that
\begin{align*}
-\sum_{\sigma \ssupface \tau} \omega(\sigma)f_\sigma(\nvect_{\sigma/\tau}) &+ f_\tau\Bigl(\sum_{\sigma \ssupface \tau} \omega(\sigma)\nvect_{\sigma/\tau}\Bigr) + \sum_{\sigma \ssupface \tau} \omega(\sigma)f_\sigma(\nvect'_{\sigma/\tau}) - f_\tau\Bigl(\sum_{\sigma \ssupface \tau} \omega(\sigma)\nvect'_{\sigma/\tau}\Bigr)\\
&= - \sum_{\sigma \ssupface \tau} \omega(\sigma)f_\sigma(\nvect_{\sigma/\tau} - \nvect'_{\sigma/\tau}) + f_\tau\Bigl(\sum_{\sigma \ssupface \tau} \omega(\sigma)\bigl(\nvect_{\sigma/\tau} - \nvect'_{\sigma/\tau}\bigr)\Bigr)\\
&= - \sum_{\sigma \ssupface \tau} \omega(\sigma)f_\tau(\nvect_{\sigma/\tau} - \nvect'_{\sigma/\tau}) + f_\tau\Bigl(\sum_{\sigma \ssupface \tau} \omega(\sigma)\bigl(\nvect_{\sigma/\tau} - \nvect'_{\sigma/\tau}\bigr)\Bigr) =0. \qedhere
\end{align*}
\end{proof}

The order of vanishing function gives a weight function $\ord(f)\colon \Sigma_{d-1} \to \Z$. If $\ord(f)$ is constant equal to zero, we associate to $f$ the empty divisor. Otherwise, we associate to $f$ the data of the pair $(\Delta, w\coloneqq\ord(f))$ with $\Delta$ the fan defined by the support of $\ord(f)$, that is, by the cones $\tau$ of dimension $d-1$ in $\Sigma$ with $\ord_\tau(f) \neq 0$, and with the weight function $w \colon \Delta_{d-1} \to \Z\setminus\{0\}$ given by $w(\tau) = \ord_\tau(f)$. We have the following result, see \cite{AR10}*{Section~3}.

\begin{prop} \label{prop:balancing_divisor}
Notations as above, $\ord(f)$ is a divisor in $\Sigma$. In the case the divisor is nontrivial, it follows that $(\Delta, w)$ is a tropical fan.
\end{prop}

\begin{defi}[Principal divisors]
Let $\Sigma$ be a tropical fan. For any meromorphic function on $\Sigma$, we denote by $\div(f)$ the divisor associated to $f$. Such divisors are called \emph{principal}. Principal divisors form a subgroup of $\Div(\Sigma)$ that we denote by $\Prin(\Sigma)$. The vector subspace of $\Div!_{\Q}(\Sigma)$ generated by $\Prin(\Sigma)$ is denoted by $\Prin!_{\Q}(\Sigma)$. Its elements are called $\Q$-principal.
\end{defi}

\begin{defi}[Holomorphic functions]
A meromorphic function $f$ on $\Sigma$ is called \emph{holomorphic} in the sequel if the principal divisor $\div(f)$ is effective.
\end{defi}

Note that if $\Sigma$ is normal, then a meromorphic function $f$ on $\Sigma$ is holomorphic if and only if it is concave in codimension one. This property might fail to hold for fans which are not normal.

\subsection{Principal and div-faithful tropical fans} \label{sec:definition_principa_divfaithful}

We now introduce tropical fans on which divisors behave nicely.

Let $\Sigma$ be a tropical fan with underlying orientation $\omega= \omega!_\Sigma$, and denote by $\omega^\sigma =\omega!_{\Sigma^\sigma}$ the induced orientation on the star fan $\Sigma^\sigma$, $\sigma\in \Sigma$.

\begin{defi} \label{defi:principal_div-faithful}
Let $\Sigma$ be a tropical fan and let $\eta$ be a cone in $\Sigma$.
\begin{itemize}
\item We say that $\Sigma$ is \emph{principal at $\eta$} if any divisor on $\Sigma^\eta$ is principal for the orientation given by $\omega^\eta$. We call the tropical fan $\Sigma$ \emph{principal} if $\Sigma$ is principal at any cone $\eta \in \Sigma$.

\item We say that $\Sigma$ is \emph{$\Q$-principal at $\eta$} if for any divisor $D$ on $\Sigma^\eta$, $D$ is $\Q$-principal, \ie, an integer multiple $aD$ for $a\in \Z \setminus\{0\}$ is principal. We say that $\Sigma$ is \emph{$\Q$-principal} if $\Sigma$ is $\Q$-principal at any cone $\eta \in \Sigma$.

\item We say that $\Sigma$ is \emph{divisorially faithful} or simply \emph{div-faithful} \emph{at $\eta$} if for any meromorphic function $f$ on $\Sigma^\eta$, if $\div(f)$ is trivial in $\Div(\Sigma^\eta)$, then $f$ is a linear function on $\Sigma^\eta$. We say the tropical fan $\Sigma$ is \emph{div-faithful} if $\Sigma$ is div-faithful at any cone $\eta\in \Sigma$. \qedhere
\end{itemize}
\end{defi}

The importance of div-faithful property in our work is based on the fact that tropical modifications behave very well on div-faithful tropical fans, \cf Section \ref{subsec:trop_mod_div-faithful}.

\begin{prop} \label{prop:div-faithful_local}
The properties of being principal, being $\Q$-principal, and being div-faithful are all local.
\end{prop}

\begin{proof}
The statement is tautological.
\end{proof}

\begin{remark}
We emphasize that being principal, \resp div-faithful, at $\conezero$ does not in general imply that the fan is principal, \resp div-faithful, see Examples \ref{ex:non-principal_fan} and \ref{ex:PD_chow_not_smooth}.
\end{remark}

\subsection{Characterization of saturation}

In the definition of div-faithfulness, we look at meromorphic functions which are linear. There is a subtle difference with integral linear functions, \ie, those functions induced by $M$, as we show in Examples \ref{ex:F1_vs_N} and \ref{ex:non_saturated_smooth_fan}.

Linear meromorphic functions on a rational fan $\Sigma$ are exactly those that have a multiple in $M$, and they are in correspondence with $((N \cap \Sigma) \otimes \Z)^\dual \supseteq M$. The last inclusion is an equality exactly when $\Sigma$ is saturated at $\conezero$. In the case $\Sigma$ is unimodular, the space $Z^1(\Sigma) = \bigoplus_{\rho\in \Sigma_1} \Z\x_\rho$ coincides with the space of meromorphic functions on $\Sigma$. Hence, $A^1(\Sigma) \simeq \rquot{\mer(\Sigma)}{M}$. The following proposition follows from this discussion.

\begin{prop} \label{prop:characterization_of_saturation}
Let $\Sigma$ be a unimodular fan. The following statements are equivalent.
\begin{enumerate}
\item $\Sigma$ is saturated at $\conezero$.
\item Meromorphic functions on $\Sigma$ that are linear are all induced by elements of $M$.
\item $A^1(\Sigma)$ has no torsion.
\end{enumerate}
\end{prop}

\subsection{Characterization of principal, $\Q$-principal, and div-faithful tropical fans} \label{subsec:characterization_principal_div-faithful}

We now provide a characterization of principal and $\Q$-principal, \resp divisorially faithful, tropical fans in the case the tropical fan is unimodular, \resp saturated and unimodular. This will be given via the cycle class map
\[ \cycl\colon A^1(\Sigma) \to \MW_{d-1}(\Sigma) \simeq \Div(\Sigma). \]

\begin{prop}[Cycle class map for divisors] \label{prop:map_cl}
Let $\Sigma$ be a unimodular fan. Consider an element $\alpha \in A^1(\Sigma)$ and take a representative $\alpha = \sum_{\zeta\in \Sigma_1} a_\zeta x_\zeta$ for coefficients $a_\zeta \in \Z$. Let $f$ be the meromorphic function on $\Sigma$ which takes value $a_\zeta$ at $\e_\zeta$, for any ray $\zeta \in \Sigma$. Then, we have $\cycl(\alpha) = -\div(f)$.
\end{prop}

\begin{proof}
Notations as above, we need to show that for each $\tau \in \Sigma_{d-1}$, we have the equality
\[ \deg(\alpha \cdot x_\tau) = -\ord_\tau(f). \]
Let $m$ be an element of $M = N^\dual$ whose restriction to $\tau$ coincides with $f\rest\tau$. The element $x_m \coloneqq \sum_{\zeta\in \Sigma_1} m(\e_\zeta)x_\zeta$ vanishes in $A^1(\Sigma)$. Replacing $\alpha$ by $\alpha -x_m$ if necessary, we can assume that $a_\zeta=0$ for any ray $\zeta$ in $\tau$. The proposition now follows by observing that
\[ \deg(\alpha \cdot x_\tau) = \sum_{\zeta\in \Sigma_1 \\ \zeta\sim \tau} a_\zeta \deg(x_\zeta x_\tau) = \sum_{\zeta\in \Sigma_1 \\ \zeta\sim \tau} a_\zeta \omega!_\Sigma(\zeta\vee \tau) = -\ord_\tau(f-m) = -\ord_\tau(f). \qedhere \]
\end{proof}

We obtain the following important result.

\begin{thm}[Characterization of principal, $\Q$-principal, and div-faithful tropical fans] \label{thm:char_div-faithful}
Let $\Sigma$ be a unimodular tropical fan. Consider the cycle class map $\cycl\colon A^1(\Sigma) \to \Div(\Sigma)$. We have
\begin{itemize}
\item $\Sigma$ is principal at $\conezero$ if and only if $\cycl$ is surjective.
\item $\Sigma$ is $\Q$-principal at $\conezero$ if and only if $\cycl!_{\Q} \colon A^1(\Sigma, \Q) \to \Div!_\Q(\Sigma) $ is surjective.
\item $\cycl$ is injective if and only if $\Sigma$ is both saturated at $\conezero$ and div-faithful at $\conezero$.
\end{itemize}
\end{thm}

\begin{proof}
The first and second part follow directly from Proposition~\ref{prop:map_cl}. For the third part, the injectivity of $\cycl$ implies $A^1(\Sigma)$ has no torsion, and the saturation follows from Proposition~\ref{prop:characterization_of_saturation}. Saturation at $\conezero$ allows to identify $A^1(\Sigma)$ with $\mer(\Sigma)$ quotiented by elements of $M$, and one can conclude thanks to Proposition~\ref{prop:map_cl}.
\end{proof}

\begin{remark} \label{rem:independent_exterior_lattice}
Note that being \textnormal{($\Q$-)}principal and being div-faithful are independent of the form of the lattice outside the support of the fan. By this, we mean that a tropical fan $\Sigma$ in $N_\R$ is \textnormal{($\Q$-)}principal, \resp div-faithful, if and only if for any full rank lattice $N'$ in $N_\R$ which verifies $\supp\Sigma \cap N = \supp\Sigma \cap N'$, the fan $\Sigma$ endowed with the integral structure induced by the new lattice $N'$ is \textnormal{($\Q$-)}principal, \resp div-faithful.

If a predicate of tropical fans verifies this property, replacing the lattice $N$ if necessary, there is no harm in assuming that the fan is saturated at $\conezero$. Note however that, in general, we cannot assume global saturation, \ie, saturation at all faces of $\Sigma$, see Example \ref{ex:non_saturated_smooth_fan}.
\end{remark}

\section{Tropical modification} \label{subsec:tropical_modification}

We recall the definition of tropical modifications, and describe their star fans. A survey of results and references related to the concept can be found in~\cites{BIMS, Kal15}.

\subsection{Definition of tropical modifications}

Let $f$ be a meromorphic function on a tropical fan $\Sigma$. Consider the principal divisor $D=\div(f)$. Denote by $\Delta$ the corresponding subfan of $\Sigma$, \ie, the support of $\ord(f)$, endowed with the orientation $\omega!_\Delta = \ord(f)$. We allow the case the divisor $\div(f)$ is trivial, in which case $\Delta$ will be empty. Note that if $f$ is holomorphic, then $\Delta$ will be effective.

We define the tropical modification of $\Sigma$ with center $\Delta$ induced by the meromorphic function $f$. This will be a fan in $\~N_\R \simeq \R^{n+1}=\R^{n} \times \R$, for the lattice $\~N \coloneqq N \times \Z$, that we will denote by $\tropmod{f}{\Sigma}$.

First assume that $\Delta$ is non-empty. Consider the graph of $f$ which is the map $\Gamma!_f$ defined as
\[ \begin{array}{rccc}
\Gamma!_f\colon & \supp\Sigma & \longrightarrow & \~N_\R=N_\R \times \R, \\
                & x           & \longmapsto     & (x, f(x)).
\end{array} \]
For each cone $\sigma$ of $\Sigma$, we consider the cone $\basetm\sigma$ in $\~N_\R$ which is the image of $\sigma$ by $\Gamma!_f$, \ie,
\[ \basetm\sigma\coloneqq \Gamma!_f(\sigma) \subset \~N_\R. \]
Moreover, to each face $\delta$ of $\Delta$, we associate the face
\[ \uptm\delta \coloneqq \basetm\delta + \R_{\geq 0} \etm \]
where $\etm = (0, 1)\in N_\R \times \R$. \\ (Here, $0$ in $\etm=(0,1)$ refers to the origin in $N_\R$).

\begin{figure}[t]
  \begin{tikzpicture}
    \begin{scope}[fill opacity=.2, z={(-.1,.1,-1)}]
      \renewcommand\a{.707}
      \renewcommand\t{{90+atan(.707)}}

      \begin{scope}[canvas is xz plane at y=0]
        \fill[blue] (0,0) -- (2,0) arc (0:90:2) -- cycle;
        \fill[green] (0,0) -- (2,0) arc (0:-135:2) -- cycle;
        \fill[violet] (0,0) -- (0,2) arc (90:225:2) -- cycle;
        \draw[very thin] (-2,0) -- (2,0);
        \draw[very thin] (0,-2) -- (0,2);
        \draw[very thin] (45:-2) -- (45:2);
        \draw[->] (0,0) -- (2.3,0) node[right, opacity=1] {$x$};
        \draw[->] (0,0) -- (0,2.3) node[pos=1.15, opacity=1] {$y$};
        \draw[very thick] (O) edge (45:-2) edge (2,0) edge (0,2);
      \end{scope}
    \end{scope}
  \end{tikzpicture} \qquad
  \begin{tikzpicture}
    \begin{scope}[yscale=.8, fill opacity=.2, z={(-.1,.1,-1)}]
      \renewcommand\a{.707}
      \renewcommand\t{{90+atan(.707)}}

      \begin{scope}[canvas is xz plane at y=0]
        \fill[blue] (0,0) -- (2,0) arc (0:90:2) -- cycle;
      \end{scope}

      \coordinate (A) at (${2/sqrt(3)}*(-1,-1,-1)$);

      \fill[red] (0,0,0) -- (2,0,0) --++ (0,2,0) -- (0,2,0) -- cycle;
      \fill[yellow, opacity=.4] (0,0,0) -- (0,0,2) --++ (0,2,0) -- (0,2,0) -- cycle;
      \fill[brown] (0,0,0) -- (A) --++ (0,{2+2/sqrt(3)},0) -- (0,2,0) -- cycle;
      \draw[very thick, -latex, opacity=1] (0,0) --++ (0,.8) node[right=-2pt, -] {$\etm$};
      \draw[] (0,0,0) --++ (0,2);

      \fill[green, y={(0,-\a,-\a)}] (0,0) -- (2,0) arc (0:\t:2) -- cycle;
      \fill[violet, y={(-\a,-\a,0)}, x={(0,0,1)}] (0,0) -- (2,0) arc (0:\t:2) -- cycle;

      \draw[very thick] (0,0) edge (A) edge (0,0,2) edge (2,0,0) edge (0,0,2);
      \draw[very thin] (0,0) edge (2*\a,0,2*\a) edge (0,-2*\a,-2*\a) edge (-2*\a,-2*\a,0);
    \end{scope}
  \end{tikzpicture}

  \caption{A tropical modification along the function $\min(x,y,0)$. The divisor is in bold. \label{fig:tropical_modification}}
\end{figure}
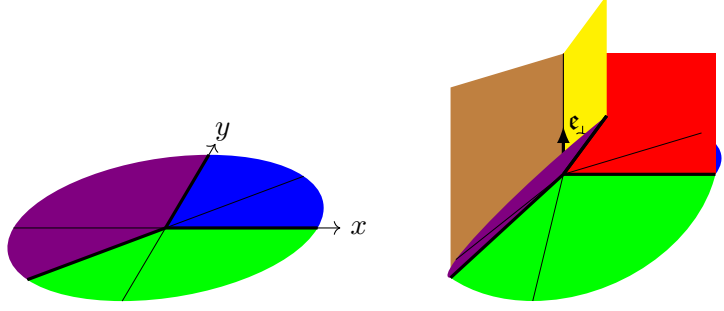

The \emph{tropical modification of\/ $\Sigma$ along $\Delta$ with respect to $f$}, or simply the tropical modification of $\Sigma$ along $\Delta$ if the other terms are understood from the context, is the fan $\~\Sigma=\tropmod{f}{\Sigma}$ in $\~N_\R\simeq \R^{n+1}$ defined as
\[ \~\Sigma=\tropmod{f}{\Sigma} \coloneqq \bigl\{\,\basetm\sigma \mid \sigma\in\Sigma\,\bigr\} \cup \bigl\{\,\uptm\delta \mid \delta\in\Delta\,\bigr\}.
\]
Note that we have $\dims{\uptm\delta} =\dims\delta+1$, for $\delta\in \Delta$, and $\dims{\basetm \sigma} = \dims\sigma$, for $\sigma \in \Sigma$. This shows that for any non-negative integer $k$,
\[ \~\Sigma_k = \bigl\{\basetm\sigma \mid \sigma\in\Sigma_k\bigr\} \cup \bigl\{\,\uptm\delta \mid \delta\in\Delta_{k-1}\,\bigr\}. \]

We endow $\~\Sigma$ with the map $\~\omega \colon \~\Sigma_d \to \Z$ defined as follows. For $\sigma \in \Sigma_d$, we set $\~\omega(\basetm\sigma) \coloneqq\omega!_\Sigma(\sigma)$. For $\delta\in \Delta_{d-1}$, we set $\~\omega(\uptm\delta) \coloneqq\omega!_\Delta(\delta)$.

\begin{prop} \label{prop:balancing_trop_modif}
Notations as above, $\~\omega$ is an orientation of\/ $\~\Sigma$, that is, the tropical modification $\tropmod{f}{\Sigma}$ endowed with $\~\omega$ is a tropical fan.

Moreover, we have a natural projection map
\[ \prtm \colon \supp{\tropmod{f}{\Sigma}} \to \supp\Sigma \]
which is conewise integral linear.
\end{prop}

Before going through the proof, we make some remarks and introduce some notations. First, we observe that the fan $\tropmod{f}\Sigma$ is rational with respect to the lattice $\~N$, and that for each face $\sigma \in \Sigma$, the lattice $\~N_{\basetm\sigma}$ can be identified with the image $\Gamma!_{f_\sigma}(N_\sigma)$ where, as before, $f_\sigma \in N_\sigma^\dual$ denotes the linear form which coincides with $f$ on $\sigma$. This shows that for an inclusion of cones $\tau \ssubface \sigma$, we can pick
\[ \nvect_{\basetm\sigma/\basetm\tau} = (\nvect_{\sigma/\tau},f_\sigma(\nvect_{\sigma/\tau})) \in \~N = N \times \Z, \]
Second, we observe that for each $\delta\in\Delta$, the lattice $\~N_{\uptm\delta}$ can be identified with $\~N_{\basetm\delta}\times\Z \simeq N_\delta \times \Z$, and that we can choose $\nvect_{\uptm\delta/\basetm\delta} = \etm$. Finally, for an inclusion of cones $\tau \ssupface \delta$ in $\Delta$, we can set $\nvect_{\uptm\tau/\uptm\delta} = (\nvect_{\tau/\delta},0) \in \~N$.

\begin{proof}[Proof of Proposition~\ref{prop:balancing_trop_modif}]
We need to prove the balancing condition around any codimension one face of $\tropmod{f}\Sigma$. These are of two kinds: namely, faces of the form $\basetm\tau$ for $\tau \in \Sigma_{d-1}$ and those of the form $\uptm\delta$ for $\delta \in \Delta_{d-2}$.

\smallskip
First, let $\tau \in \Sigma_{d-1}$, and consider the codimension one face $\basetm\tau$ of $\tropmod{f}\Sigma$. Two cases happen:

\begin{itemize}[leftmargin=1em]
\item Either, $\tau \in \Delta$, in which case, the $d$-dimensional faces around $\basetm\tau$ are of the form $\basetm\sigma$ for $\sigma \ssupface \tau$ in $\Sigma$ as well as the face $\uptm\tau$. The balancing condition in this case amounts to showing that the vector
\begin{align*}
\omega!_\Delta(\tau)\nvect_{\uptm\tau/\tau} + \sum_{\sigma\ssupface \tau\\ \sigma\in \Sigma} \omega!_\Sigma(\sigma)\nvect_{\basetm\sigma/\basetm\tau}
  &= \omega!_\Delta(\tau)\etm + \sum_{\sigma\ssupface \tau\\ \sigma\in \Sigma} \omega!_\Sigma(\sigma)(\nvect_{\sigma/\tau}, f_\sigma(\nvect_{\sigma/\tau}))\\
  &= (0,\omega!_\Delta(\tau))+\Bigl(\,\sum_{\sigma\ssupface \tau\\ \sigma\in \Sigma} \omega!_\Sigma(\sigma)\nvect_{\sigma/\tau}, \sum_{\sigma\ssupface \tau\\ \sigma\in \Sigma} \omega!_\Sigma(\sigma)f_\sigma(\nvect_{\sigma/\tau}) \Bigr)
\end{align*}
belongs to $\~N_{\basetm\tau}$. Since $\omega!_\Delta(\tau)= \ord_\tau(f)$, the term on the right hand side of the above equality becomes equal to
\[ \Bigl(\,\sum_{\sigma\ssupface \tau\\ \sigma\in \Sigma} \omega!_\Sigma(\sigma)\nvect_{\sigma/\tau}, f_\tau\bigl(\sum_{\sigma\ssupface \tau\\ \sigma\in \Sigma}\omega!_\Sigma(\sigma) \nvect_{\sigma/\tau}\bigr) \Bigr), \]
which by balancing condition at $\tau$ in $\Sigma$ belongs to $\~N_{\basetm\tau}= \Gamma!_{f_\tau}(N_\tau)$.

\item Or, we have $\tau \notin \Delta$, \ie, $\ord_\tau(f)=0$. In this case, the facets around $\basetm\tau$ are of the form $\basetm\sigma$ for $\sigma \ssupface \tau$ in $\Sigma$, and we get
\begin{align*}
\sum_{\sigma\ssupface \tau\\ \sigma\in \Sigma} \omega!_\Sigma(\sigma)\nvect_{\basetm\sigma/\basetm\tau}
  &= \sum_{\sigma\ssupface \tau\\ \sigma\in \Sigma} \omega!_\Sigma(\sigma)(\nvect_{\sigma/\tau}, f_\sigma(\nvect_{\sigma/\tau}))=\Bigl(\sum_{\sigma\ssupface \tau \\ \sigma\in \Sigma} \omega!_\Sigma(\sigma)\nvect_{\sigma/\tau}, \sum_{\sigma\ssupface \tau \\ \sigma\in \Sigma} \omega!_\Sigma(\sigma)f_\sigma(\nvect_{\sigma/\tau}) \Bigr)\\
  &= \Bigl(\sum_{\sigma\ssupface \tau\\ \sigma\in \Sigma} \omega!_\Sigma(\sigma)\nvect_{\sigma/\tau}, f_\tau\bigl(\sum_{\sigma\ssupface \tau\\ \sigma\in \Sigma} \omega!_\Sigma(\sigma)\nvect_{\sigma/\tau}\bigr) \Bigr),
\end{align*}
which again belongs to $\~N_{\basetm\tau}= \Gamma!_{f_\tau}(N_\tau)$ by the balancing condition.
\end{itemize}

It remains to check the balancing condition around a codimension one face of the form $\uptm\delta$ in $\tropmod{f}\Sigma$ with $\delta \in \Delta_{d-2}$. Facets around $\uptm\delta$ are all the cones $\uptm\tau$ for $\tau \in \Delta_{d-1}$ and $\tau \ssupface \delta$. Using the balancing condition in $\Delta =\div(f)$ around $\delta$, Proposition~\ref{prop:balancing_divisor}, we see that the sum
\[ \sum_{\tau\ssupface \delta\\
\tau\in \Delta} \omega!_\Delta(\tau)\nvect_{\uptm\tau/\uptm\delta} = \Bigl(\sum_{\tau\ssupface \delta \\ \tau \in \Delta} \omega!_\Delta(\tau)\nvect_{\tau/\delta}, 0\Bigr) \]
belongs to $N_\delta \times \Z = \~N_{\uptm\delta}$, and the assertion follows.

\smallskip
The second statement is straightforward.
\end{proof}

We now consider the case where the divisor is trivial, that is, $\Delta =\emptyset$. The tropical modification of $\Sigma$ with respect to $f$ still has a meaning and is equal to the graph of $f$. We call this case \emph{degenerate}. The faces of the tropical modification will be in one-to-one correspondence with faces of $\Sigma$, and the orientation is preserved. However, unless $f$ is a integral linear form on $\Sigma$, the tropical modification is not isomorphic to $\Sigma$. We refer to Examples \ref{ex:nontrivial_degenerate_tropical_modification} and \ref{ex:nontrivial_linear_tropical_modification} in Section~\ref{sec:examples} which explains this phenomenon. On the contrary, if $f$ is integral linear, then we just obtain the image of $\Sigma$ by a linear map, and in this case, the fans $\Sigma$ and $\tropmod{f}{\Sigma}$ are isomorphic. In particular, if $\Sigma$ is div-faithful and saturated at $\conezero$, then the vanishing of $\div(f)$ implies that $f$ is integral linear, and in this case, the tropical modification becomes isomorphic to $\Sigma$. Unless otherwise stated, in this article we allow tropical modifications to be degenerate.

By the description above, we get the following result.

\begin{prop}
Assume $\Sigma$ is unimodular. The tropical modification $\~\Sigma$ remains unimodular.
\end{prop}

\subsection{The case of tropical fans which are principal, div-faithful and saturated at $\conezero$}

Assume $\Sigma$ is a tropical fan which is principal, div-faithful and saturated at the cone $\conezero \in \Sigma$. In this case, any divisor $D$ is the divisor of a meromorphic function $f$ on $\Sigma$, and in addition, if $g$ is another meromorphic function on $\Sigma$ such that $\div(g) = D$, then $g-f$ is integral linear, that is, it is the restriction to $\supp{\Sigma}$ of an element $\ell\in M$. Therefore, the two tropical modifications $\tropmod{f}{\Sigma}$ and $\tropmod{g}{\Sigma}$ are isomorphic via the affine map which sends the point $x \in \R^{n+1}$ to the point $x + \ell(\prtm(x))\etm$. This means that, working modulo isomorphisms, we can talk about \emph{the tropical modification of\/ $(\Sigma,\omega)$ along a divisor $D$}. Replacing $D$ with the tropical fan $\Delta$ on its support, we denote this by $\tropmod{\Delta}{\Sigma}$.

\subsection{Star fans of a tropical modification} \label{sec:star_fans_tropical_modifications}

For future use, we record here a description of the star fans of a tropical modification. Let $\Sigma$ be a tropical fan and let $\sigma \in \Sigma$. Let $f$ be a meromorphic function on $\Sigma$. Then, $f$ induces a meromorphic function $f^\sigma$ on the star fan $\Sigma^\sigma$ defined as follows. Let $m \in M$ be a linear map that coincide with $f$ on $\sigma$. Then $f-m$ is zero on $\sigma$. This means that we can restrict $f-m$ on $\Sigma^\sigma$. We denote the restriction by $f^\sigma$. The function $f^\sigma$ is meromorphic on $\Sigma^\sigma$. Note that although $f^\sigma$ depends on the choice of $m$, it is well-defined up to an element of $M^{\sigma} =(N^\sigma)^\dual$. This is enough for our purpose, that is why, abusing the terminology, we sometimes call $f^\sigma$ \emph{the} meromorphic function induced by $f$ on $\Sigma^\sigma$. The Minkowski weight $\div(f)^\sigma$ induced by $\div(f)$ on $\Sigma^\sigma$ coincides with $\div(f^\sigma)$. This implies that if $f$ is holomorphic on $\Sigma$, then $f^\sigma$ is holomorphic on $\Sigma^\sigma$.

\begin{prop} \label{prop:star_fans_tropical_modifications}
Let $\Sigma$ be a tropical fan and let $f$ be meromorphic on $\Sigma$. Let $\Delta$ be the tropical fan defined by $\div(f)$. Set $\~\Sigma = \tropmod{f}{\Sigma}$. Then, we have the following description of the star fans of\/ $\~\Sigma$.
\begin{itemize}
\item If $\delta \in \Delta$, then $\~\Sigma^{\uptm\delta} \simeq \Delta^\delta$.
\item If $\sigma \in \Delta$, then $\~\Sigma^{\basetm\sigma} \simeq \tropmod{f^\sigma}{\Sigma^\sigma}$ where $f^\sigma$ is the meromorphic function induced by $f$ on $\Sigma^\sigma$.
\item If $\sigma \in \Sigma\setminus\Delta$, then once again we have $\~\Sigma^{\basetm\sigma} \simeq \tropmod{f^\sigma}{\Sigma^\sigma}$. However, this time the tropical modification is degenerate. In particular, if\/ $\Sigma$ is div-faithful and saturated (at $\sigma$), then $\~\Sigma^{\basetm\sigma}$ is isomorphic to $\Sigma^\sigma$.
\end{itemize}
\end{prop}

\begin{proof}
The proof is a direct verification.
\end{proof}

\section{Behavior of Chow rings under tropical modifications} \label{sec:chow_tropmod}

We describe the behavior of Chow rings with respect to tropical modifications.

\subsection{Surjection}

Let $\Sigma$ be a simplicial tropical fan of dimension $d$ and let $\Delta$ be the tropical fan associated to a meromorphic function on $\Sigma$. Let $\~\Sigma = \tropmod{f}{\Sigma,\omega}$ be the tropical fan obtained by the tropical modification of $\Sigma$ along $\Delta$ with respect to $f$. Denote by $\~\omega$ the induced orientation of $\~\Sigma$.

\begin{prop} \label{prop:chow_ring_tropical_modification}
Notations as above, we have a natural surjective morphism of rings
\[ \Psi\colon A^\bul(\Sigma) \twoheadrightarrow A^\bul(\~\Sigma). \]
By Duality Theorem~\ref{thm:duality_A_MW}, this induces an injective morphism
\[ \MW_\bul(\~\Sigma) \hookrightarrow \MW_\bul(\Sigma). \]
\end{prop}

The map $\MW_k(\~\Sigma) \to \MW_k(\Sigma)$ is obtained as follows: to each face of $\Sigma_k$, we associate the weight of the corresponding face in $\~\Sigma_k$. The map $\Psi$ is introduced below.

We need some preparation before giving the proof. Following the notations of Section \ref{subsec:tropical_modification}, for any ray $\zeta \in \Sigma_1$, we denotes by $\basetm\zeta$ the corresponding ray in $\~\Sigma_1$. We assume for now that $\Delta$ is nontrivial and denote by $\rho= \R_{\geq 0}\etm$ the special ray of $\~\Sigma$, which is the unique ray of $\~\Sigma$ that does not come from $\Sigma$.

We define first the map $\Psi$ on the level of polynomial rings by
\begin{align*}
\Psi\colon \Z[\x_\zeta \mid \zeta\in\Sigma_1] &\to \Z[\x_{\zeta} \mid \zeta\in\~\Sigma_1]\\
\x_\zeta &\mapsto \x_{\basetm\zeta}.
\end{align*}
Let $\ssI_\Sigma$, $\ssJ_\Sigma$, $\ssI_{\~\Sigma}$ and $\ssJ_{\~\Sigma}$ be the ideals appearing in the definition of the Chow rings of $\Sigma$ and $\~\Sigma$.

\begin{lemma} \label{lem:map_psi}
We have $\Psi(\ssI_\Sigma) \subseteq \ssI_{\~\Sigma}$ and $\Psi(\ssJ_\Sigma) \subseteq \ssJ_{\~\Sigma}$.
\end{lemma}

\begin{proof}
Clearly, $\Psi(\ssI_{\Sigma}) \subseteq \ssI_{\~\Sigma}$. If $m \in M=\ssM_{\Sigma}$ is a linear form on the ambient space of $\Sigma$, and if $\prtm$ denotes the projection associated to the tropical modification, since $\~m\coloneqq\prtm^*(m) \in \~M$ is zero on the special ray $\rho$, we get that
\[ \Psi\Bigl(\sum_{\zeta \in \Sigma_1} m(\e_\zeta)\x_\zeta\Bigr) = \sum_{\zeta \in \Sigma_1} \~m(\e_{\basetm\zeta}) \x_{\basetm\zeta} = \sum_{\zeta \in \~\Sigma_1} \~m(\e_{\zeta}) \x_{\zeta}. \]
This implies the inclusion $\Psi(\ssJ_{\Sigma}) \subseteq \ssJ_{\~\Sigma}$, and the lemma follows.
\end{proof}

Let $\e_\rho^*$ be the element of $\~M$ which takes value one on $\e_\rho$, and which vanishes on $N$. Note that for any ray $\basetm\zeta$, we have $\e_{\basetm\zeta} = (\e_\zeta,0) + f(\e_\zeta) \e_\rho$. It follows that $\e_\rho^*(\e_{\basetm\zeta}) = f(\e_\zeta)$.

In the Chow ring $A^\bul(\~\Sigma)$, we have the following equation
\[ \sum_{\zeta \in \~\Sigma_1} \e_\rho^*(\e_\zeta) x_\zeta =0. \]
This implies the following result.

\begin{lemma} \label{lem:fundamental_equality_chow_ring_tropical_modification}
In the Chow ring $A^\bul(\~\Sigma)$, we have
\[ x_\rho = -\sum_{\zeta \in \Sigma_1} f(\e_\zeta) x_{\basetm \zeta}. \]
\end{lemma}

\begin{proof}[Proof of Proposition~\ref{prop:chow_ring_tropical_modification}]
Applying Lemma~\ref{lem:map_psi}, we get a well-defined morphism of ring $\Psi$ from $A^\bul(\Sigma)$ to $A^\bul(\~\Sigma)$. It remains to prove the surjectivity. By the definition of $\Psi$, we just need to find a preimage for $x_\rho \in A^\bul(\~\Sigma)$. This follows from Lemma~\ref{lem:fundamental_equality_chow_ring_tropical_modification}. In the case $\Delta=\div(f)$ is trivial, the proof is similar.

The statement for Minkowski weights then follows by duality using Theorem \ref{thm:duality_A_MW}.
\end{proof}

\subsection{Div-faithfulness and stability of Chow rings under tropical modifications} \label{subsec:trop_mod_div-faithful}

We follow the preceding notations and denote by $\~\Sigma$ the tropical modification of $\Sigma$ along the divisor $\Delta$ given by the meromorphic function $f$. We denote by $\~\omega$ the orientation of $\~\Sigma$. The following theorem essentially implies the invariance of the Chow ring under tropical modifications under the assumption that the underlying tropical fan $\Sigma$ is div-faithful.

\begin{thm}[Stability of the Chow ring under tropical modifications] \label{thm:invariant_chow_tropical_modification}
Notations as above, let $\~\Sigma$ be the tropical modification of the tropical fan $\Sigma$ along $\Delta$ with respect to $f$. Assume furthermore that $\Sigma$ is div-faithful. Then, we get an isomorphism
\[ A^\bul(\~\Sigma,\Q) \simeq A^\bul(\Sigma,\Q) \]
between the Chow rings with rational coefficients, and an isomorphism
\[\MW_{\bul}(\~\Sigma) \simeq \MW_{\bul}(\Sigma).\]
between the Minkowski weights with integral coefficients. Moreover, if either $\Sigma$ is saturated, or $A^\bul(\Sigma)$ has no torsion, then we get an isomorphism between the Chow rings with integral coefficients
\[ A^\bul(\~\Sigma) \simeq A^\bul(\Sigma). \]
The isomorphisms are all compatible with the degree maps.
\end{thm}

The map $\MW_{\bul}(\Sigma, \Z) \to \MW_{\bul}(\~\Sigma, \Z)$ is given as follows. Let $\alpha \colon \Sigma_k \to \Z$ be an element of $\MW_k(\Sigma, \Z)$ and denote by $\Theta$ the subfan of $\Sigma$ defined by the support of $\alpha$. The pair $(\Theta, \alpha)$ is a tropical fan of dimension $k$ and the restriction of $f$ to $\Theta$ defines a meromorphic function $h = f\rest{\Theta}$ in $\mer(\Theta)$. By div-faithful property of $\Sigma$, the tropical modification $\~\Theta$ of $\Theta$ along the divisor $\div(h)$ is a Minkowski weight $\~\alpha$ in $\MW_{k}(\~\Sigma, \Z)$. The map $\MW_k(\Sigma, \Z) \to \MW_{k}(\~\Sigma, \Z)$ sends $\alpha$ to $\~\alpha$. The composition of this map and the one given in Proposition~\ref{prop:chow_ring_tropical_modification} is identity. Note that without the div-faithful assumption, the tropical modification $\~\Theta$ is not necessary a subfan of $\~\Sigma$. Indeed, the statement of the theorem does not necessarily hold in the absence of this assumption, see Example~\ref{ex:nontrivial_degenerate_tropical_modification}.

\smallskip
Our Example~\ref{ex:necessity_invariant_chow_tropical_modification} shows that, dealing with Chow rings with integral coefficients, the assumption made in the theorem on saturation or torsion freeness is needed.

\smallskip
We prove the statement in the theorem that concerns Chow rings with integral coefficients. The proof of the statements for Chow rings with rational coefficient can be obtained using a similar argument. The isomorphism for Minkowski weights with rational coefficients is then obtained by duality, Theorem \ref{thm:duality_A_MW}. The fact that the isomorphism restricts to Minkowski weights with integral coefficients stems from the fact that maps in both directions preserve integral weights.

If $\Delta$ is empty, since $\Sigma$ is div-faithful and saturated at $\conezero$ (either by assumption or by Proposition \ref{prop:characterization_of_saturation}), $f$ is globally integral linear and we have the isomorphism $\~\Sigma \simeq \Sigma$ from which the result follows. Hence, in what follows, we assume $\Delta\neq\emptyset$. In this case, we denote by $\rho = \R_{\geq 0}\etm$ the special ray of $\~\Sigma$. We already know there exists a surjective map $\Psi\colon A^\bul(\Sigma) \twoheadrightarrow A^\bul(\~\Sigma)$ given in Proposition \ref{prop:chow_ring_tropical_modification}. We now construct a surjective map $\Phi\colon A^\bul(\~\Sigma) \twoheadrightarrow A^\bul(\Sigma)$ and show that $\Phi\circ\Psi=\id$, from which we get the result.

\smallskip
We first define a surjective map $\Phi$ on the level of polynomial rings using Lemma~\ref{lem:fundamental_equality_chow_ring_tropical_modification}. This is the map
\[ \Phi \colon \Z[\x_\zeta \mid \zeta\in\~\Sigma_1] \to \Z[\x_{\zeta} \mid \zeta\in \Sigma_1] \]
defined on the level of generators as follows. Take the linear form $\e_\rho^*$ on $\~N$ that takes value one on the primitive vector $\e_\rho = \uptm\e$ of the special ray $\rho$ and which vanishes on $N$. Each ray $\zeta$ of $\~\Sigma$ is either of the form $\basetm\zeta$ for $\zeta \in \Sigma_1$ or is equal to $\rho$. We set
\[ \Phi(\x_{\basetm \zeta}) \coloneqq \x_\zeta, \qquad \textrm{and} \qquad \Phi(\x_\rho) \coloneqq - \sum_{\zeta \in\Sigma_1} f(\e_{\zeta})\x_\zeta. \]

As before, $\ssI_{\Sigma}$, $\ssJ_{\Sigma}$, $\ssI_{\~\Sigma}$ and $\ssJ_{\~\Sigma}$ are the ideals appearing in the definition of the Chow rings of $\Sigma$ and $\~\Sigma$.

\begin{lemma} Notations as above, we have
\begin{enumerate}
\item $\Phi(\ssJ_{\~\Sigma}) \subseteq \ssJ_{\Sigma}$, and
\item if\/ $\Sigma$ is div-faithful, then we have $\Phi(\ssI_{\~\Sigma}) \subseteq \ssI_{\Sigma} +\ssJ_{\Sigma}$.
\end{enumerate}
\end{lemma}

\begin{proof}
To show that $\Phi(\ssJ_{\~\Sigma}) \subseteq \ssJ_{\Sigma}$, consider an integral linear form $l$ on $\~N$. Let $\~ m \coloneqq l - l(\e_\rho)\e_\rho^*$. We have $\~m(\e_\rho)=0$ and so $\~m$ gives an integral linear form on $N \simeq \rquot{\~N}{\Z \e_\rho}$. We denote this by~$m$. We have, using $\e_{\basetm\zeta} = (\e_{\zeta},0) + f(\e_{\basetm\zeta})\e_\rho$,
\begin{align*}
\Phi\Bigl(\sum_{\zeta\in\~\Sigma_1} l(\e_\zeta)\x_\zeta\Bigr) &= l(\e_\rho)\Phi(\x_{\rho}) + \sum_{\zeta\in \Sigma_1} l(\e_{\basetm \zeta})\Phi(\x_{\basetm\zeta})=-\sum_{\zeta\in \Sigma_1} l(\e_\rho)f(\e_{\zeta}) \x_{\zeta} + \sum_{\zeta\in \Sigma_1} l(\e_{\basetm \zeta})\x_{\zeta}\\
&=-\sum_{\zeta\in \Sigma_1} l(\e_\rho)f(\e_{\zeta}) \x_{\zeta} + \sum_{\zeta\in \Sigma_1} \bigl(\~m(\e_{\basetm\zeta})+ l(\e_\rho)f(\e_\zeta)\bigr)\x_{\zeta}=\sum_{\zeta\in \Sigma_1} m(\e_\zeta)\x_\zeta.
\end{align*}
This shows that $\Phi(\ssJ_{\~\Sigma}) \subseteq \ssJ_{\Sigma}$.

\smallskip
We now consider the image of $\ssI_{\~\Sigma}$. Consider a collection of distinct rays $\tilde\zeta_1, \dots, \tilde\zeta_k$ of $\~\Sigma$, $k\in \N$, and suppose they are not comparable so that we have $\x_{\tilde\zeta_1} \dots \x_{\tilde\zeta_k} \in \ssI_{\~\Sigma}$. Two cases can happen.
\begin{enumerate}
\item \label{enum:invariant_chow:different_rho} Either, $\tilde\zeta_1, \dots, \tilde\zeta_k$ are different from $\rho$.
\item \label{enum:invariant_chow:contains_rho} Or, one of the rays, say $\tilde\zeta_1$, is equal to $\rho$.
\end{enumerate}

Consider the case \ref{enum:invariant_chow:different_rho}. We have $\tilde\zeta_j = \basetm{{\zeta_j}}$ for $j=1, \dots, k$ and rays $\zeta_j$ in $\Sigma$. Moreover, the rays $\zeta_1, \dots, \zeta_k$ do not form a cone in $\Sigma$ and we get $\Phi(\x_{\tilde\zeta_1} \dots \x_{\tilde\zeta_k}) = \x_{\zeta_1} \dots \x_{\zeta_k} \in \ssI_{\Sigma}$, as desired.

\smallskip
Consider now the case \ref{enum:invariant_chow:contains_rho}. Let as in the previous case, $\tilde\zeta_j = \basetm{{\zeta_j}}$ for $j=2, \dots, k$ and rays $\zeta_j$ in $\Sigma$. At this point, two cases can happen:

Either, these rays do not form a cone in $\Sigma$ in which case we get
\[ \Phi(\x_{\tilde\zeta_1}\x_{\tilde\zeta_2} \dots \x_{\tilde\zeta_k}) = \Phi(\x_{\rho}) \x_{\zeta_2} \dots \x_{\zeta_k} \in \ssI_{\Sigma}.
\]

Or, $\zeta_2, \dots, \zeta_k$ form a cone $\tau$ in $\Sigma$. Note that $\tau \not\in\Delta$, as otherwise, $\basetm\tau$ and $\rho$ would be comparable which would contradict our assumption. This implies that the divisor $\div(f)^\tau$ induced by $\div(f)$ on $\Sigma^\tau$ is trivial. In particular, the induced holomorphic function $f^\tau$ on $\Sigma^\tau$ is obtained by the restriction of $f-m$ on $\Sigma^\tau$, for an element $m\in M$ which verifies $f\rest\tau = m\rest\tau$. It verifies moreover $\div(f^\tau)=0$.

Assume first that $\Sigma$ is saturated, we deal with the case $A^\bul(\Sigma)$ has no torsion later. Since $\Sigma$ is div-faithful and saturated at $\tau$, this implies that $f^\tau$ coincides with an element $m^\tau \in M^\tau$. Altogether, this means $f$ coincides with the restriction of the linear function $l\coloneqq m+\pi_{\conezero\subfaceeq\tau}^*(m^\tau) \in M$ on all faces $\sigma$ of $\Sigma$ with $\sigma \supfaceeq \tau$, where $\pi_{\conezero\subfaceeq\tau}\colon N_\R \to N_\R^\tau$ denotes the projection. Set $\hat f = f-l$, and note that $\hat f$ is zero on every ray comparable with $\tau$.

We now observe that
\[ \Phi(\x_\rho) = -\sum_{\zeta\in\Sigma_1} f(\e_{\zeta}) \x_\zeta = -\sum_{\zeta\in\Sigma_1} \hat f(\e_{\zeta}) \x_\zeta - \sum_{\zeta\in\Sigma_1} l(\e_\zeta) \x_\zeta. \]
Using the notation $\x_\tau = \x_{\zeta_2} \cdots \x_{\zeta_k}$, we get
\[ \Phi(\x_{\tilde\zeta_1}\x_{\tilde\zeta_2} \dots \x_{\tilde\zeta_k}) = \Phi(\x_{\rho})\x_{\tau}
  = -\sum_{\zeta\in\Sigma_1} \hat f(\e_\zeta)\x_\zeta\x_{\tau} - \Bigl(\sum_{\zeta\in\Sigma_1} l(\e_\zeta) \x_\zeta\Bigr) \x_{\tau}. \]
Since $\hat f(\e_\zeta)$ is trivial if $\zeta$ is comparable with $\tau$, the first sum is in $\ssI_\Sigma$. Since $l \in M$, the second term is in $\ssJ_\Sigma$. This finishes the proof in the case $\Sigma$ is saturated.

If $\Sigma$ is not saturated but has a Chow ring without torsion, then we work as above but tensoring everything with $\Q$ because $m^\tau$ might not be integral. This way we get that $\Phi(\x_{\tilde\zeta_1}\cdots\x_{\tilde\zeta_k})$ belongs to $A^k(\Sigma) \cap (I_\Sigma\otimes\Q + J_\Sigma\otimes\Q)$. Hence, some multiple of this element belongs to $I_\Sigma + J_\Sigma$. Since the Chow ring is torsion-free, the element itself must be zero in $A^k(\Sigma)$.

In any case, we conclude that $\phi(\x_{\rho_1}\x_{\rho_2} \cdots \x_{\rho_k}) \in \ssI_\Sigma + \ssJ_\Sigma$, as desired. This finishes the proof of the lemma.
\end{proof}

\begin{proof}[Proof of Theorem~\ref{thm:invariant_chow_tropical_modification}]
Applying the above lemma, we obtain $\Phi(\ssJ_{\~\Sigma}) \subseteq \ssJ_\Sigma$ and $\Phi(\ssI_{\~\Sigma}) \subseteq \ssI_\Sigma +\ssJ_\Sigma$. This induces a map
\[ \Phi \colon A^\bul(\~\Sigma) \to A^\bul(\Sigma). \]
This map is surjective, and verifies $\Phi\circ\Psi=\id$. The theorem follows.
\end{proof}

\section{$\mT$-stability} \label{sec:operations}

In the previous sections, we defined three types of operations on tropical fans: products, tropical modifications and stellar subdivisions/assemblies. In this section, we study a notion of \emph{$\mT$-stability} for subclasses of a class $\Csh$ of tropical fans, using these operations. The idea behind the definition is that if a property $P$ holds for some basic tropical fans and if, in addition, this property happens to be preserved by the above three operations, a wide collection of tropical fans verify $P$. This happens in practice for various geometric properties that will be discussed later in the paper.

This leads to the definition of quasilinear fans: roughly speaking, a tropical fan is quasilinear if it can be obtained from the most basic tropical fans, the point, with arbitrary orientation, and the line by using only the above three operations.

As in the previous sections, $\Sigma$ will be a tropical fan of pure dimension $d$ in $N_\R \simeq \R^n$ for some natural number $n$, and its orientation will be denoted by $\omega!_\Sigma$.

\subsection{$\mT$-stability} \label{sec:t-stability}

By our convention from Section~\ref{sec:prel}, we work with fans modulo isomorphisms. This allows to talk about the set of isomorphism classes of rational fans.

\begin{defi} \label{def:t-stability}
Let $\Csh$ be a class of tropical fans (or more precisely, a set of isomorphism classes of tropical fans). A subclass $\Sh \subseteq \Csh$ is called \emph{$\mT$-stable in $\Csh$}, or simply \emph{$\mT$-stable} if $\Csh$ is the class of all tropical fans, if it verifies the following properties:
\begin{itemize}
\item (Stability under products) If $\Sigma, \Sigma' \in \Sh$ and the product $\Sigma \times \Sigma'$ belongs to $\Csh$, then $\Sigma \times \Sigma' \in \Sh$.

\item (Stability under tropical modifications along a divisor in the subclass) If $\Sigma \in \Sh$ and if $f$ is a meromorphic function on $\Sigma$ such that either $\div(f)$ is trivial or $\div(f)\in \Sh$, and if in addition $\tropmod{f}{\Sigma} \in \Csh$, then $\tropmod{f}{\Sigma} \in \Sh$.

\item (Stability under blow-ups and blow-downs with center in the subclass) If $\Sigma \in \Csh$, for any cone $\sigma \in \Sigma$ and any ray $\rho$ in the relative interior of $\sigma$ which verify $\Sigma^\sigma \in \Sh$ and $\Sigma_{(\rho)} \in \Csh$, we have
\[ \Sigma \in \Sh \Longleftrightarrow \Sigma_{(\rho)} \in \Sh. \qedhere\]
\end{itemize}
\end{defi}

We would like to make a comment on terminology. One way of looking at $\mT$-stability is to think of $\mT$ as a set of operations on tropical fans, namely products, tropical modifications, blow-ups and blow-downs. A precise meaning of this idea is given in Appendix \ref{sec:multimagmoid}, where we introduce multi-magmoids and provide an equivalent reformulation of $\mT$-stability.

The class $\Csh$ itself is clearly $\mT$-stable in $\Csh$. Moreover, an intersection of $\mT$-stable subclasses of $\Csh$ remains $\mT$-stable in $\Csh$. This implies that the subclass $\gst{\Bsh}{\Csh} \subseteq \Csh$ in the following definition exists.

\begin{defi}[$\mT$-stable subclass generated by a base set] \label{def:stability}
Let $\Csh$ be a class of tropical fans. Let $\Bsh$ be a subset of $\Csh$ that we call the \emph{base set}. The $\mT$-stable subclass of $\Csh$ generated by $\Bsh$, denoted by $\gst{\Bsh}{\Csh}$ is defined as the smallest $\mT$-stable subclass of $\Csh$ which contains $\Bsh$. If $\Csh$ is the class of all tropical fans, we just write $\gst{\Bsh}{}$.
\end{defi}

The main examples of classes $\Csh$ of tropical fans which are of interest to us are \emph{all}, \resp \emph{simplicial}, \resp \emph{unimodular}, \resp \emph{unimodular quasi-projective} tropical fans. We call them the \emph{standard classes} of tropical fans. These classes only constrain blow-ups and blow-downs since they are all closed by products and by tropical modifications (see Theorem \ref{thm:standard_properties}).

\smallskip
An important example of the base set is the set $\Bsho$ consisting of two simple fans: the fan $\conezero$ with an arbitrary orientation, and the unique complete fan $\Lambda$ in $\R$, with lattice $N=\Z$, with three cones $\conezero, \R_{\geq 0}$, and $\R_{\leq 0}$ and with the orientation constant equal to one. The four \emph{basic sets} that we will consider are $\Bsho, \Bsho_+, \Bsho_1$ and $\Bsho_{\pm1}$, which are the set of all, \resp effective, \resp reduced, \resp unitary, elements of $\Bsho$. We see later in Section~\ref{subsec:Bergman_fans_stable} that $\gst{\Bsho_1}{}$, and so $\gst{\Bsho}{}$, contains many interesting fans.

\begin{defi}[Quasilinear tropical fans]
Let $\Csh$ be a class of tropical fans and $\Bsh$ a subset of $\Csh$. A tropical fan $\Sigma$ is called \emph{$\mT$-generated by $\Bsh$ in $\Csh$} if $\Sigma$ belongs to $\gst{\Bsh}{\Csh}$.

If $\Csh$ is the class of all tropical fans and $\Bsh=\Bsho$, we simply say that $\Sigma$ is \emph{quasilinear}.
\end{defi}

The following result is a consequence of Theorem~\ref{thm:stability_factorization} proved later in this section.

\begin{prop}
Let $\Sigma$ be a quasilinear tropical fan in $N_\R$. Then any tropical fan $\Sigma'$ in $N_\R$ with the same support and the same integral structure is quasilinear.
\end{prop}

\begin{defi}[$\mT$-stable properties]
Let $\Csh$ be a class of tropical fans and let $P$ be a predicate on elements of $\Csh$. Then, $P$ is called \emph{$\mT$-stable in $\Csh$} if the subclass of $\Csh$ consisting of those tropical fans that verify $P$ is $\mT$-stable in $\Csh$. If $\Csh$ is the class of all tropical fans, we simply say $P$ is \emph{$\mT$-stable}.
\end{defi}

\subsection{Star-stability}

A class $\Ssh$ of tropical fans is called \emph{star-stable} if for any $\Sigma \in \Ssh$ and any $\sigma \in \Sigma$, the star fan $\Sigma^\sigma$ also belongs to $\Ssh$. For instance, the four standard classes as well as the four basic sets are star-stable (\cf Theorem \ref{thm:standard_properties}). A predicate $P$ on fans is called \emph{star-stable} or \emph{local} if $P(\Sigma)$ implies $P(\Sigma^\sigma)$ for any $\sigma \in \Sigma$.

There is a natural way to construct a star-stable predicate from an arbitrary one. If $P$ is any predicate on tropical fans, then we denote by $\Ploc P$ the predicate
\[ \Ploc P(\Sigma)\colon\quad \forall\sigma\in\Sigma,\ P(\Sigma^\sigma). \]
For instance, local and $\Q$-local irreducibility, being div-faithful, and being principal and $\Q$-principal are all defined in this way relative to the corresponding property holding at $\conezero$. As we will show later in Lemma \ref{lem:stability_meta_lemma}, viewing them as such will allow to considerably simplify the proofs of their $\mT$-stability.

\subsection{Strong $\mT$-stability and intersection property}

\begin{defi}[Strong $\mT$-stability]
Let $\Csh$ be a class of tropical fans, and let $\Sh \subseteq \Csh$ be a subclass. We say that $\Sh$ is \emph{strongly $\mT$-stable in $\Csh$} if $\Sh$ is $\mT$-stable in $\Csh$ and in addition, we have
\begin{itemize}
  \item for any pair of fans $\Sigma$, $\Sigma'$ in $\Csh$, $\Sigma\times\Sigma'\in\Sh$ implies $\Sigma$ and $\Sigma'$ are in $\Sh$,
  \item for any $\Sigma \in \Csh$ and any meromorphic function $f$ on $\Sigma$ such that either $\div(f)$ is trivial or $\div(f)\in\Csh$, if $\tropmod{f}{\Sigma} \in \Sh$, then both $\Sigma$ and $\div(f)$ are in $\Sh$,
  \item if $\Sigma \in \Csh$, for any cone $\sigma \in \Sigma$ and any ray $\rho$ in the relative interior of $\sigma$ which verify $\Sigma^\sigma \in \Csh$ and $\Sigma_{(\rho)} \in \Csh$, we have
  \[ \Sigma \in \Sh \Longleftrightarrow \Sigma_{(\rho)} \in \Sh \Longrightarrow \Sigma^\sigma \in \Sh. \qedhere\]
\end{itemize}
\end{defi}

Examples are given by the following proposition.

\begin{prop}
Let\/ $\Csh$ be any of the standard classes. The subclass of reduced, \resp effective, \resp unitary tropical fans of\/ $\Csh$ is strongly $\mT$-stable in $\Csh$.
\end{prop}

The proof is direct and omitted.

We have the following intersection property, see Proposition~\ref{prop:basic_properties_multi-magmoid} in Section~\ref{sec:multimagmoid} for the proof.

\begin{prop}\label{prop:intersection-property}
Let $\Csh$ be a class of tropical fans. If $\Sh\subseteq \Csh$ is strongly $\mT$-stable in $\Csh$, then for any base set $\Bsh$, we have
\[ \gst{\Bsh}{\Csh}\cap\Sh = \gst{\Bsh \cap \Sh}{\Csh}. \]
\end{prop}

Being effective, \resp reduced, \resp unitary, is strongly $\mT$-stable. Hence, a direct consequence of Proposition \ref{prop:intersection-property} is that $\gst{\Bsho_+}{\Csh}$, \resp $\gst{\Bsho_1}{\Csh}$, \resp $\gst{\Bsho_{\pm1}}{\Csh}$, is precisely the set of effective, \resp reduced, \resp unitary, quasilinear tropical fans.

\subsection{Properties of standard classes}

The following theorem summarizes several nice properties enjoyed by the four standard classes of tropical fans that we introduced in Section~\ref{sec:t-stability}.

\begin{thm} \label{thm:standard_properties}
Let $\Csh$ be one of the four standard classes of tropical fans. Then $\Csh$ verifies the following properties.
\begin{itemize}
\item \textnormal{(Closedness under products)} If\/ $\Sigma$ and $\Sigma'$ belong to $\Csh$, then we have $\Sigma \times \Sigma' \in \Csh$.
\item \textnormal{(Closedness under containment)} If\/ $\Sigma$ is in $\Csh$, any subfan $\Delta$ of\/ $\Sigma$ endowed with arbitrary orientation is in $Csh$.
\item \textnormal{(Closedness under tropical modifications)} If\/ $\Sigma \in \Csh$ and if $f$ is a meromorphic function on $\Sigma$, then $\tropmod{f}{\Sigma} \in \Csh$.
\item \textnormal{(Containment of the basic fans)} Fans of $\Bsho$ are in $\Csh$.
\item \textnormal{(Star-stability)} The class $\Csh$ is star-stable.
\item \textnormal{(Existence of unimodular quasi-projective subdivisions)} Any fan in $\Csh$ has a subdivision in $\Csh$ which is unimodular and quasi-projective.
\item \textnormal{(Weak factorization)} Let $\Sigma$ and $\Sigma'$ be two fans in $\Csh$ with the same support and compatible orientations. Then, there exists a sequence of fans $\Sigma=\Sigma^0, \Sigma^1, \Sigma^2, \dots, \Sigma^k=\Sigma'$ in $\Csh$ such that for any $i\leq k-1$, $\Sigma^{i+1}$ is obtained from $\Sigma^i$ by performing a blow-up or a blow-down.
\end{itemize}
\end{thm}

\begin{proof}
We only sketch the proof here.
\begin{itemize}[leftmargin=0pt, itemindent=1em]
\item (Closedness under products) For the case of all tropical fans this is stated by Proposition~\ref{prop:orientable_stable}. Moreover, simpliciality and unimodularity are preserved by taking products. To see that the product of two quasi-projective fans $\Sigma$ and $\Sigma'$ is quasi-projective, consider the projection maps $\pi\colon \Sigma \times \Sigma' \to \Sigma$ and $\pi'\colon \Sigma \times \Sigma' \to \Sigma'$. For two strictly convex conewise linear functions $f$ and $f'$ on $\Sigma$ and $\Sigma'$, respectively, the sum $\pi^*(f)+\pi'^*(f')$ is a strictly convex conewise linear function on $\Sigma \times \Sigma'$.

\item (Closedness under containment) For the case of all, \resp simplicial, \resp unimodular, tropical fans, this is trivial. For quasi-projectivity, notice that the restriction of a strictly convex function to a subfan remains strictly convex.

\item (Closedness under tropical modifications) For the case of all tropical fans, this is Proposition~\ref{prop:balancing_trop_modif}. The new cones of the form $\uptm\delta$ with $\delta$ in the divisor are simplicial, \resp unimodular, provided $\delta$ is simplicial, \resp unimodular. Hence, being simplicial and being unimodular are preserved by tropical modifications. For quasi-projectivity, let $g$ be a strictly convex conewise linear function on a fan $\Sigma$, and let $\prtm\colon \~\Sigma \to \Sigma$ be the projection associated to the tropical modification $\~\Sigma = \tropmod{f}{\Sigma}$. We show the conewise linear function $\prtm^*(g)$ is strictly convex on $\~\Sigma$. Let $\delta \in \Delta$, we prove that $\prtm^*(g)$ is strictly convex around $\basetm\delta$. Similar arguments show strict convexity around other faces of $\~\Sigma$. Since $g$ is strictly convex around $\delta$, there exists an element $m \in M_\R$ such that $g-m$ is zero on $\delta$ and is strictly positive on all the incident rays $\rho\sim\delta$ in $\Sigma$. Let $l$ be a linear form in $\~M_\R$ which takes values one on $\etm$ and vanishes on $\basetm\delta$. For a small enough positive real number $\varepsilon$, $\prtm^*(g)-\prtm^*(m)+\varepsilon l$ is zero on $\basetm\delta$, takes value $\varepsilon>0$ on $\etm$, and is strictly positive on rays $\basetm\rho$ for $\rho\sim\delta$. This proves that $\prtm^*(g)$ is strictly convex around $\basetm\delta$.

\item (Containment of the basic fans) This is trivial.

\item (Star-stability) For the case of all tropical fans, star-stability follows from Proposition \ref{prop:orientable_stable}. Being simplicial and being unimodular are both local properties. For quasi-projectivity, a strictly convex conewise linear function on a fan induces strictly convex conewise linear functions on the star fans around its faces.

\item (Existence of unimodular quasi-projective subdivisions) This is a well-known fact. We refer to Section 4 of \cite{AP-tht} for more details.

\item (Weak factorization) This last property is far from being trivial. For the class of simplicial tropical fans and the class of unimodular tropical fans, this follows from Theorem A of~\cite{Wlo97}, proved independently by Morelli~\cite{Mor96} and expanded by Abramovich-Matsuki-Rashid, see~\cite{AMR}. For unimodular quasi-projective tropical fans, this can be obtained from relevant parts of~\cites{Wlo97, Mor96, AKMW} as discussed and generalized by Abramovich and Temkin in~\cite{AT19}*{Section 3}. Note that we are requiring the orientations to be consistent, so the statement here is only about the fan structure. \qedhere
\end{itemize}
\end{proof}

\subsection{$\mT$-stability, support and factorization}

The following shows that in some cases of interest, $\mT$-stability is only a property of the support with a fixed ambient lattice.

\begin{thm} \label{thm:stability_factorization}
Let $\Csh$ be a standard class. Let $\Ssh$ be subclass of\/ $\Csh$ which is both star-stable and $\mT$-stable in $\Csh$. Then, a tropical fan $\Sigma$ of\/ $\Csh$ is in $\Ssh$ if and only if any tropical fan of\/ $\Csh$ with the same support $\supp{\Sigma}$ considered with the same lattice is in $\Ssh$.
\end{thm}

\begin{remark}
We note a subtle point here worth emphasizing. In the statement of the theorem, we remember the ambient lattice. This is weaker than being a property of the support. For instance, if $\Sigma$ is not saturated in $N_\R$ and if $\Sigma'$ is the same fan considered with a different lattice $N' = (\Sigma\cap N)\otimes\Z$, then the above theorem does not imply that $\Sigma\in\Ssh$ if and only if $\Sigma'\in\Ssh$, even though $\Sigma\cap N = \Sigma'\cap N'$. We refer to Example~\ref{ex:F1_vs_N} related to the $\mT$-stability statement \ref{thm:examples_stable:PD_Chow_ring} of Section \ref{sec:stable_examples} in the study of the Chow ring, proved in Section \ref{sec:PD_stable}.

For those properties that depend only on $\Sigma\cap N$ in the sense of Remark \ref{rem:independent_exterior_lattice}, this remark is irrelevant.
\end{remark}

\begin{proof}
Let $\Sigma\in\Ssh$ and let $\Sigma'$ be another fan of $\Csh$ with the same support and the same lattice. By the weak factorization property of Theorem \ref{thm:standard_properties}, there exists a sequence of fans $\Sigma = \Sigma^0, \Sigma^1, \dots, \Sigma^{k-1}, \Sigma^k = \Sigma'$ all belonging to $\Csh$ such that $\Sigma^{i+1}$ is obtained from $\Sigma^i$ by performing a blow-up or a blow-down.

We prove that $\Sigma^1 \in \Ssh$. If $\Sigma^1$ is obtained from $\Sigma$ by blowing up a face $\eta \in \Sigma$, then we have $\Sigma^\eta\in\Ssh$ by star-stability. Since $\Ssh$ is $\mT$-stable in $\Csh$, the blow-up along $\eta$ of $\Sigma$ will be in $\Ssh$, and we get $\Sigma^1 \in \Ssh$. Now, if $\Sigma$ is obtained from $\Sigma^1$ by blowing up along a ray $\rho$ which is in the relative interior of a face $\eta'\in\Sigma^1$, then, we note that $(\Sigma^1)^{\eta'} = \Sigma^{\tau\vee\rho}$ where $\tau$ is any face of codimension one in $\eta'$. Once again, $(\Sigma^1)^{\eta'}$ belongs to $\Ssh$. Since $\Ssh$ is $\mT$-stable in $\Csh$, it is closed under blow-down along $\rho$, and we get $\Sigma^1 \in \Ssh$.

Proceeding this way step by step, we obtain that $\Sigma' = \Sigma^k \in \Ssh$.
\end{proof}

\begin{remark} \label{rem:stability_factorization}
In the proof we only used that $\Ssh$ is star-stable and closed under blow-ups and blow-downs along faces whose star fans belong to $\Ssh$.
\end{remark}

\subsection{Examples of $\mT$-stable geometric properties} \label{sec:stable_examples}

Here is a list of geometric properties which are $\mT$-stable in relevant classes of tropical fans. The first two points are easy to verify. We will prove the other ones later in this article.
\begin{enumerate}
\item Connectedness through codimension one is $\mT$-stable.
\item Being effective, \resp reduced, \resp unitary, is $\mT$-stable.
\item \label{thm:examples_stable:normal} $\Q$-normality is $\mT$-stable.
\item \label{thm:examples_stable:irreducible} Irreducibility and $\Q$-local irreducibility are $\mT$-stable.
\item \label{thm:examples_stable:div-faithful} Div-faithfulness is $\mT$-stable.
\item \label{thm:examples_stable:principal} Being $\Q$-principal is $\mT$-stable in the class of $\Q$-locally irreducible fans.
\item \label{thm:examples_stable:PD_Chow_ring} Poincaré duality with $\Q$-coefficients, \resp $\Z$-coefficients, for the Chow ring in the sense of Theorem \ref{thm:PD_stable} is $\mT$-stable in the class of div-faithful unimodular fans.
\item \label{thm:examples_stable:HR_Chow_ring} Being Chow-Kähler in the sense of Theorem \ref{thm:chow-KP_stable} is $\mT$-stable in the class of effective quasi-projective unimodular fans.
\end{enumerate}

Since the statements listed above are verified by elements of $\Bsho_1$, we infer that they are true in $\gst{\Bsho_1}{}$. That is, reduced quasilinear tropical fans verify all these properties. More generally, the different statements hold for quasilinear tropical fans in the adequate class.

\subsection{A tool to prove $\mT$-stability} \label{subsec:stability_meta_lemma}

Checking all the axioms of $\mT$-stability for a class $\Csh$ can be somewhat tedious in general. The following lemma is helpful in practice to simplify the verification of these different points.

\begin{lemma}[$\mT$-stability meta lemma] \label{lem:stability_meta_lemma}
Let $\Csh$ be a star-stable and\/ $\mT$-stable subclass of one of the four standard classes such that $\Bsho_1 \subseteq \Csh$. Let $P$ be a predicate on elements of\/ $\Csh$. Assume that the elements of $\Bsho_1$ verify $P$.

Let $\Sigma$ be an arbitrary tropical fan in $\Csh$ such that for any face $\sigma \neq \conezero$ in $\Sigma$, the star fan $\Sigma^\sigma$ verifies $\Ploc P$, and such that at least one of the following points is verified.
\begin{enumerate}[label=(ML\arabic*)]
\item \label{lem:ml1} $\Sigma$ is the product of two unimodular quasi-projective tropical fans in $\Csh$ verifying $\Ploc P$.
\item \label{lem:ml2} $\Sigma$ is the tropical modification of a unimodular quasi-projective tropical fan $\Sigma'\in\Csh$ verifying $\Ploc P$ with respect to some meromorphic function $f$ on $\Sigma'$ such that either $\div(f)$ is trivial, or the tropical fan $\Delta$ associated to $\div(f)$ is an element of $\Csh$ which is unimodular, quasi-projective, and verifies $\Ploc P$.
\item \label{lem:ml3} $\Sigma$ is the blow-up along some ray of a tropical fan in $\Csh$ verifying $\Ploc P$.
\item \label{lem:ml4} $\Sigma$ is the blow-down along some ray of a tropical fan in $\Csh$ verifying $\Ploc P$.
\end{enumerate}
If for any tropical fan $\Sigma$ as above, the property $P$ is verified, then $\Ploc P$ is $\mT$-stable in $\Csh$.

\smallskip
Moreover, if $P$ is a predicate only depending on the support of the tropical fan in the sense of Theorem~\ref{thm:stability_factorization}, we can restrict ourselves to tropical fans $\Sigma$ verifying one of the two first points.
\end{lemma}

Using this lemma, we can reduce the proof of $\mT$-stability statements to verifying that in each of the four cases enumerated above, the fan $\Sigma$ verifies $P$. For instance, we have to prove that if $\Sigma$ is a product of two fans verifying $P^\star$, then $\Sigma$ verifies $P$. To do so, we can assume without loss of generality that the two factors (and thus $\Sigma$ itself) are unimodular and quasi-projective, and that every proper star fan of $\Sigma$ verifies $P$.

\smallskip
Before proving the lemma, we state some consequences.

\begin{prop} \label{prop:loc_P_stable}
Notation as in Lemma \ref{lem:stability_meta_lemma}, if a property $P$ is $\mT$-stable in $\Csh$, then the property $\Ploc{P}$ is $\mT$-stable in $\Csh$.
\end{prop}

\begin{proof}
The property $\Ploc{P}$ implies $P$. The four conditions of the lemma are verified since $P$ is $\mT$-stable in $\Csh$. Therefore, $\Ploc{P}$ is $\mT$-stable in $\Csh$.
\end{proof}

\begin{cor} \label{cor:generated_in_cap}
Notation as in Lemma \ref{lem:stability_meta_lemma}, let $\Bsh \subseteq \Csh$ be a star-stable subclass of\/ $\Csh$ containing $\Bsho_1$. Then $\gst{\Bsh}{\Csh}$ is star-stable. In particular it verifies Theorem~\ref{thm:stability_factorization}. Moreover, $\gst{\Bsh}{\Csh} = \gst{\Bsh}{} \cap \Csh$.
\end{cor}

\begin{proof}
For $\Sigma \in \Csh$, let $P(\Sigma)$ be the predicate \enquote{$\Sigma$ belongs to $\gst{\Bsh}{\Csh}$}. By the previous proposition, we infer that $\Ploc{P}$ is $\mT$-stable in $\Csh$. Since elements of $\Bsh$ verify $\Ploc{P}$, we deduce that every element of $\gst{\Bsh}{\Csh}$ verifies~$\Ploc{P}$, \ie, $\gst{\Bsh}{\Csh}$ is star-stable. In particular, Theorem~\ref{thm:stability_factorization} applies.

\smallskip
For the last statement, consider the class $\Ssh$ of all tropical fans having the same support as an element of $\gst{\Bsh}{\Csh}$ in the sense of Theorem \ref{thm:stability_factorization}. From the star-stability of $\gst{\Bsh}{\Csh}$ we deduce the star-stability of $\Ssh$. We can apply Lemma \ref{lem:stability_meta_lemma} to the predicate \enquote{$\Sigma$ belongs to $\Ssh$} inside the class of all tropical fans. The first two points \ref{lem:ml1}-\ref{lem:ml2} follow directly from the $\mT$-stability of $\gst{\Bsh}{\Csh}$, with the help of Theorem~\ref{thm:stability_factorization}. And the last two points \ref{lem:ml3}-\ref{lem:ml4} are direct by the definition of $\Ssh$. Hence, $\Ssh$ is $\mT$-stable.

Notice that $\Bsh \subseteq \Ssh$. Moreover, by the first part of the corollary applied to the whole class of tropical fans, $\gst{\Bsh}{}$ verifies Theorem \ref{thm:stability_factorization}. Since $\gst{\Bsh}{\Csh} \subseteq \gst{\Bsh}{}$, we deduce that $\Ssh \subseteq \gst{\Bsh}{}$. By $\mT$-stability, we get that $\Ssh = \gst{\Bsh}{}$. Now we know that any element of $\gst{\Bsh}{} \cap \Csh$ has the same support than an element of $\gst{\Bsh}{\Csh}$. Applying once more Theorem \ref{thm:stability_factorization}, we get that $\gst{\Bsh}{} \cap \Csh \subseteq \gst{\Bsh}{\Csh}$. The other inclusion is clear and the result follows.
\end{proof}

\begin{remark}
From Proposition~\ref{prop:intersection-property} and Corollary~\ref{cor:generated_in_cap} we deduce for instance the following useful results. Let $\Bsh$ be $\Bsho$, \resp $\Bsho_{+}$, \resp $\Bsho_{\pm1}$, \resp $\Bsho_{1}$, and let $\Csh$ be as in Lemma \ref{lem:stability_meta_lemma}. Assume $\Bsh \subseteq \Csh$. Then $\gst{\Bsh}{\Csh}$ is the class of all, \resp effective, \resp unitary, \resp reduced, quasilinear tropical fans that are in $\Csh$.
\end{remark}

\smallskip
The rest of this section is devoted to the proof of Lemma~\ref{lem:stability_meta_lemma}. First, notice that the last statement in Lemma~\ref{lem:stability_meta_lemma} is clear since blow-ups and blow-downs do not change the support of the fan.

Let $\Csh_0$ be the standard class in which $\Csh$ is $\mT$-stable. Extending trivially the predicate $P$ in $\Csh_0$, \ie, such that $P(\Sigma)$ is false if $\Sigma\in\Csh_0\setminus\Csh$, we can assume without loss of generality that $\Csh = \Csh_0$ is a standard class.

Let $\Ssh \subseteq \Csh$ be the class of fans verifying $\Ploc P$. Since $\Ploc P$ is star-stable, so is $\Ssh$. For any integer $n$, we use the notation $\Csh^{<n}$, \resp $\Csh^{\leq n}$, to denote the subset of $\Csh$ of fans of dimension less than $n$, \resp at most $n$. We define $\Ssh^{<n}$ and $\Ssh^{\leq n}$ similarly.

We will prove that under the assumption of the lemma, $\Ssh$ is $\mT$-stable in $\Csh$. We prove this by induction on $n$. Assume that, for some integer $n$, $\Ssh^{< n}$ is $\mT$-stable in $\Csh^{<n}$. For $n=1$, this is a consequence of \ref{lem:ml1} in the lemma. We prove that $\Ssh^{\leq n}$ is $\mT$-stable in $\Csh^{\leq n}$.

\subsubsection{Closedness under blow-ups and blow-downs} We verify that $\Ssh^{\leq n}$ is closed under blow-ups and blow-downs. Let $\Sigma$ be a tropical fan in $\Csh^{\leq n}$. Let $\eta \in \Sigma$ be a face and let $\rho$ be a ray in the relative interior of $\eta$. Assume that $\Sigma_{(\rho)} \in \Csh^{\leq n}$. We need to show the equivalence
\[ \Sigma \in \Ssh^{\leq n} \Longleftrightarrow \Sigma_{(\rho)} \in \Ssh^{\leq n}. \]

We first prove if $\Sigma_{(\rho)} \in \Ssh^{\leq n}$, then $\Sigma \in \Ssh^{\leq n}$. For this, we compare the star fans of $\Sigma$ and $\Sigma_{(\rho)}$ as follows. Consider a face $\sigma$ of $\Sigma$ different from $\conezero$. There are three cases.

\begin{itemize}
\item First, assume that $\sigma$ is not comparable with $\eta$. In this case, the two star fans $\Sigma_{(\rho)}^\sigma$ and $\Sigma^\sigma$ are identical. Since $\Ssh$ is star-stable, the first one is in $\Ssh$ by assumption. Hence $\Sigma^\sigma\in\Ssh$.

\item Second, we assume that $\sigma \supfaceeq \eta$. In this case, we have $\Sigma^\sigma \simeq \Sigma_{(\rho)}^{\tau\vee\rho}$ where $\tau$ is any face of codimension one in $\sigma$ such that $\tau\wedge\eta \ssubface \eta$. Once again $\Sigma^\sigma \in \Ssh$.

\item Finally, assume that $\sigma$ and $\eta$ are comparable but $\sigma \not\supfaceeq \eta$. Denote by $\eta^\sigma$, \resp $\rho^\sigma$, the cone corresponding to $\eta$, \resp to $\rho$, in $\Sigma^\sigma$. Then $\rho^\sigma$ is a ray in the relative interior $\eta^\sigma$. The star fan $\Sigma_{(\rho)}^\sigma$ is then naturally isomorphic to the blow-up star fan $\bigl(\Sigma^\sigma\bigr)_{(\rho^\sigma)}$. By assumption, $\Sigma_{(\rho)}^\sigma \in \Ssh^{<n}$. Moreover $\bigl(\Sigma^\sigma\bigr)^{\eta^\sigma} = \Sigma^{\sigma\vee\eta}$ which is in $\Ssh^{<n}$ by the second point above. Hence, applying the closedness by blow-down of $\Ssh^{<n}$, we deduce that $\Sigma^\sigma$ is also in $\Ssh^{<n} \subseteq \Ssh$.
\end{itemize}

In any case, $\Sigma^\sigma \in \Ssh$ for any face $\sigma\neq\conezero$. One can apply the assumption of the lemma to deduce that $\Sigma$ verifies $P$. Hence, $\Sigma$ verifies $\Ploc P$, thus $\Sigma \in \Ssh^{\leq n}$. This proves the direction $\Leftarrow$.

\smallskip
To prove the direction $\Rightarrow$, assume $\Sigma \in \Ssh^{\leq n}$. Take a cone $\nu$ in $\Sigma_{(\rho)}$ different from $\conezero$. Apart from the faces which already appeared in the above case analysis for which the reversed argument applies, it remains to consider those faces $\nu$ with $\rho \subfaceeq \nu$. Denote by $\nu-\rho \in \Sigma$ the face of $\nu$ of codimension one which does not contain $\rho$. Set $\sigma=(\nu-\rho)\vee\eta$. Then we get
\[ \supp{\Sigma_{(\rho)}^\nu} \simeq \supp{\Sigma^\sigma} \times \R^k, \]
where $k = \dims{\sigma}-\dims{\nu} = \dims\eta-\dims{\nu\cap\eta} \in \{0, \dots, n-1\}$. Since $\Ssh^{<n}$ is star-stable and closed under blow-ups and blow-downs, we can apply Remark \ref{rem:stability_factorization} and Theorem \ref{thm:stability_factorization}: $\Sigma_{(\rho)}^\nu$ is in $\Ssh^{<n}$ if and only if there exists a fan with the same support in $\Ssh^{<n}$. This is the case. Indeed, $\Sigma^\sigma \in \Ssh^{<n}$. Moreover, $\Lambda$ and $\{\conezero\}$ belong to $\Ssh$ by assumption. Since $\Ssh^{<n}$ is $\mT$-stable in $\Csh^{<n}$, we deduce that $\Sigma^\sigma \times \Lambda^k \in \Ssh^{<n}$. This last fan has the same support as $\Sigma_{(\rho)}^\nu$. Hence, we infer that $\Sigma_{(\rho)}^\nu \in \Ssh^{<n} \subseteq \Ssh$.

We have proved that $\Sigma_{(\rho)}^\nu \in \Ssh$ for any nontrivial face $\nu \in \Sigma_{(\rho)}$. As before, we apply the assumption of the lemma to deduce that $\Sigma_{(\rho)} \in \Ssh^{\leq n}$. Hence, $\Ssh^{\leq n}$ is closed under blow-ups and blow-downs. In particular, we can apply Theorem \ref{thm:stability_factorization} in $\Ssh^{\leq n}$ for the rest of this proof.

\subsubsection{Closedness under products} We whish to prove that $\Ssh^{\leq n}$ is closed under products which remain inside $\Csh^{\leq n}$. Let $\~\Sigma\in\Csh^{\leq n}$ be the product of two fans in $\Ssh^{\leq n}$. We prove that $\~\Sigma\in\Ssh^{\leq n}$.

Denote this two factors by $\~\Sigma^1$ and $\~\Sigma^2$. Let $\Sigma^1$, \resp $\Sigma^2$, be a unimodular quasi-projective subdivision in $\Csh$ of $\Sigma^1$, \resp of $\Sigma^2$, which exists by Theorem \ref{thm:standard_properties}. Then, clearly $\Sigma^1$ has the same support as $\~\Sigma^1$, and by Theorem \ref{thm:stability_factorization}, we get $\Sigma^1\in\Ssh^{\leq n}$. In the same way, we obtain $\Sigma^2\in\Ssh^{\leq n}$. Set $\Sigma=\Sigma^1\times\Sigma^2$. Let $\sigma^1\times\sigma^2$ be a nontrivial face of $\Sigma$. Then
\[ \Sigma^{\sigma^1\times\sigma^2} \simeq (\Sigma^1)^{\sigma^1} \times (\Sigma^2)^{\sigma^2}. \]
By star-stability of $\Ssh$, both factors belong to $\Ssh^{<n}$. Hence, the $\mT$-stability in $\Ssh^{<n}$ implies that the product belongs to $\Ssh^{<n}$.

Therefore, for any nontrivial face $\sigma$ of $\Sigma$, we get $\Sigma^\sigma\in\Ssh$. Applying the assumption of the lemma, we deduce that $\Sigma$ verifies $P$ and thus $\Sigma \in \Ssh^{\leq n}$. Since $\Sigma$ and $\~\Sigma$ have the same support, we can apply Theorem \ref{thm:stability_factorization} to deduce that $\~\Sigma\in\Ssh^{\leq n}$ as well. Thus, $\Ssh^{\leq n} \subset \Csh^{\leq n}$ is closed by products.

\subsubsection{Closedness under tropical modifications} Let $\Sigma'$ be a fan in $\Ssh^{\leq n}$. Let $f$ be a meromorphic function on $\Sigma'$. Set $\~\Sigma' = \tropmod{f}{\Sigma'}$. Let $\Delta'$ be the tropical fan associated to $\div(f)$ and assume that $\Delta'$ is in $\Ssh^{\leq n}$ (by convention, we assume in this proof that $\emptyset \in \Ssh^{\leq n}$). We wish to prove that $\~\Sigma'$ is in $\Ssh^{\leq n}$.

Let $\Sigma$ be a unimodular quasi-projective subdivision of $\Sigma'$. Then $f$ is a meromorphic function on $\Sigma$. Moreover, the tropical fan $\Delta$ associated to $\div(f)$, taken in $\Sigma$, is a unimodular quasi-projective subdivision of $\Delta'$. As for the case of the product, Theorem \ref{thm:stability_factorization} implies that both $\Sigma$ and $\Delta$ are in $\Ssh^{\leq n}$. Set $\~\Sigma = \tropmod{f}{\Sigma}$. A face $\~\sigma$ of dimension $k>0$ in $\~\Sigma$ is of two kinds, either it is equal to $\uptm\delta$ for $\delta \in \Delta_{k-1}$ or it coincides with $\basetm\sigma$ for $\sigma \in \Sigma_{k}$. By Proposition \ref{prop:star_fans_tropical_modifications}, in the first case, the star fan $\~\Sigma^{\uptm\delta}$ is isomorphic to $\Delta^\delta$ and so belongs to $\Ssh$. So we can now assume that $\~\sigma = \basetm\sigma$ for a cone $\sigma \in \Sigma$. Then, $\~\Sigma^{\basetm\sigma}$ is the tropical modification of $\Sigma^\sigma$ along $\Delta^\sigma$ with respect to the function $f^\sigma$. By convention here we set $\Delta^\sigma = \emptyset$ if $\sigma\not\in\Delta$. Note that $\Sigma^\sigma$ and $\Delta^\sigma$ are in $\Ssh^{<n}$. Since $\Ssh^{<n}$ is $\mT$-stable in $\Csh^{<n}$, $\tropmod{f^\sigma}{\Sigma^\sigma} \in \Ssh^{<n}$. We infer again that $\~\Sigma^{\basetm\sigma} \in \Ssh$, as desired.

At this point we have verified that for any nontrivial cone $\~\sigma$ in $\~\Sigma$, the star fan $\~\Sigma^{\~\sigma}$ is in $\Ssh$. Using the assumption of the lemma, we deduce that $\~\Sigma$ verifies $P$ and so $\~\Sigma \in \Ssh^{\leq n}$. By Theorem \ref{thm:stability_factorization}, we deduce that $\~\Sigma' \in \Ssh^{\leq n}$. Therefore, $\Ssh^{\leq n} \subseteq \Csh^{\leq n}$ is closed under tropical modifications.

\subsubsection{End of the proof} We have proved that $\Ssh^{\leq n}$ is $\mT$-stable in $\Csh^{\leq n}$. By induction, we deduce that $\Ssh$ is $\mT$-stable in $\Csh$, \ie, $\Ploc P$ is $\mT$-stable in $\Csh$. \qed

\section{Bergman fans} \label{subsec:Bergman_fans_stable}

This section is devoted to recalling basic definitions and properties regarding matroids and their Bergman fans. Bergman fans of matroids are quasi-projective since they can be realized as subfans of the Bergman fan of a free matroid that is obtained by the barycentric subdivision of the fan of a projective space, see Remark~\ref{rem:permutohedron}.

It is easy to see that complete unimodular fans are quasilinear (see Section~\ref{sec:proof-matroid-stability}). In this section we prove the following generalization of this statement.

\begin{thm} \label{thm:bergman_stable}
The Bergman fan and the augmented Bergman fan of a matroid are quasilinear. More generally, any generalized Bergman fan is quasilinear.
\end{thm}

\subsection{Matroids} \label{sec:matroids}

We start by briefly recalling basic definitions involving matroids and refer to relevant part of \cite{Oxl06} for more details. A matroid can be defined in different equivalent ways, for example by specifying what is called its \emph{collection of independent sets}, or its \emph{collection of bases}, or its \emph{collection of flats}, or its \emph{collection of circuits}, or still by giving its \emph{rank function}. The data of any of these collections determine all the others.

\begin{defi}[Matroid: definition with respect to the family of independent sets] \label{defi:matroid}
A \emph{matroid} $\Ma$ is a pair $(E, \ind)$ consisting of a finite set $E$ called the \emph{ground set} and a collection $\ind$ of subsets of $E$ called the \emph{family of independent sets of\/ $\Ma$} which verifies the following axiomatic properties:
\begin{enumerate}
\item The empty set is an independent set: $\emptyset\in\ind$.
\item (Hereditary property) $\ind$ is closed under inclusion: if $J\subseteq I$ and $I\in\ind$, then $J\in\ind$.
\item (Augmentation property) for two elements $I, J\in\ind$, if $\Card J<\Card I$, then one can find an element $i$ in $I\setminus J$ such that $J\cup\{i\}\in\ind$. \label{defi:matroid:augmentation_property} \qedhere
\end{enumerate}
\end{defi}

\begin{remark}
An example of a such a pair $=(E, \ind)$ is given by a collection of vectors $v_1, \dots, v_m$ in a finite dimensional vector space $H$ over a field $\k$. The ground set $E$ is $[m]$ and the collection $\ind$ of independent sets consists of all subsets $I \subseteq [m]$ verifying that the corresponding vectors $v_i$ for $i\in I$ are linearly independent. A matroid $\Ma$ of this form is called \emph{representable} (over $\k$). Nelson shows in~\cite{Nel18} that \emph{almost any matroid is non-representable over any field}.
\end{remark}

To a given matroid $\Ma=(E,\ind)$ we can associate the so-called \emph{rank function} $\rkm\colon 2^E \to \Z_{\geq 0}$ which is defined as follows. For a subset $A \subseteq E$, the \emph{rank of $A$} is defined as the maximum of $\Card I$ over all independent sets $I$ which are included in $A$. The integer $\rkm(E)$ is called the \emph{rank of\/ $\Ma$}. The rank function satisfies the \emph{submodularity property}
\[ \forall A, A' \subseteq E, \qquad \rkm(A\cap A') + \rkm(A \cup A') \leq \rkm(A) + \rkm(A'). \]

A \emph{basis} of $\Ma$ by definition is a maximal independent set. The collection of bases of $\Ma$ is denoted by $\bases(\Ma)$. A \emph{circuit} of $\Ma$ is a minimal \emph{dependent} set. A set is called dependent if it is not independent. The collection of circuits of $\Ma$ is denoted by $\crct(\Ma)$.

The \emph{closure} $\cl(A)$ of a subset $A$ in $E$ is defined as
\[ \cl(A)\coloneqq\{i \in E \textrm{ such that } \rkm(A \cup \{i\})=\rkm(A)\}. \]
A \emph{flat} of $\Ma$ is a subset $F \subseteq E$ with $\cl(F) = F$. Flats are also equally called \emph{closed sets} and the collection of flats of the matroid $\Ma$ is denoted by $\Cl(\Ma)$. A flat of $\Ma$ is called \emph{proper} if it is different from $E$. The set of nonempty proper flats of $\Ma$ is denoted $\prCl(\Ma)$.

Let $i$ be an element of $E$. The matroid $\Ma\del\{i\}$, \resp $\Ma\contr\{i\}$, is the matroid obtained by deleting $i$, \resp contracting $i$, in $\Ma$. Both these matroids have ground set $E \setminus \{i\}$, and their set of flats is characterized by
\begin{gather*}
  \Cl(\Ma\del\{i\}) = \left\{ F \st F \in \Cl(\Ma)\text{ or } F\cup\{i\} \in \Cl(\Ma) \right\}, \quad
  \Cl(\Ma\contr\{i\}) = \left\{ F \st F \cup\{i\} \in \Cl(\Ma) \right\}.
\end{gather*}

An element $i$ with $\rkm(\{i\})=0$ is called a \emph{loop}. Two elements $i$ and $i'$ of $E$ are called \emph{parallel} if $\rkm(\{i, i'\}) = \rkm(\{i\}) = \rkm(\{i'\}) = 1$. A matroid is \emph{simple} if it neither contains loops nor parallel elements.

An element $i\in E$ is called a \emph{coloop} if $\rkm(E\setminus\{i\})=\rkm(E)-1$. This is equivalent to $E\setminus\{i\}$ being a flat.

\subsection{Bergman fans}

Let $\Ma$ be a simple matroid of rank $r$ on a ground set $E$. The Bergman fan of $\Ma$ denoted by $\ssSigma_{\Ma}$ is defined as follows.

For a subset $A \subseteq E$, we denote by $\e!_A$ the sum $\sum_{i\in A} \e_i$. Here $\{\e_i\}_{i\in E}$ is the standard basis of $\R^E$. Let $N = \rquot{\Z^E}{\Z \e!_E}$ and denote by $M$ the dual of $N$. By an abuse of the notation, we denote by the same notation $\e!_A$ the projection in $N_\R$ of $\sum_{i\in A}\e_i$. Note in particular that $\e!_E =0$.

The Bergman fan $\ssSigma_\Ma$ of $\Ma$ is the rational fan in $N_\R$ of dimension $r -1$ defined as follows. First, a \emph{flag of nonempty proper flats} $\Fl$ of $\Ma$ is a collection
\[ \Fl\colon \quad \emptyset \neq F_1 \subsetneq F_2 \subsetneq \dots \subsetneq F_{\ell} \neq E \]
consisting of flats $F_1, \dots, F_\ell$ of $\Ma$. The number $\ell$ is called the \emph{length} of $\Fl$.

To such a flag $\Fl$, we associate the rational cone $\sssigma_\Fl \subseteq N_\R$ generated by the vectors $\e!_{F_1}, \e!_{F_2}, \dots, \e!_{F_\ell}$, that is,
\[ \sssigma_\Fl \coloneqq \Bigl\{ \lambda_1 \e!_{F_1} + \dots + \lambda_\ell \e!_{F_\ell} \Bigst \lambda_1, \dots, \lambda_\ell \geq0\Bigr\}. \]
The dimension of $\sssigma_\Fl$ is equal to the length of $\Fl$.

The \emph{Bergman fan of\/ $\Ma$} is the fan consisting of all the cones $\sssigma_\Fl$, $\Fl$ a flag of nonempty proper flats of $\Ma$, \ie,
\[ \ssSigma_\Ma \coloneqq \Bigl\{\, \sssigma_\Fl \Bigst \Fl \textrm{ flag of nonempty proper flats of }\Ma\,\Bigr\}. \]
The fan $\ssSigma_\Ma$ (endowed with weight function taking value 1 on each facet) is a tropical fan of pure dimension $\rkm(\Ma)-1$. A \emph{generalized Bergman fan} is any fan isomorphic to a fan with support $\supp{\ssSigma_\Ma}$ (considered with the same lattice) for some matroid $\Ma$.

\begin{example}
For non-negative integers $r$ and $n$ with $r+1\leq n$, the uniform matroid $U_{r+1, n}$ has ground set $E=[n]$ with collection of independent sets $\ind$ consisting of all subsets of $E$ of size bounded by $r+1$. The Bergman fan of $U_{3,4}$ is depicted in Figure~\ref{fig:Bergman_fan_U34}.
\end{example}

\begin{remark} \label{rem:permutohedron}
The support of $\ssSigma_{U_{r+1, r+1}}$ is the full space $\R^r \simeq \rquot{\R^{r+1}}{\R \e!_{E}}$. The fan itself is the normal fan of the permutohedron. Moreover, for any matroid $\Ma$ on $r+1$ elements, $\ssSigma_\Ma$ is a subfan of $\ssSigma_{U_{r+1,r+1}}$.
\end{remark}

\begin{remark} \label{rem:bergman_parallel_elements}
If a matroid $\Ma$ has parallel elements but no loops, then we can still define the Bergman fan $\ssSigma_\Ma$ of $\Ma$. This fan is isomorphic to the Bergman fan of the matroid $\~\Ma$ obtained from $\Ma$ by deleting all but one element in each set of parallel elements.
\end{remark}

Let $\Ma$ be a simple matroid on the ground set $E=[m]$ and denote by $\crct(\Ma)$ the circuits of $\Ma$. The support of $\ssSigma_\Ma$ is described by a theorem of Ardila-Klivans~\cite{AK06} as follows. Let $\~\Sigma$ be the set of points $x=(x_i)_{i \in E}\in \R^E$ such that for every circuit $C\in\crct(\Ma)$, the minimum of $x_i$ for $i\in C$ is achieved at least twice. Note that if $x\in\~\Sigma$, then $x+\lambda\e!_E\in\~\Sigma$ for all $\lambda\in\R$.

\begin{thm}[Ardila-Klivans~\cite{AK06}] \label{thm:support_bergman}
The support of\/ $\ssSigma_\Ma$ coincides with the projection of\/ $\~\Sigma$ in $\rquot{\R^E\!}{\R\e!_E}$.
\end{thm}

The augmented Bergman fan of $\Ma$ defined in \cite{BHMPW} is the fan $\sshatSigma_\Ma$ which lives in the space $\R^E$, and which has the following description.

An independent set $I$ in $\Ma$ is called \emph{compatible} with a flag of (possibly empty) proper flats $\Fl\colon \, F_1 \subsetneq \dots \subsetneq F_l$ if $I \subseteq F_1$. In this case, the pair $(I, \Fl)$ is called compatible. To any compatible pair $(I, \Fl)$ of an independent set $I$ and a flag of proper flats $\Fl$ in $\Ma$ we associate the unimodular cone $\sssigma_{I,\Fl}$ in $\R^E$ defined by
\[\sssigma_{I,\Fl} \coloneqq \sum_{i\in I} \R_{\geq 0}\e!_i + \sum_{F\in \Fl} \R_{\geq 0}(-\e!_{E\setminus F}). \]
The augmented Bergman fan $\sshatSigma_\Ma$ is the collection of all the cones $\sssigma_{I,\Fl}$, $(I, \Fl)$ compatible pair in $\Ma$. The augmented Bergman fan $\sshatSigma_\Ma$ is a generalized Bergman fan. Indeed, the augmented Bergman fan $\sshatSigma_\Ma$ has the same support as the Bergman fan $\ssSigma_{\~\Ma}$, where $\~\Ma$ is obtained as a free coextension of $\Ma$ by a single element. This means $\~\Ma$ has ground set $E\sqcup\{e\}$ and bases of $\~\Ma$ are either $B \sqcup\{e\}$ for $B$ a basis of $\Ma$, or $A$ where $A$ is a subset of $E$ of size $\rk(\Ma)+1$ and of maximal rank in $\Ma$. Equivalently, the dual $\~\Ma^\dual$ is a free extension of $\Ma^\dual$ (with the same rank).

\subsection{Products of generalized Bergman fans}

\begin{prop} \label{prop:productBergman}
The product of two generalized Bergman fans is again generalized Bergman.
\end{prop}

More precisely, for two matroids $\Ma$ and $\Ma'$, we show the relation
\[ \supp{\ssSigma_{\Ma}} \times \supp{\ssSigma_{\Ma'}} \simeq \supp{\ssSigma_{\Ma\vee\Ma'}} \]
where $\Ma\vee \Ma'$ is any \emph{parallel connection} of $\Ma$ and $\Ma'$. We give the definition below and refer to~\cite{Oxl06}*{Chapter 7} for more details. See also Remark \ref{rem:parallel_connection} to get an intuition.

A \emph{pointed matroid} is a matroid $\Ma$ on a ground set $E$ with a choice of a distinguished element $\distel = \distel!_\Ma$ in $E$. Let now $\Ma$ and $\Ma'$ be two pointed matroids on the ground sets $E$ and $E'$ with distinguished elements $\distel!_\Ma \in E$ and $\distel!_{\Ma'} \in E'$, respectively. The parallel connection of $\Ma$ and $\Ma'$ denoted by $\Ma \vee \Ma'$ is by definition the pointed matroid on the ground set $E \vee E' = \rquot{E \sqcup E'}{(\distel!_\Ma=\distel!_{\Ma'})}$, the wedge sum of the two pointed sets $E$ and $E'$, with distinguished element $\distel = \distel!_\Ma = \distel!_{\Ma'}$, and with the following collection of bases:
\begin{align*}
\bases( \Ma \vee \Ma') =& \Bigl\{ \,B \cup B' \st B \in \bases(\Ma), B' \in \bases(\Ma') \,\textrm{ with }\,\distel!_\Ma\in B\, \textrm{ and }\,  \distel!_{\Ma'}\in B'\Bigr \} \\
& \cup  \Bigl\{ \,B \cup B' \setminus\{\distel\} \st B \in \bases (\Ma), B' \in \bases(\Ma') \, \textrm{ with }\, \distel!_\Ma\in B \, \textrm{ and }\, \distel!_{\Ma'}\notin B'\Bigr \}\\
& \cup  \Bigl\{ \,B \cup B' \setminus\{\distel\}\st B \in \bases(\Ma), B' \in \bases(\Ma')\, \textrm{ with }\, \distel!_\Ma \notin B\, \textrm{ and }\,  \distel!_{\Ma'}\in B'\Bigr \}.
\end{align*}
The circuits of $\Ma \vee \Ma'$ are given by
\begin{align*}
\crct(\Ma\vee \Ma') = &\ \crct(\Ma) \cup \crct(\Ma') \\
& \ \cup \Bigl\{ (C \cup C')\setminus \{\distel\} \st C \in\crct(\Ma) \textrm{ with } \distel!_\Ma \in C \myand C'\in\crct(\Ma')\textrm{ with } \distel!_{\Ma'} \in C' \Bigr\}.
\end{align*}

A parallel connection of two matroids $\Ma$ and $\Ma'$ is a wedge sum of the form $(\Ma, i) \vee (\Ma', i')$ for some choices of elements $i\in \Ma$ and $i'\in \Ma'$ that turn them into pointed matroids.

\begin{remark} \label{rem:parallel_connection}
The analogous operation for graphs consists in gluing two different graphs along distinguished oriented edges as illustrated below. One can check that the relations between the bases, \resp the circuits, of three matroids involved in a parallel connection mimics the relations between the spanning trees, \resp the circuits, of the three graphs.

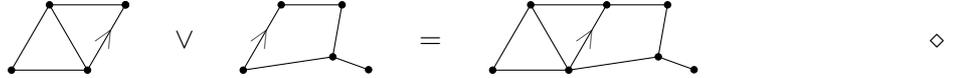
\[ \tikz[vcenter]{\draw (0,0)[dot] --++(60:1)[dot] --+(0:1)[dot] ++(0,0) --++(-60:1)[dot] --(0,0); \draw (0:1)--++(60:1)[dot] node[midway, sloped] {>};}
\quad \vee \quad
\tikz[vcenter]\draw (0,0)[dot]--++(60:1)[dot] node[midway, sloped] {>} --++(0:.8)[dot] --++(80:-.7)[dot] --+(-20:.5)[dot] ++(0,0)--(0,0);
\quad = \quad
\tikz[vcenter]{\draw (0,0)[dot] --++ (60:1)[dot] --+(0:1) ++(0,0) --++(-60:1)[dot] --(0,0); \draw (0:1)--++(60:1)[dot] node[midway, sloped] {>} --++(0:.8)[dot] --++(80:-.7)[dot] --+(-20:.5)[dot] +(0,0) --(0:1);} \qedhere \]
\end{remark}

\begin{proof}[Proof of Proposition~\ref{prop:productBergman}]
Let $\Ma$ and $\Ma'$ be two matroids on ground sets $E=[m]$ and $E'=[m']$, respectively. We show that the support of $\ssSigma_{\Ma}\times\ssSigma_{\Ma'}$ is isomorphic, by an integral linear isomorphism on the ambient spaces, to $\ssSigma_{\Ma\vee \Ma'}$.

Consider the following maps
\[ \R^{E\vee E'} \xrightarrow{\ \phi\times \phi'\ } \R^{E} \times \R^{E'} \xrightarrow{\ \pi \times \pi' \ } \rquot{\R^{E}\!}{\R\e!_E} \times \rquot{\R^{E'}\!}{\R\e!_{E'}}. \]
Here, $\pi$, $\pi'$, $\phi$ and $\phi'$ are the natural projections.

We follow the notations introduced in Theorem~\ref{thm:support_bergman} and consider the subsets $\sshatSigma_{\Ma \vee \Ma'}$, $\sshatSigma_{\Ma}$ and $\sshatSigma_{\Ma'}$ of $\R^{E\vee E'}$, $\R^E$ and $\R^{E'}$, respectively.

Since $\crct(\Ma\vee \Ma')$ contains $\crct(\Ma)$ and $\crct(\Ma')$, $\phi\times \phi'$ restricts to a map $\hat\phi\times \hat \phi'$ from $\sshatSigma_{\Ma\vee \Ma'}$ to $\sshatSigma_{\Ma} \times \sshatSigma_{\Ma'}$. Denote by $\psi\colon \sshatSigma_{\Ma\vee\Ma'} \to \supp{\ssSigma_{\Ma}}\times \supp{\ssSigma_{\Ma'}}$ the composition of this map with the projection $\pi\times\pi'$. Clearly, $\psi$ is a linear map and $\ker(\psi)=\R\e_E$. It remains to prove that $\psi$ is surjective. This can be checked directly by applying Theorem~\ref{thm:support_bergman}, using the description given above of the circuits in $\Ma \vee \Ma'$.
\end{proof}

\subsection{Proof of Theorem~\ref{thm:bergman_stable}} \label{sec:proof-matroid-stability}

Let $\Ma$ be a simple matroid with Bergman fan $\ssSigma_\Ma$. Denote by $r+1$ the rank of $\Ma$ with $r$ a non-negative integer. We need to show that any unimodular fan $\Sigma$ with $\supp{\Sigma} = \supp{\ssSigma_\Ma}$ is in $\gst{\Bsho}{}$. Applying Theorem~\ref{thm:stability_factorization}, it will be enough to produce one such fan $\Sigma \in \gst{\Bsho}{}$.

First, we observe that since $\Lambda \in \Bsho$, the product $\Lambda^n$ is in $\gst{\Bsho}{}$. This implies that complete fans are all quasilinear.

We infer that unimodular complete fans are quasilinear. We can therefore assume that $\Ma$ is not a free matroid. We now proceed by induction on the size of $\Ma$.

There exists an element $i$ in the ground set $E$ of $\Ma$ such that $\Ma \del \{i\}$ has the same rank as $\Ma$. By induction, and applying Remark \ref{rem:bergman_parallel_elements} (since $\Ma\contr\{i\}$ may have parallel elements), the Bergman fans $\ssSigma_{\Ma\contr\{i\}}$ of $\Ma\contr\{i\}$ and $\ssSigma_{\Ma \del \{i\}}$ of $\Ma \del \{i\}$ are both quasilinear. It follows from Section \ref{sec:stable_examples} that $\ssSigma_{\Ma \del \{i\}}$ is principal and div-faithful. Moreover, by Lemma~\ref{lem:contraction-divisor-deletion} below, $\ssSigma_{\Ma\contr\{i\}}$ is a divisor in $\ssSigma_{\Ma \del \{i\}}$. Therefore, $\ssSigma_{\Ma\contr\{i\}}$ is the divisor of a holomorphic function in $\mer(\ssSigma_{\Ma\del\{i\}})$, unique up to addition by an integral linear function on $\ssSigma_{\Ma \del \{i\}}$. The tropical modification of $\ssSigma_{\Ma \del \{i\}}$ along the divisor $\ssSigma_{\Ma\contr\{i\}}$ is a well-defined unimodular quasi-projective fan $\Sigma$. By Lemma~\ref{lem:shaw} below, the supports of $\Sigma$ and $\ssSigma_\Ma$ are the same. We have produced a quasilinear fan with support $\supp{\ssSigma_\Ma}$ and the theorem follows. \qed

\begin{lemma} \label{lem:contraction-divisor-deletion} Notations as above, $\ssSigma_{\Ma\contr\{i\}}$ is a divisor in $\ssSigma_{\Ma \del \{i\}}$.
\end{lemma}

\begin{proof}
We know that $\Cl(\Ma\contr\{i\}) \subset \Cl(\Ma\del\{i\})$. By the definition of the Bergman fans, we infer that $\ssSigma_{\Ma\contr\{i\}}$ is a subfan of $\ssSigma_{\Ma \del \{i\}}$. Since $\ssSigma_{\Ma\contr\{i\}}$ is tropical and has codimension one in $\ssSigma_{\Ma \del \{i\}}$, the result follows.
\end{proof}

Since $\Ma\del\{i\}$ is quasilinear by induction hypothesis, $\Cl(\Ma\contr\{i\})$ is the divisor of a holomorphic function that is unique up to addition of an integral linear function.

\begin{lemma}[Shaw~\cite{Sha13a}]\label{lem:shaw} Notations as above, let $\Sigma$ be the tropical modification of\/ $\ssSigma_{\Ma \del \{i\}}$ along the divisor $\ssSigma_{\Ma\contr\{i\}}$. We have $\supp\Sigma = \supp{\ssSigma_\Ma}$.
\end{lemma}

We provide a short proof of this result written in our framework.

\begin{proof}
We denote by $N$ the ambient lattice of $\ssSigma_{\Ma}$ and by $N'$ the ambient lattice of both $\ssSigma_{\Ma\del\{i\}}$ and $\ssSigma_{\Ma\contr\{i\}}$. We have a projection $\pi\colon N \to N'$. We choose an arbitrary section $\i\colon N' \to N$ and identify $N$ with $N' \times \Z$ via this section. Let $f$ be the meromorphic functions on $\Sigma_{\Ma\del\{i\}}$ such that $(\e!_F, f(\e!_F)) = \e!_{\cl_\Ma(F)}$ for any $F \in \prCl(\Ma\del\{i\})$ and where the closure operator is taken with respect to $\Ma$. Let $\Gamma!_f \coloneqq \id \times f$ be the graph of $f$ as in Section \ref{subsec:tropical_modification}.

To a flag $\Fl\colon\, F_1 \subsetneq \dots \subsetneq F_k$ of nonempty proper flats of $\Ma\del\{i\}$ corresponds the cone $\sssigma_\Fl \in \Sigma_{\Ma\del\{i\}}$. We denote by $\cl_\Ma(\Fl)$ the flag $\cl_\Ma(F_1) \subsetneq \dots \subsetneq \cl_\Ma(F_k)$ of nonempty proper flats of $\Ma$. From the definition of $f$, $\Gamma!_f(\sssigma_\Fl) = \sigma_{\cl_\Ma(\Fl)} \in \Sigma_\Ma$.

A complete flag of nonempty proper flats of $\Ma$ which is not of the form $\cl_\Ma(\Fl)$ for $\Fl$ a complete flag of $\Ma\del\{i\}$ is necessarily of the form
\[ F_1 \subsetneq \dots \subsetneq F_l \subsetneq F_l\cup \{i\} \subsetneq F_{l+1}\cup\{i\} \subsetneq \dots \subsetneq F_{r-1} \cup \{i\}. \]
Hence facets $\eta$ of $\Sigma_\Ma$ which are not in $\Gamma!_f(\Sigma_{\Ma\del\{i\}})$ are parallel to $\e_i$, \ie, $\e_i \in N_\eta$, and they live above the codimension one face $\Gamma!_f(\sssigma_{\Fl'})$ where $\Fl'$ is the flag $F_1 \subsetneq \cdots \subsetneq F_{r-1}$. One can prove that $F_k\cup\{i\}$ are all in $\prCl(\Ma)$ for $k \leq l$, therefore $\Fl'$ is a flag of nonempty proper flats of $\Ma\contr\{i\}$.

We prove that $\Sigma_\Ma$ and $\tropmod{f}{\Sigma_{\Ma\del\{i\}}}$ are equal as weighted fans, after performing some subdivisions if necessary. Subdividing these two fans if necessary, one can find a complete fan containing both of them as subfans. We can see the two fans as Minkowski weights in this complete fan. Their difference has support of dimension $r$ included in
\[ \bigcup_{\delta \in \Sigma_{\Ma\del\{i\},r-1}} \hspace{-.65cm}\Gamma!_f(\delta) + \R_{\geq 0} \e_i. \]
It is easy to see the only such Minkowski weight is zero.

It remains to prove that $\div(f) = \Sigma_{\Ma\contr\{i\}}$. The divisor of $f$ can be identified with $\Sigma_\Ma^{\rho_i}$, the star fan around the ray $\rho_i = \R_{\geq 0}\e_i$. This star fan is generated by the faces $\sssigma_\Fl$ with $\Fl$ a flag of nonempty proper flats of $\Ma$ of the form
\[ \{i\} \subsetneq F_1\cup\{i\} \subsetneq \cdots \subsetneq F_k\cup\{i\}. \]
The flags $F_1 \subsetneq \cdots \subsetneq F_k$ form precisely the set of flags of nonempty proper flats of $\Ma\contr\{i\}$.
\end{proof}

\begin{remark}
In Example \ref{ex:stable_not_Bergman}, we show that the property of \emph{being a generalized Bergman fan} is not $\mT$-stable. That example is a quasilinear tropical fan which is not a generalized Bergman fan. Given the results we are about to prove about the $\mT$-stability of several geometric properties, such as Poincaré duality and the property of being Chow-Kähler, this means that the tropical fans to which our theorems apply go beyond the setting of matroids and fans with the same support as their Bergman fans.
\end{remark}

\section{Basic $\mT$-stability results}

The aim of this section is to prove the $\mT$-stability of several geometric properties introduced in the previous sections. These properties will be crucial in the treatment of the Kähler package in the following sections.

\subsection{Normality is $\mT$-stable} \label{sec:normal_stable}

We first prove the following theorem.

\begin{thm} \label{thm:normal_stable}
Being $\Q$-normal is $\mT$-stable. Being normal is $\mT$-stable in the class of unitary tropical fans.
\end{thm}

\begin{proof}
We prove the first statement. The second statement follows from $\mT$-stability of the class of unitary tropical fans, and the fact that $\Q$-normality and normality are equivalent for unitary tropical fans by Proposition \ref{prop:trop_normal_support_local_product}.

Recall that a tropical fan $(\Sigma, \omega)$ is normal if and only if for any face $\eta$ of codimension one in $\Sigma$, $(\Sigma^\eta, \omega^\eta)$ is a generalized tropical line. In particular, both normality and $\Q$-normality are star-stable, and $\Q$-normality only depends on the support (see Proposition \ref{prop:trop_normal_support_local_product}).

Denote by $P$ the predicate for a tropical fan to be $\Q$-normal. Therefore, $P = \Ploc P$. Note that the elements of $\Bsho_1$ are $\Q$-normal. Let $\Sigma$ be a tropical fan verifying the condition of the $\mT$-stability meta lemma, Lemma \ref{lem:stability_meta_lemma}. In particular, if the dimension of $\Sigma$ is at least two, then for any face $\eta$ of codimension one in $\Sigma$, the star fan $\Sigma^\eta$ is $\Q$-normal. As a consequence $\Sigma$ verifies $P$.

It remains to treat the case where $\Sigma$ is of dimension one. In this particular situation, the only non-trivial case to deal with is the one where $(\Sigma, \omega)$ is the tropical modification of a generalized tropical line in $\R^k$ for some integer $k$. There are two cases depending on whether the tropical modification is degenerate or not.

\begin{itemize}
\item First, assume $\Sigma$ is a non-degenerate tropical modification of a generalized tropical line $(L, \omega)$ in $\R^k$ with $k+1$ rays. Rays in $\Sigma$ are of the form $\basetm\zeta$, $\zeta$ a ray of $L$, and the special ray $\rho = \R_{\geq 0}\etm$. Denote the primitive vectors of the rays of $\Sigma$ by $\e_0, \dots, \e_k$, and $\etm$. Consider the projection map $\prtm \colon \tropmod{f}{\Sigma} \to \Sigma$ associated to the tropical modification. We have $\prtm(\etm)=0$. Let $a_0, \dots, a_k, \uptm a$ be scalar coefficients such that we have $a_0\e_0 + \dots +a_k \e_k+ \uptm a\etm = 0$. Applying $\prtm$, we deduce that $a_0\prtm(\e_0) + \dots + a_k\prtm(\e_k)=0$. The $\Q$-normality of $(L, \omega)$ implies that $a_0 = \lambda \omega(\zeta_0), \dots, a_k = \lambda \omega(\zeta_k)$ for some rational number $\lambda$. Moreover, by the balancing condition in $\Sigma$, we get $\omega(\zeta_0) \e_0 + \dots + \omega(\zeta_k)\e_k = -\omega(\rho)\etm$. Thus, $(-\lambda \omega(\rho) + \uptm a)\etm = 0$. So we should have $\uptm a = \lambda \omega(\rho)$, which shows that $(a_0, \dots, a_k, \uptm a)$ is a rational multiple of the weight vector. This proves that $\Sigma$ is $\Q$-normal.

\item In the second case, $\Sigma$ is a degenerate tropical modification of a generalized tropical line. We can prove in the same way as above that $\Sigma$ is again $\Q$-normal in this case. We omit the details.
\end{itemize}

We have proved that $\Q$-normality fulfills the conditions of Lemma \ref{lem:stability_meta_lemma}, and so, it is $\mT$-stable.
\end{proof}

\subsection{$\Q$-local irreducibility is $\mT$-stable} \label{sec:local_irreducibility_stable}

From the $\mT$-stability of $\Q$-normality, we can deduce the following theorem.

\begin{thm} \label{thm:irreducible_stable}
$\Q$-local irreducibility is $\mT$-stable. Local irreducibility is $\mT$-stable in the class of unitary tropical fans.
\end{thm}

\begin{proof}
Again we only prove the first statement. Denote by $Q$ the property of being connected through codimension one. A tropical fan $\Sigma$ is $\Q$-locally irreducible if and only if it is $\Q$-normal and verifies $\Ploc Q$. This is the analogue with rational coefficients of Proposition \ref{thm:characterization_irreducible}. By statement (1) in Section~\ref{sec:stable_examples}, $\Ploc Q$ is $\mT$-stable. We conclude by $\mT$-stability of $\Q$-normality.
\end{proof}

\subsection{$\mT$-stability of the principality} \label{sec:stability_principality}

In this section we prove the following theorem.

\begin{thm}[$\mT$-stability of the principality for locally irreducible tropical fans] \label{thm:stability_principal}
The property for a tropical fan to be $\Q$-principal is $\mT$-stable in the class of\/ $\Q$-locally irreducible tropical fans. Being principal is $\mT$-stable in the class of locally irreducible and unitary tropical fans.
\end{thm}

Before proceeding with the proof, we would like to make a few comments.

First, we note that, as shown in Example~\ref{ex:q-principal}, $\mT$-stability might fail for principality (instead of $\Q$-principality, stated in the theorem). Moreover, ($\Q$-)local irreducibility assumption is needed in the theorem, as explained by Example \ref{ex:div-faithful_not_irreducible}.

Finally, unlike the properties of normality and local irreducibility, treated in Theorems \ref{thm:normal_stable} and \ref{thm:irreducible_stable}, the equivalence between principality and $\Q$-principality does not hold in the class of unitary tropical fan. Example \ref{ex:cube} provides a reduced tropical fan that is not principal, but only $\Q$-principal.

\smallskip
We now turn to the proof of the theorem~\ref{thm:stability_principal}. We will only prove the first statement, as the second follows by the same reasoning. We use Lemma \ref{lem:stability_meta_lemma} with $\Csh$ the class of $\Q$-locally irreducible tropical fans and with $P$ the predicate of being $\Q$-principal at $\conezero$. We have already proved that $\Csh$ verifies the condition required in the lemma. The fact that elements of $\Bsho_1$ are $\Q$-principal is trivial. We show that the four properties \ref{lem:ml1}-\ref{lem:ml2}-\ref{lem:ml3}-\ref{lem:ml4} stated in the lemma are verified.

\subsubsection{Closedness under products \ref{lem:ml1}}

Consider two $\Q$-locally irreducible $\Q$-principal tropical fans $(\Sigma,\omega!_\Sigma)$ and $(\Sigma', \omega!_{\Sigma'})$ of dimension $d$ and $d'$, respectively. We prove that $\Sigma\times\Sigma'$ is $\Q$-principal at $\conezero$. By Künneth formula for Minkowski weights, we get
\[ \MW_{d+d'-1}(\Sigma\times \Sigma') = \MW_d(\Sigma) \otimes \MW_{d'-1}(\Sigma')\ \oplus\ \MW_{d-1}(\Sigma) \otimes \MW_{d'}(\Sigma'). \]
Using the above decomposition, and arguing by symmetry, we only need to verify that divisors of the form $D \times C'$ for $D \in \MW_{d-1}(\Sigma)$ and $C'\in \MW_{d'}(\Sigma')$ are $\Q$-principal.

Since $\Sigma'$ is $\Q$-locally irreducible, we can suppose that $C'$ is equal to the element $[\Sigma']$ in $\MW_{d'}(\Sigma')$ given by $\Sigma'$. Since $\Sigma$ is $\Q$-principal, we have $aD = \div(f)$ for a conewise integral linear function $f$ on $\Sigma$ and a positive integer $a$. Denote by $\pi\colon \Sigma\times\Sigma' \to \Sigma$ the natural projection. We have $\div(\pi^*f) = \div(f)\times[\Sigma'] = aD\times C'$, from which the result follows.

\subsubsection{Closedness under tropical modifications \ref{lem:ml2}}

Let $(\Sigma, \omega!_\Sigma)$ be a $\Q$-principal $\Q$-locally irreducible tropical fan, and let $(\~\Sigma, \omega!_{\~\Sigma})$ be a tropical modification of $\Sigma$. By Lemma~\ref{lem:stability_meta_lemma}, we can assume without loss of generality that $\Sigma$ is unimodular.

By Proposition~\ref{prop:chow_ring_tropical_modification}, the application $A^1(\Sigma) \to A^1(\~\Sigma)$ is surjective and the application $\MW_{d-1}(\~\Sigma) \to \MW_{d-1}(\Sigma)$ is injective. From the explicit description of the map $\cycl$ given in Proposition~\ref{prop:map_cl}, we get the following commutative diagram
\begin{equation} \label{diag:tropical_modification_chow_minkowski}
\begin{tikzcd}
  A^1(\~\Sigma)  \rar{\cycl}                        & \MW_{d-1}({\~\Sigma}) \dar[hook] \\
  A^1(\Sigma)    \uar[twoheadrightarrow]\rar{\cycl} & \MW_{d-1} (\Sigma).
\end{tikzcd}
\end{equation}
Since $\Sigma$ is $\Q$-principal, the map $A^1(\Sigma, \Q) \to \Div!_{\Q}(\Sigma)=\MW_{d-1}(\Sigma, \Q)$ is surjective. Using the above diagram, we obtain the surjectivity of $A^1(\~\Sigma, \Q)\to \MW_{d-1}(\~\Sigma, \Q)$. By Theorem~\ref{thm:char_div-faithful}, this implies that $\~\Sigma$ is $\Q$-principal at $\conezero$, as desired.

\subsubsection{Closedness under blow-ups and blow-downs} \label{sec:closedness_bu_bd} Consider a tropical fan $\Sigma$ of dimension $d$ and let $\sigma$ be a cone of $\Sigma$ and $\rho$ a ray inside $\sigma$. Set $\Sigma' = \Sigma_{(\rho)}$. Assume that $\Sigma$ and $\Sigma'$ are $\Q$-principal at any nontrivial face. We need to prove $\Sigma$ is $\Q$-principal at $\conezero_{\Sigma}$ if and only if $\Sigma'$ is $\Q$-principal at $\conezero_{\Sigma'}$.

\ref{lem:ml3} Assume first that $\Sigma$ is $\Q$-principal. Let $D$ be a divisor in $\Sigma'$. Let $D^\rho$ be the induced divisor on $\Sigma'^\rho$. We infer the existence of a meromorphic function $f^\rho$ on $\Sigma'^\rho$ such that $\div(f^\rho) = a D^\rho$ for a non-zero integer $a$. Via the projection map $N_\R \to N^\rho_\R$, we view $f^\rho$ as a meromorphic function on the fan $\Sigma'' \subset \Sigma'$ consisting of the cones $\tau$ comparable with $\rho$, and extend it to a meromorphic function on full $\Sigma'$. By an abuse of the notation, we denote this function by $f^\rho$. The divisor $aD - \div(f^\rho)$ on $\Sigma'$ does not have any of the faces $\tau \supfaceeq \rho$ with $\tau \in \Sigma'_{d-1}$ in its support. This means it can be viewed as a divisor in $\Sigma$. Since $\Sigma$ is $\Q$-principal, we can find a meromorphic function $f$ on $\Sigma$ with $\div(f) = b(aD - \div(f^\rho))$ for a non-zero integer $b$. We infer that $ba D$ is the divisor of the meromorphic function $f + bf^\rho$ on $\Sigma'$, which shows that $D$ is $\Q$-principal, as desired.

\ref{lem:ml4} Assume now that $\Sigma'$ is $\Q$-principal. Let $D$ be a divisor in $\Sigma$. Viewing $D$ as a divisor in $\Sigma'$, we find a meromorphic function $f$ on $\Sigma'$ such that $\div(f) = aD$ for a non-zero integer $a$. Take a facet $\eta\supfaceeq\sigma$ in $\Sigma$. Denote by $\Sigma'\rest\eta$ the fan of support $\eta$ induced by $\Sigma'$. If $f$ is not linear on $\eta$, then one can find a face $\tau$ of dimension $d-1$ in $\Sigma'\rest\eta \setminus \Sigma$ such that $\ord_\tau(f)\neq0$. This is impossible since $D$ is supported on faces of $\Sigma$. Hence $f$ is linear on all the facets of $\Sigma$ containing $\sigma$. Therefore, $f$ is meromorphic on $\Sigma$ which proves that $aD$ is principal in $\Sigma$.

\subsubsection{Proof of Theorem~\ref{thm:stability_principal}} At this point, we have verified all the cases of Lemma \ref{lem:stability_meta_lemma} and this concludes the proof of Theorem~\ref{thm:stability_principal}. \qed

\subsection{$\mT$-stability of div-faithfulness} \label{sec:stability_div-faithful}

In this section we prove the following.

\begin{thm} \label{thm:stability_div-faithful}
The property for a tropical fan to be div-faithful is $\mT$-stable.
\end{thm}

Again, we use Lemma \ref{lem:stability_meta_lemma} with $\Csh$ the class of all tropical fans and with $P$ the predicate of being div-faithful at $\conezero$. Elements of $\Bsho_1$ are trivially div-faithful. We show that the four properties stated in the lemma are verified.

\subsubsection{Closedness under products \ref{lem:ml1}} Consider two tropical fans $(\Sigma, \omega!_\Sigma)$ and $(\Sigma',\omega!_{\Sigma'})$ that are div-faithful. We need to show that $(\Sigma \times \Sigma', \omega!_{\Sigma\times \Sigma'})$ is div-faithful at $\conezero$. We can assume without loss of generality that $\Sigma$ and $\Sigma'$ are unimodular (by Lemma \ref{lem:stability_meta_lemma}) and saturated at $\conezero$ (see Remark \ref{rem:independent_exterior_lattice}).

We use the decomposition
\[ \MW_{d+d'-1}(\Sigma\times \Sigma') = \MW_{d-1}(\Sigma) \otimes \MW_{d'}(\Sigma')\ \oplus\ \MW_d(\Sigma) \otimes \MW_{d'-1}(\Sigma'). \]
We have injections $\Z=A^0(\Sigma)\hookrightarrow \MW_d(\Sigma)$ and $A^0(\Sigma') \hookrightarrow \MW_{d'}(\Sigma')$. We have moreover a decomposition
\[ A^1(\Sigma \times \Sigma') \simeq A^1(\Sigma)\otimes A^0(\Sigma')\ \oplus\ A^0(\Sigma)\otimes A^1(\Sigma'). \]

The injectivity of the map $A^1(\Sigma \times \Sigma')\to\MW_{d+d'-1}(\Sigma\times \Sigma')$ now follows from the injectivity of the corresponding maps for $\Sigma$ and $\Sigma'$, given by Theorem~\ref{thm:char_div-faithful}, which shows that $\Sigma\times \Sigma'$ is div-faithful at $\conezero$, as desired.

\subsubsection{Closedness under tropical modifications \ref{lem:ml2}} Let $(\Sigma, \omega)$ be a tropical fan that is div-faithful. Let $(\~\Sigma,\~\omega)$ be a tropical modification of $\Sigma$ along a divisor $\Delta$ which is div-faithful. We need to prove that $\~\Sigma$ is div-faithful at $\conezero$. We can assume without loss of generality that $\Sigma$ is saturated at $\conezero$ and unimodular.

We use the commutative diagram \eqref{diag:tropical_modification_chow_minkowski}. This time, Theorem \ref{thm:char_div-faithful} implies that the bottom map is injective, from which we deduce that the top map is also injective. The same theorem then gives the result.

\subsubsection{Closedness under blow-ups and blow-downs} Consider a tropical fan $(\Sigma, \omega)$, and let $\sigma$ be a cone of $\Sigma$ and $\rho$ a ray inside $\sigma$. Set $\Sigma' = \Sigma_{(\rho)}$. Assume that $\Sigma$ and $\Sigma'$ are div-faithful at any of their non-zero faces. We need to show that $\Sigma$ is div-faithful at $\conezero$ if and only if $\Sigma'$ is div-faithful at $\conezero$.

\ref{lem:ml3} First, assume that $\Sigma$ is div-faithful at $\conezero$. Let $f$ be a meromorphic function on $\Sigma'$ whose divisor is trivial. Then, the induced meromorphic function $f^\rho$ on $\Sigma^\rho$ (see Section~\ref{sec:closedness_bu_bd}) verifies $\div(f^\rho)=0$. By assumption, $\Sigma^\rho$ is div-faithful. Thus, $f^\rho$ is linear, and we deduce that $f$ is linear on each face $\eta\supfaceeq\sigma$ of $\Sigma$. Thus, $f$ is conewise linear on $\Sigma$. We infer that $f$ is globally linear. This shows that $\Sigma'$ is div-faithful at $\conezero$.

\ref{lem:ml4} Second, assume that $\Sigma'$ is div-faithful. Let $f$ be a meromorphic function on $\Sigma$ such that $\div(f)=0$. Then, $f$ is also a meromorphic function on $\Sigma'$, and since its divisor is trivial, it is linear. This shows that $\Sigma$ is div-faithful at $\conezero$, as required.

\subsubsection{Proof of Theorem~\ref{thm:stability_div-faithful}} We have studied all the cases of Lemma \ref{lem:stability_meta_lemma}. This concludes the proof of Theorem~\ref{thm:stability_div-faithful}. \qed

\section{$\mT$-stability of Poincaré duality for div-faithful unimodular fans} \label{sec:PD_stable}

Let $(\Sigma,\omega)$ be a unimodular tropical fan of dimension $d$. Recall that $A^\bul(\Sigma)$ denotes the Chow ring of $\Sigma$ with coefficients in $\Z$, and that $A^\bul(\Sigma, \Q) \coloneqq A^\bul(\Sigma) \otimes \Q$. The degree map $\deg\colon A^d(\Sigma) \to \Z$ is given by the canonical class $\omega=\omega!_\Sigma\in\MW_d(\Sigma)\simeq A^d(\Sigma)^\dual$, and sends any element $x_\eta$ of $A^d(\Sigma)$, $\eta \in \Sigma_d$, to $\omega(\eta)$. Consider the pairing
\begin{align*}
  \Phi \colon A^\bul(\Sigma) \times A^{d-\bul}(\Sigma) & \to \Z, \\
  (\alpha, \beta) \in A^k(\Sigma) \times A^{d-k'}(\Sigma) & \mapsto
  \begin{cases}
    \deg(\alpha\cdot \beta) & \text{if $k = k'$,} \\
    0 & \text{otherwise.}
  \end{cases}
\end{align*}
We denote by $\Phi_\Q$ the induced bilinear pairing on $A^\bul(\Sigma, \Q)$.

\begin{defi}[Poincaré duality]
We say that $(\Sigma, \omega)$ verifies the \emph{Poincaré duality with $\Z$-coefficients} denoted $\PD_\Z$ if $\Phi$ is a perfect pairing. Similarly, we say $(\Sigma, \omega)$ verifies the \emph{Poincaré duality with $\Q$-coefficients} denoted $\PD_\Q$ if $\Phi_\Q$ is a perfect pairing.
\end{defi}

We have the following theorem.

\begin{thm} \label{thm:PD_stable}
Properties $\PD_\Z$ and $\PD_\Q$ are both $\mT$-stable in the class of div-faithful unimodular tropical fans.
\end{thm}

\begin{remark}
In this section, the orientation is allowed to take negative values on facets. We require tropical fans to be effective only in the treatment of positivity in the next section.
\end{remark}

Before providing the proof, we make a remark and state some consequences of this theorem.

\begin{remark} \label{remark:coprime}
If $(\Sigma, \omega)$ is a fan which verifies $\PD_\Z$, then, the weights $\omega(\eta)$, $\eta \in \Sigma_d$, are coprime. In particular, if\/ $\Sigma$ verifies $\Ploc{\PD_\Z}$, then $\Sigma$ is unitary. To see this, note that since $A^0(\Sigma) \simeq \Z$, Poincaré duality $\PD_\Z$ implies that $\deg\colon A^d(\Sigma) \to \Z$ is an isomorphism. It follows from the definition that the image of the degree map is $\gcd\left((\omega(\eta))_{\eta\in\Sigma_d}\right)\cdot\Z$. Therefore, the weights must be coprime. For any facet $\eta$, $\Sigma^\eta$ is a point with multiplicity $\omega(\eta)$, which verifies $\PD_\Z$ if and only if $\omega(\eta) = \pm1$. In particular, if $\Sigma$ verifies $\Ploc{\PD_\Z}$, then $\Sigma$ is unitary.
\end{remark}

\begin{thm}
$\Ploc{\PD_\Q}$ and $\Ploc{\PD_\Z}$ are both $\mT$-stable in the class of unimodular tropical fans.
\end{thm}

\begin{proof}
Let $\Sigma$ be a fan verifying $\Ploc{\PD_\Z}$. By Poincaré duality at $\conezero$, the map $A^1(\Sigma) \to \Div(\Sigma) \simeq \MW_{d-1} \simeq (A^{d-1}(\Sigma))^\dual$ is injective. Hence, $\Sigma$ is div-faithful at $\conezero$. Applying the same argument to all star fans, $\Ploc{\PD_\Z}$ implies div-faithfulness. From the definition of $\mT$-stability we infer that the $\mT$-stability of $\Ploc{\PD_\Z}$ in the class of unimodular tropical fans is equivalent to the $\mT$-stability of $\Ploc{\PD_\Z}$ in the class of div-faithful unimodular tropical fans. The result follows from Theorem~\ref{thm:PD_stable} and Proposition~\ref{prop:loc_P_stable}. The same argument works with rational coefficients.
\end{proof}

We directly get the following corollary.

\begin{cor}
Any unimodular quasilinear tropical fan verifies $\Ploc{\PD_\Q}$. If the fan is additionally unitary, then it verifies $\Ploc{\PD_\Z}$.
\end{cor}

The rest of this section is devoted to the proof of Theorem~\ref{thm:PD_stable}. We focus on the $\mT$-stability of $\PD_\Z$, as the proof with rational coefficients will be identical.

We check the properties required in Definition~\ref{def:t-stability} one by one.

\subsection{Closedness under products}

Let $(\Sigma,\omega)$ and $(\Sigma',\omega')$ be two unimodular tropical fans of dimension $d$ and $d'$, respectively. By Künneth decomposition, we have a ring isomorphism
\[ A^\bul(\Sigma \times \Sigma') \simeq A^\bul(\Sigma) \otimes A^\bul(\Sigma'). \]
The induced isomorphism between dual spaces sends $\deg!_{\Sigma\times\Sigma'}$ to $\deg!_\Sigma\otimes\deg!_{\Sigma'}$. We infer that $\PD_\Z(\Sigma \times \Sigma')$ holds provided that $\PD_\Z(\Sigma)$ and $\PD_\Z(\Sigma')$ both hold.

\subsection{Closedness under tropical modifications}

Let $(\Sigma,\omega)$ be a div-faithful tropical fan and let $(\~\Sigma,\omega)$ be a tropical modification of $(\Sigma,\omega)$. By Poincaré duality, $A^\bul(\Sigma)$ has no torsion. Applying Theorem~\ref{thm:invariant_chow_tropical_modification}, we obtain the Chow stability of tropical modification, namely that the natural ring morphism $A^\bul(\Sigma) \to A^\bul(\~\Sigma)$ is an isomorphism. By compatibility of the degree maps $\deg!_\Sigma$ and $\deg!_{\~\Sigma}$, we conclude that $\PD_\Z(\Sigma)$ implies $\PD_\Z(\~\Sigma)$. Note that, if we are working with rational coefficients instead, we can still apply Theorem~\ref{thm:invariant_chow_tropical_modification} and the result follows in a similar way.

\begin{remark}
In the proof, we did not require that the divisor along which we performed the tropical modification verify Poincaré duality. We thus get a slightly stronger result.
\end{remark}

\subsection{Closedness under blow-ups and blow-downs} \label{sec:Keel}

Consider a unimodular fan $\Sigma$ of dimension $d$ and let $\sigma$ be a cone in $\Sigma$. Let $\Sigma' = \Bl\Sigma\sigma$ be the unimodular stellar subdivision of $\Sigma$ obtained by stellar subdividing the cone $\sigma$. Denote by $\rho$ the new ray in $\Sigma'$, that is, $\rho=\R_{\geq 0}(\e_1+\dots+\e_{\dims{\sigma}})$ with $\e_1, \dots, \e_{\dims{\sigma}}$ the primitive vectors of the rays of $\sigma$. Recall that for any $\sigma\in\Sigma$, we have a surjective restriction map $\i^*_{\conezero\subfaceeq\sigma}\colon A^\bul(\Sigma) \to A^\bul(\Sigma^\sigma)$ described in Section~\ref{subsec:restr_gys}. We have the following key result on the relation between the Chow rings of $\Sigma$ and $\Sigma'$.

\begin{thm}[Keel's lemma] \label{thm:keel}
Let $\J$ be the kernel of the surjective map $\i^*_{\conezero\subfaceeq\sigma}\colon A^\bul(\Sigma)\to A^\bul(\Sigma^\sigma)$ and let
\[ P(T)\coloneqq\prod_{\zeta\subfaceeq \sigma \\ \dims\zeta=1}(x_\zeta+T). \]
We have an isomorphism
\[ \chi\colon \rquot{A^\bul(\Sigma)[T]}{(\J T+P(T))} \simto A^\bul(\Sigma') \]
which sends $T$ to $-x_\rho$ and which verifies
\[ \forall \zeta \in \Sigma_1, \qquad
\chi(x_\zeta) = \begin{cases}
x_\zeta+x_\rho & \text{if $\zeta \subfaceeq \sigma$,} \\
x_\zeta & \text{otherwise.}
\end{cases}
\]
In particular, we obtain a decomposition of $A^\bul(\Sigma')$ as
\begin{equation} \label{eq:keel}
A^\bul(\Sigma')\simeq A^\bul(\Sigma)\oplus A^{\bul-1}(\Sigma^\sigma)T \oplus \dots \oplus A^{\bul-\dims{\sigma}+1}(\Sigma^\sigma)T^{\dims{\sigma}-1}.
\end{equation}
\end{thm}

\begin{proof}
This follows from~\cite{Kee92}*{Theorem 1 in the appendix} for the map of toric varieties $\P_{\Sigma'} \to \P_{\Sigma}$. Here $P(T)$ is the polynomial in $A^\bul(\P_\Sigma)$ whose restriction in $A^\bul(\P_{\Sigma^\sigma})$ is the Chern polynomial of the normal bundle for the inclusion of toric varieties $\P_{\Sigma^\sigma} \hookrightarrow \P_{\Sigma}$.

We will provide an elementary proof of this lemma in Section~\ref{sec:keel_proof}.
\end{proof}

Let $(\Sigma,\omega)$ be a unimodular tropical fan and let $(\Sigma', \omega')$ be the tropical fan obtained as the result of the unimodular blow-up of the cone $\sigma$ in $\Sigma$. By Keel's lemma, we have $A^d(\Sigma) \simeq A^d(\Sigma')$. Moreover, as the proof of Keel's lemma given in Section~\ref{sec:keel_proof} shows, this isomorphism is compatible with the degree maps. This implies that $\deg!_{\Sigma} \colon A^d(\Sigma) \to \Z$ is an isomorphism if and only if $\deg!_{\Sigma'}\colon A^d(\Sigma') \to \Z$ is an isomorphism.

Let $\sigma$ be a face in $\Sigma$. Assume that $\PD_\Z(\Sigma^\sigma)$ holds. We need to show the equivalence of $\PD_\Z(\Sigma)$ and $\PD_\Z(\Sigma')$. By the preceding discussion, we can assume in the following that the degree maps are both isomorphisms.

\begin{lemma} \label{lem:lifting}
We have the following commutative diagram
\[ \begin{tikzcd}
A^{d-\dims\sigma}(\Sigma^\sigma) \arrow[r, "\cdot \left(-T^{\dims{\sigma}}\right)", "\sim"'] \arrow[rd, "\sim"{sloped, above}, "\deg"'] & A^d(\Sigma') \arrow[d, "\deg", "\sim"{sloped, below}] & \lar{\sim} A^d(\Sigma) \arrow[ld, "\deg", "\sim"{sloped, above}] \\
& \Z
\end{tikzcd} \]
\end{lemma}

\begin{proof}
Let $\beta$ be a top-degree element in $A^{d-\dims\sigma}(\Sigma^\sigma)$, and let $\alpha \in A^\bul(\Sigma')$ be a lifting of $\beta$ for the restriction map $\i^*_{\conezero \subface \sigma}$, which exists by \eqref{eqn:i_surjective_homeo-local} in Proposition \ref{lem:i_gys_basic_properties-local}.

Using the identity $P(T) = 0$, we get
\[ -T^{\dims \sigma} = \sum_{j=1}^{\dims \sigma} S_{j} T^{\dims \sigma -j} \]
with $S_j$ referring to the $j$-th symmetric function in the variables $x_\zeta\in A^1(\Sigma')$ for $\zeta$ a ray in $\sigma$. Therefore,
\[ -T^{\dims \sigma} \alpha = \sum_{j=1}^{\dims \sigma} S_{j} \alpha T^{\dims \sigma -j} = x_\sigma \alpha + \sum_{j=1}^{\dims \sigma -1} S_{j} \alpha T^{\dims \sigma -j}. \]
Since $\beta=\i^*_\sigma(\alpha)$ lives in the top-degree part of $A^\bul(\Sigma^\sigma)$, the products $S_j\alpha$ all belong to $\J$. We infer that the terms of the sum are all vanishing for $j=1,\dots, \dims \sigma -1$, and using Propositions~\ref{lem:i_gys_basic_properties-local} and \ref{lem:i_gys_basic_properties-local-tropical}, we infer that
\[ \deg_{\Sigma'}\circ\chi(-T^{\dims \sigma} \alpha) = \deg_\Sigma(x_\sigma \alpha) = \deg_\Sigma(\gys_{\sigma \supface \conezero} \circ \i^*_{\conezero \subface \sigma} (\alpha)) = \deg_{\Sigma^\sigma}(\beta). \]
The fact that $\cdot(-T^{\dims\sigma})$ is an isomorphism then follows from the rest of the diagram.
\end{proof}

Denote by $\Psi_k(\Sigma')$ the pairing $A^k(\Sigma') \times A^{d-k}(\Sigma') \to \Z$. We use similar notations for $\Sigma$ and $\Sigma^\sigma$. The above diagram allows to describe the pairing induced by $\Psi_k(\Sigma')$ on the different parts of the decomposition given by Keel's lemma.
\begin{itemize}
\item Between $A^k(\Sigma)$ and $A^{d-k}(\Sigma)$, we get the pairing $\Psi_k(\Sigma)$.
\item For any positive integer $i$, between $A^{k-i}(\Sigma^\sigma)T^i$ and $A^{d-\dims\sigma-k+i}(\Sigma^\sigma)\,T^{\dims\sigma-i}$ we get the pairing $-\Psi_{k-i}(\Sigma^\sigma)$.
\item For positive integers $i<j$, the pairing between $A^{k-i}(\Sigma^\sigma)T^i$ and $A^{d-\dims\sigma-k+j}(\Sigma^\sigma)\,T^{\dims\sigma-j}$ is trivial. Indeed, $A^{d-\dims\sigma+j-i}(\Sigma^\sigma)$ is trivial.
\item For any positive integer $i<\dims\sigma$, the pairing between $A^k(\Sigma)$ and $A^{d-k-i}(\Sigma^\sigma)T^i$ is trivial. Indeed, the product lives in $A^{d-i}(\Sigma^\sigma)T^i$ which is trivial since $d-i > d-\dims\sigma$.
\end{itemize}
We do not need to compute the pairing between the other parts. With respect to the decomposition given by Keel's lemma, the bilinear form $\Psi_k(\Sigma')$ can be written in the form of a block matrix consisting of bilinear maps, as follows.

\smallskip

\begin{equation}\label{eq:matrix}
\begin{tikzpicture}
\matrix[matrix of math nodes]{
  & \quad & A^{d-k}(\Sigma) & A^{d-k-\dims\sigma+1}(\Sigma^\sigma) T^{\dims\sigma-1} & \quad\cdots\quad & A^{d-k-1}(\Sigma^\sigma)T \\[ 1ex]
A^k(\Sigma) && \Psi_k(\Sigma) & 0 && \\
A^{k-1}(\Sigma^\sigma)T && 0 & -\Psi_{k-1}(\Sigma^\sigma) & \node (zero) {}; & \\
\vdots && \vdots & \node (ast) {}; & \rotatebox{-30}{$\cdots$} & \\
A^{k-\dims\sigma+1}(\Sigma^\sigma) T^{\dims\sigma-1} && 0 &&& -\Psi_{k-\dims\sigma+1}(\Sigma^\sigma) \\};
\draw[gray, thin, dotted] ($(zero)+(-1.7,1)$) --++ (4.2,0) --++ (0,-2.1) --++ (-1,0) -- cycle;
\path (zero) node[shift={(1.4,0)}] {\Large0};
\draw[gray, thin, dotted] ($(ast)+(2.6,-1)$) --++ (-3.4,0) --++ (0,1.4) --++ (1.2,0) -- cycle;
\path (ast) node[shift={(.2,-.5)}] {\LARGE*};
\path (-2.8,-.4) node {$\left(\rule[-1.4cm]{0pt}{2.8cm}\right.$};
\path (7.2,-.4) node {$\left.\rule[-1.4cm]{0pt}{2.8cm}\right)$};
\end{tikzpicture}
\end{equation}

The matrix is lower triangular. For the terms appearing on the diagonal, since we assume $\PD_\Z(\Sigma^\sigma)$, all the bilinear maps except the first one are non-degenerate. Thus, $\Psi_k(\Sigma')$ is non-degenerate if and only if $\Psi_k(\Sigma)$ is non-degenerate. This shows the properties $\PD_\Z(\Sigma)$ and $\PD_\Z(\Sigma')$ are equivalent provided that $\PD_\Z(\Sigma^\sigma)$ is verified.

The above reasoning also applies to Chow rings with rational coefficients.

\subsection{Proof of Theorem~\ref{thm:PD_stable}}

At this point, we have verified all the needed properties for the $\mT$-stability of the property $\PD_\Z$ in the class of unimodular tropical fans, and this concludes the proof of Theorem~\ref{thm:PD_stable}. The proof of the statement for $\PD_\Q$ is similar. \qed

\subsection{Proof of Keel's lemma} \label{sec:keel_proof}

In this section we give an elementary proof of Theorem~\ref{thm:keel}. First, we define a map
\[ \Phi \colon \Z[\x_\zeta]_{\zeta\in \Sigma_1} \to \Z[\x_{\zeta}]_{\zeta\in \Sigma'_1} \]
on the level of polynomial rings by setting
\[ \forall \zeta \in \Sigma_1, \qquad
\Phi(\x_\zeta) \coloneqq \begin{cases}
\x_\zeta+\x_\rho & \text{if $\zeta \subfaceeq \sigma$,}\\
\x_\zeta & \text{otherwise.}
\end{cases}
\]

\begin{claim}
We have
\[ \Phi(\ssI_{\Sigma}) \subset \ssI_{\Sigma'} \qquad \textrm{and} \qquad \Phi(\ssJ_{\Sigma}) \subset \ssJ_{\Sigma'}. \]
\end{claim}

\begin{proof}
The first assertion is clear from the definition of the blow-up. The second inclusion $\Phi(\ssJ_{\Sigma}) \subset \ssJ_{\Sigma'}$ follows from the relation $\e_\rho = \sum_{\zeta \subfaceeq \sigma\\\dims\zeta=1} \e_\zeta$, which implies for any $m\in M$, the equality
\[ \Phi\Bigl(\sum_{\zeta\in \Sigma_1} m(\e_\zeta)\x_\zeta\Bigr) = \sum_{\zeta \in \Sigma'_1} m(\e_\zeta)\x_\zeta. \qedhere \]
\end{proof}

We infer that the map $\Phi$ descends to a morphism of Chow rings
\[ \Phi\colon A^\bullet(\Sigma) \to A^{\bullet}(\Sigma'). \]
We now extend the morphism $\Phi$ to a morphism of rings
\[ \chi \colon A^\bullet(\Sigma)[T] \to A^{\bullet}(\Sigma') \]
by sending $T$ to $-x_\rho$.

Recall that $P(T) \coloneqq \prod_{\zeta\in \sigma\\\dims{\zeta}=1} (T+x_\zeta)$ and $\mathfrak J$ is the kernel of the restriction map $\i^*\colon A^\bullet(\Sigma) \to A^\bullet(\Sigma^\sigma)$.

\begin{claim}
We have $\chi(P(T)) = 0$ and $\chi\rest{\J T} =0$.
\end{claim}

\begin{proof}
By definition, $\chi\rest{A^\bullet(\Sigma)}$ coincides with $\Phi$. It follows that
\[ \chi(P(T)) = \prod_{\zeta\in \sigma\\\dims{\zeta}=1} x_\zeta \]
which is zero since the cone $\sigma$ does not belong to $\Sigma'$.

To see the second assertion, we use the commutative diagram
\[ \begin{tikzcd}
A^{\bullet}(\Sigma') \arrow[r, "\i^*"]  & A^\bullet(\Sigma'^{\rho}) \arrow[r, "\gys"] & A^{\bullet+1}(\Sigma') \\
A^{\bullet}(\Sigma) \arrow[r, "\i^*"] \arrow[u, "\Phi"] & A^\bullet(\Sigma^\sigma) \arrow[u]
\end{tikzcd} \]
where the second vertical map is obtained from the natural inclusion $\Sigma^\sigma_1 \hookrightarrow \Sigma'^\rho_1$. The image of an element $a$ in $\J$ in $A^{\bullet+1}(\Sigma')$ following first the bottom line, then going up and ending in $A^{\bullet+1}(\Sigma')$ is zero. On the other hand, the image of the same element going up first and then following the upper line is $\Phi(a)x_\rho = -\chi(aT)$. We infer that $\chi\rest{\J T} =0$.
\end{proof}

This implies that $\chi$ descends to a morphism of rings
\[ \chi\colon \rquot{A^\bullet(\Sigma)}{\bigl(\J T + (P(T))\bigr)} \to A^\bullet(\Sigma'). \]
Since $x_\zeta$ for $\zeta \in \Sigma'$ are all in the image of $\chi$, $\chi$ is surjective.

In order to prove $\chi$ is an isomorphism, we construct a morphism
\[ \Psi \colon A^\bullet(\Sigma') \to \rquot{A^\bullet(\Sigma)}{\bigl(\J T + (P(T))\bigr)}, \]
and from the definition, we get $\Psi \circ \chi =\id$ and $\chi\circ \Psi =\id$.

We first construct
\[ \Psi \colon \Z[\x_{\zeta}]_{\zeta \in \Sigma'_1} \to \Z[\x_{\zeta}]_{\zeta \in \Sigma_1}[T] \]
on the level of polynomial rings by
\[ \Psi(\x_\zeta) =\begin{cases}
  \x_{\zeta} & \textrm{if $\zeta$ is not a ray of $\sigma$ and $\zeta \neq \rho$,} \\
  \x_{\zeta}+ T & \textrm{if $\zeta$ is a ray of $\sigma$,} \\
  -T & \textrm{if $\zeta=\rho$.}
\end{cases} \]

\begin{claim}
We have $\Psi(\ssJ_{\Sigma'}) \subseteq \ssJ_{\Sigma}$.
\end{claim}

\begin{proof}
Let $m\in M$. We have, using $\e_\rho = \sum_{\zeta\subfaceeq \sigma \\\dims\zeta=1}\e_\zeta$,
\begin{align*}
\sum_{\zeta\in \Sigma'_1}m(\e_\zeta) \x_\zeta
  &= -m(\e_\rho)T + \sum_{\zeta \subfaceeq \sigma \\\dims\zeta=1} m(\e_\zeta) (\x_\zeta + T) + \sum_{\zeta \in \Sigma_1\\\zeta\not\subfaceeq \sigma} m(\e_\zeta)\Psi(\x_\zeta)\\
  & = \sum_{\zeta\in \Sigma_1}m(\e_\zeta) \x_\zeta \in \ssJ_{\Sigma}. \qedhere
\end{align*}
\end{proof}

\begin{claim}
We have $\Psi(\ssI_{\Sigma'}) \subseteq \ssI_{\Sigma} + \J T + (P(T))$.
\end{claim}

\begin{proof}
Let $S = \{\zeta_1, \dots, \zeta_k\}$ be a set of distinct rays in $\Sigma'$ which do not form a cone in $\Sigma'$ so that $\x_S = \prod_{\zeta\in S} \x_\zeta \in \ssI_{\Sigma'}$. We proceed by a case analysis.
\begin{itemize}
\item If all the rays of $\sigma$ are in $S$, then $\Psi(\x_S) \in (P(T))$.
\item If $S$ does not contain $\rho$ and it contains not all rays in $\sigma$. Let $\tau$ be the (possibly zero) face of $\sigma$ generated by those rays which belong to $S$ and let $S'$ be those rays in $S$ which are not in $\sigma$. We can write $\x_S = \x_\tau \x_{S'}$. Moreover, the rays in $S'$ do not form a cone with $\tau$, nor with $\sigma$, in $\Sigma$. We have
\[ \Psi(\x_S) = \x_{S'} \prod_{\zeta \subfaceeq \tau \\ \dims\zeta=1}(T+\x_\zeta). \]
Developing the product, we see that $\x_{S'}\x_\tau \in \ssI_{\Sigma}$, and all the other terms are in $\J T$. That is, $\Psi(\x_S) \in \ssI_{\Sigma}+\J T $.
\item In the remaining case, $S$ contains $\rho$ but not all rays of $\sigma$. Let $\tau$ be the (possibly zero) face of $\sigma$ generated by those rays in $S$ which belong to $\sigma$, and $S'$ all the rays of $S$ which are different from $\rho$ and do not belong to $\sigma$. Again, $S'$ and rays of $\sigma$ do not form a cone in $\Sigma$. It follows that
\[ \Psi(\x_S) = - \x_{S'} T \prod_{\zeta\in \tau \\\dims\zeta=1}(T+\x_\zeta) \in \J T. \qedhere \]
\end{itemize}
\end{proof}

We infer from the above claims that $\Psi$ descends to a morphism of rings
\[ \Psi \colon A^\bullet(\Sigma') \to \rquot{A^\bullet(\Sigma)}{\J T + (P(T))}. \]
We have $\Psi \circ \chi =\id$ and $\chi\circ \Psi =\id$ by verifying them on the level of generators, and this concludes the proof of Keel's lemma. \qed

\section{$\mT$-stability of being Chow-Kähler} \label{sec:kahler}

All through this section, the orientations appearing in the tropical fans will be effective.

In order to apply topological arguments in few places, we will work with the Chow ring $A^\bullet(\Sigma, \R)$ with real coefficients (but drop the mention of $\R$). If the considered ample classes are rational, then the stated results hold for the Chow ring with rational coefficients.

\smallskip
Let $(\Sigma, \omega)$ be an effective unimodular tropical fan of dimension $d$.

\begin{defi}[Hard Lefschetz property] Let $\ell$ be an element of $A^1(\Sigma)$. We say that $(\Sigma, \omega)$ verifies the
\emph{Hard Lefschetz property} $\HL(\ell)$ if $\Sigma$ verifies $\PD_\Q$, and for any non-negative integer $k\leq d/2$, the multiplication map $\ell^{d-2k}\colon A^k(\Sigma) \to A^{d-k}(\Sigma)$ is an isomorphism.
\end{defi}

For an element $\ell\in A^1(\Sigma)$, and $k\leq d/2$, we consider the bilinear form
\[ Q^k_{\ell} \colon A^k(\Sigma)\times A^k(\Sigma) \to \R \]
defined by
\begin{equation} \label{eqn:def_Q_l}
Q^k_{\ell}(a,b)\coloneqq\deg(\ell^{d-2k}ab), \qquad a,b \in A^k(\Sigma).
\end{equation}

Note that if $\PD_\Q$ holds, then $(\Sigma, \omega)$ verifies $\HL(\ell)$ if and only if $Q^k_\ell$ is perfect.

\begin{defi}[Hodge-Riemann] \label{defi:HR} A tropical fan $(\Sigma, \omega)$ verifies the \emph{Hodge-Riemann bilinear relations} $\HR(\ell)$ for an element $\ell \in A^1(\Sigma)$ if it verifies $\HL(\ell)$, and in addition, for any non-negative integer $k\leq d/2$, the signature of the perfect pairing $Q^k_\ell$ is given by the sum
\[ \sum_{i=0}^k (-1)^i\Bigl(\dim(A^i(\Sigma))-\dim(A^{i-1}(\Sigma))\Bigr). \qedhere \]
\end{defi}

If $(\Sigma, \omega)$ of dimension $d$ verifies $\HL(\ell)$, then for any $k\leq d/2$, we get the Lefschetz decomposition
\[ A^k(\Sigma) = P^k_\ell(\Sigma)\oplus \ell P_\ell^{k-1}(\Sigma) \oplus \dots \oplus \ell^kP_\ell^{0}(\Sigma) \]
with the \emph{primitive part} $P^k_\ell(\Sigma)$, dependent on $\ell$, defined by
\[ P_\ell^k(\Sigma) \coloneqq \ker\Bigl(\ell^{d-2k+1} \colon A^{k}(\Sigma) \to A^{d-k+1}(\Sigma)\Bigr). \]
Moreover, this is an orthogonal decomposition with respect to the bilinear form $Q^k_\ell$.

\begin{prop} \label{prop:equivalence_signature_definiteness}
Notations as in Definition~\ref{defi:HR}, assume that $\HL(\ell)$ holds. The following properties are equivalent.
\begin{itemize}
\item $\HR(\ell)$ holds.
\item For each $k\leq \frac d2$, the restriction of $(-1)^kQ_\ell^k$ to the primitive part $P_\ell^k(\Sigma)$ is positive definite.
\end{itemize}
\end{prop}

\begin{proof}
This follows by induction using the Lefschetz decomposition. We omit the details.
\end{proof}

\smallskip
Let $f$ be a conewise linear function on $\Sigma$ which takes real value $c_\zeta$ on the primitive vector $\e_\zeta$ and let $\ell = \sum_{\zeta} c_\zeta x_\zeta$ be the corresponding element of $A^1(\Sigma)$. We call $\ell$ \emph{ample} provided that $f$ is strictly convex on $\Sigma$.

\begin{defi}[Kähler package for the Chow ring] \label{def:KP-Chow}
We say that a effective quasi-projective unimodular tropical fan $(\Sigma, \omega)$ verifies the Kähler package for the Chow ring if it verifies $\HR(\ell)$ for any ample element $\ell \in A^1(\Sigma)$.
\end{defi}

\begin{defi}[Chow-Kähler fans]
We say that a effective quasi-projective unimodular tropical fan $(\Sigma, \omega)$ is Chow-Kähler if for each $\sigma \in \Sigma$, the tropical fan $(\Sigma^\sigma, \omega^\sigma)$ verifies the Kähler package for the Chow ring.
\end{defi}

This is the theorem we prove in this section.

\begin{thm} \label{thm:chow-KP_stable}
The property of being Chow-Kähler is $\mT$-stable in the class of effective quasi-projective unimodular tropical fans.
\end{thm}

As a corollary, using Proposition~\ref{prop:intersection-property}, we get the following result.

\begin{cor}
Any effective unimodular quasi-projective quasilinear fan is Chow-Kähler.
\end{cor}

In particular, since we proved that generalized Bergman fans are quasilinear, we get the following result.

\begin{cor}[Kähler package for matroids~\cites{AHK, ADH, BHMPW, BHMPW20b}]
The Chow ring of any quasi-projective generalized Bergman fan verifies the Kähler package. In particular, the Chow ring and the augmented Chow ring of matroids verify the Kähler package.
\end{cor}

The rest of this section is devoted to the proof of Theorem~\ref{thm:chow-KP_stable}. Sections~\ref{subsec:restriction_ample} and~\ref{subsec:ample_cone} gather some basic properties of ample classes. In Section~\ref{sec:HR_local_global}, we prove a local to global property for the Hodge-Riemann bilinear relations. Section~\ref{subsec:local_ascent_descent} establishes important ascent and descent properties for the Hodge-Riemann bilinear relations which allow to control stellar subdivisions and stellar assemblies. Finally, Section \ref{subsec:proof_chow-KP_stable} establishes the $\mT$-stability of being Chow-Kähler.

\subsection{Restriction of ample classes}\label{subsec:restriction_ample}

Let $(\Sigma, \omega)$ be a unimodular tropical fan. Let $f$ be a strictly convex conewise linear function on $\Sigma$. We denote by $\ell(f)$ the element of $A^1(\Sigma)$ defined by
\[ \ell(f)\coloneqq\sum_{\zeta\in\Sigma_1}f(\e_\zeta)x_\zeta. \]
Let $\sigma$ be a cone of $\Sigma$ and let $\phi$ be a linear form on $N_\R$ which coincides with $f$ on $\sigma$. The function $f-\phi$ induces a conewise linear function on $\Sigma^\sigma$ which we denote by $f^\sigma$.

\begin{prop} \label{prop:ell_f_independent}
Notations as above, we have
\[ \ell(f^\sigma)=\i^*_{\conezero\subfaceeq\sigma}(\ell(f)) \]
in the Chow ring of\/ $\Sigma^\sigma$. In particular, $\ell(f^\sigma)$ does not depend on the choice of the linear form~$\phi$. Moreover, if $f$ is strictly convex on $\Sigma$, then $f^\sigma$ is strictly convex on $\Sigma^\sigma$, that is, $\i^*$ sends ample classes in $A^1(\Sigma)$ to ample classes in $A^1(\Sigma^\sigma)$.
\end{prop}

\begin{proof}
Using the notations we introduced previously, we write $\zeta\sim\sigma$ if $\zeta$ is a ray in the link of $\sigma$, \ie, $\zeta\not\subfaceeq\sigma$ and $\sigma+\zeta$ is a face of $\Sigma$. Such rays are in one-to-one correspondence with rays of $\Sigma^\sigma$. Recall the following facts:
\[ \ell(\phi)=0, \qquad \i^*_{\conezero\subfaceeq\sigma}(x_\zeta)=0 \text{ if $\zeta\not\sim\sigma$ and $\zeta\not\subfaceeq\sigma$}, \qquad f(\e_\zeta)=\phi(\e_\zeta) \text{ if $\zeta\subfaceeq\sigma$}. \]
We have
\begin{align*}
\i^*_{\conezero\subfaceeq\sigma}(\ell(f))
  &= \i^*_{\conezero\subfaceeq\sigma}(\ell(f-\phi)) = \sum_{\zeta\in\Sigma_1}(f-\phi)(\e_\zeta)\i^*_{\conezero\subfaceeq\sigma}(x_\zeta) \\
  &= \sum_{\zeta\in\Sigma_1 \\ \zeta\sim\sigma}(f-\phi)(\e_\zeta)\i^*_{\conezero\subfaceeq\sigma}(x_\zeta) = \sum_{\zeta\in\Sigma^\sigma\\ \dims\zeta=1}f^\sigma(\e_{\zeta})x_{\zeta} = \ell(f^\sigma),
\end{align*}
which proves the first statement. The second claim is straightforward.
\end{proof}

\subsection{Ample cone and $\HR(\texorpdfstring{\ell}{l})$} \label{subsec:ample_cone}

The following is straightforward.

\begin{prop}
The set of ample elements in $A^1(\Sigma)$ is an open convex cone.
\end{prop}

\begin{defi}
The subset of $A^1(\Sigma)$ consisting of ample elements is called the \emph{ample cone} of $\Sigma$.
\end{defi}

\begin{prop} \label{prop:HRbis}
The property $\HR(\ell)$ is an open condition in $\ell \in A^1(\Sigma)$. Moreover, $\HR(\ell)$ is both a closed and an open condition in the space of all $\ell$ for which $\HL(\ell)$ is verified.
\end{prop}

\begin{proof}
The proposition directly follows from the following fact: the signature of non-degenerate bilinear maps remains constant under small deformations.
\end{proof}

\subsection{Hodge-Riemann for star fans of rays implies Hard Lefschetz} \label{sec:HR_local_global}

Let $\ell \in A^1(\Sigma)$ be an ample element. In this case, $\ell$ has a representative in $A^1(\Sigma)$ with strictly positive coefficients, i.e.,
\[ \ell =\sum_{\zeta \in \Sigma_1} c_\zeta x_\zeta \]
for scalars $c_\zeta >0$ in $\R$. For each $\zeta \in \Sigma_1$, define $\ell^\zeta = \i^*_{\conezero \subface \zeta}(\ell) \in A^1(\Sigma^\zeta)$ and let $\deg_\zeta$ be the degree map of $\Sigma^\zeta$.

\smallskip
We need the following result. It is used in~\cite{CM05} and~\cite{AHK}*{Proposition 7.15}.

\begin{prop} \label{prop:local_HR}
Assume that $(\Sigma,\omega)$ verifies Poincaré duality. If\/ $\HR(\Sigma^\zeta, \ell^\zeta)$ holds for all rays $\zeta\in\Sigma_1$, then we have $\HL(\Sigma, \ell)$.
\end{prop}

\begin{proof}
Let $k \leq \frac d2$. By Poincaré duality for $A^\bul(\Sigma)$, it will be enough to show that the multiplication map
\[ \ell^{d-2k}\cdot - \colon A^{k}(\Sigma) \longrightarrow A^{d-k}(\Sigma) \]
is injective. Let $a \in A^k(\Sigma)$ be such that $\ell^{d-2k}\cdot a =0$. We have to show that $a=0$.

\smallskip
There is nothing to prove if $k=d/2$, so assume $2k < d$. For each $\zeta \in \Sigma_1$, define
\[ a^\zeta \coloneqq\i^*_{\conezero \subface \zeta}(a). \]
It follows that
\[ (\ell^\zeta)^{d-2k} \cdot a^\zeta = \i^*_{\conezero \subface \zeta}(\ell^{d-2k}a) =0, \]
and so $a^\zeta$ lives in the primitive part $P^{k}_{\ell^\zeta}(\Sigma^\zeta) \subseteq A^k(\Sigma^\zeta)$, for each ray $\zeta \in \Sigma_1$.

\smallskip
Using now Propositions~\ref{lem:i_gys_basic_properties-local} and \ref{lem:i_gys_basic_properties-local-tropical}, for each $\zeta\in \Sigma_1$, we get
\[ \deg_\zeta\bigl((\ell^\zeta)^{d-2k-1}\cdot a^\zeta \cdot a^{\zeta}\bigr) = \deg_\zeta\bigl(\i^*_{\conezero \subface \zeta}(\ell^{d-2k-1}\cdot a) \cdot a^{\zeta}\bigr) = \deg\bigl(\ell^{d-2k-1}\cdot a \cdot x_\zeta \cdot a\bigr). \]

We infer that
\[ \sum_{\zeta \in \Sigma_1} c_\zeta \deg_\zeta\bigl((\ell^\zeta)^{d-2k-1}\cdot a_\zeta \cdot a_{\zeta}\bigr) = \deg\bigl(\ell^{d-2k-1}\cdot a \cdot \bigl(\,\sum_{\zeta }c_\zeta x_\zeta\bigr) \cdot a\bigr) = \deg(\ell^{d-2k}\cdot a\cdot a)=0. \]

By $\HR(\Sigma^\zeta,\ell^\zeta)$, since $a^\zeta\in P^{k}_{\ell^\zeta}(\Sigma^\zeta)$, we have $(-1)^k \deg_\zeta\bigl((\ell^\zeta)^{d-2k-1}\cdot a^\zeta \cdot a^{\zeta}\bigr) \geq 0$ with equality if and only if $a^\zeta =0$. Since $c_\zeta >0$ for all $\zeta$, we conclude that $a^\zeta =0$ for all $\zeta \in \Sigma_1$.

\smallskip
Applying Proposition~\ref{lem:i_gys_basic_properties-local} once more, we infer that
\[ x_\zeta a = \gys_{\zeta \supface \conezero}\circ \,\i^*_{\conezero \subface \zeta} (a) = \gys_{\zeta \supface \conezero}(a^\zeta) =0. \]

Since the elements $x_\zeta$ generate the Chow ring, Poincaré duality for $A^\bul(\Sigma)$ implies that $a=0$, and the proposition follows.
\end{proof}

We deduce the following result.

\begin{prop} \label{prop:HR-oneall}
Let $(\Sigma, \omega)$ be an effective and unimodular tropical fan of dimension $d$ which verifies the property $\Ploc \PD_\Q$. Assume for each cone $\sigma \in \Sigma$, there exists an ample element $a(\sigma) \in A^1(\Sigma^\sigma)$ such that $\HR(a(\sigma))$ holds. Then, $(\Sigma, \omega)$ is Chow-Kähler.
\end{prop}

\begin{proof}
Proceeding by induction, we can assume $(\Sigma^\sigma, \omega^\sigma)$ is Chow-Kähler for all $\sigma \neq \conezero$. We need to show that for any ample element $\ell \in A^1(\Sigma)$, $\HR(\ell)$ holds.

If $\ell \in A^1(\Sigma)$ is an ample element, for any ray $\zeta\in \Sigma$, we get an ample element $\ell^\zeta$ in $A^1(\Sigma^\zeta)$. By the hypothesis of our induction, $\HR(\ell^\zeta)$ holds. Applying Proposition~\ref{prop:local_HR}, we deduce that $\HL(\Sigma, \ell)$ holds for any element $\ell$ of the ample cone.

By Proposition~\ref{prop:HRbis}, the set of $\ell$ which verify $\HR(\ell)$ is both open and closed in the set of all $\ell$ which verify $\HL(\ell)$. By assumption, there is an ample element $a(\conezero)$ in the ample cone of $\Sigma$ such that $\HR(a(\conezero))$ holds. Since any element $\ell$ in the ample cone of $\Sigma$ verifies $\HL(\ell)$, we infer that $\HR(\ell)$ holds for any element of the ample cone of $\Sigma$. This proves that $(\Sigma, \omega)$ is Chow-Kähler.
\end{proof}

\subsection{Ascent and Descent} \label{subsec:local_ascent_descent}

Let $(\Sigma, \omega)$ be an effective unimodular tropical fan. Let $h$ be a strictly convex conewise linear function on $\Sigma$ and denote by $\ell$ the corresponding ample element of $A^1(\Sigma)$. Let $\Sigma'$ be the fan obtained from $\Sigma$ by unimodular blow-up of a cone $\sigma\in\Sigma$. Denote by $\rho$ the new ray in $\Sigma'$. The function $h$ defines a conewise linear function on $\Sigma'$ that we denote by $h'$. Since $\e_\rho =\sum_{\zeta\subface \sigma\\ \dims\zeta=1} \e_\zeta$, we get $h'(\e_\rho) = \sum_{\zeta \subface \sigma \\ \dims\zeta=1}h(\e_\zeta)$. The corresponding element $\ell' \coloneqq \ell(h')$ in $A^1(\Sigma')$ is thus given by
\[ \ell'=\sum_{\zeta \in \Sigma_1} h(\e_\zeta)x_\zeta+\bigl(\sum_{\zeta \subface \sigma \\ \dims\zeta=1}h(\e_\zeta)\bigr)x_\rho. \]

Notice that the definition of $\ell'$ only depends on the class $\ell$, and not on the chosen representative $h$.

\begin{thm} \label{thm:ascent_descent}
We have the following properties.
\begin{itemize}
\item \emph{(Ascent)} Assume the property $\HR(\Sigma^\sigma,\ell^\sigma)$ holds. Then, $\HR(\Sigma, \ell)$ implies $\HR(\Sigma',\ell'-\epsilon x_\rho)$ for any small enough $\epsilon>0$,
\item \emph{(Descent)} We have the following partial inverse: if both the properties $\HR(\Sigma^\sigma,\ell^\sigma)$ and $\HL(\Sigma,\ell)$ hold, and if we have the property $\HR(\Sigma',\ell'- \epsilon x_\rho)$ for any small enough $\epsilon>0$, then we have $\HR(\Sigma, \ell)$.
\end{itemize}
\end{thm}

Before proving the theorem, we need to introduce some basic results about graded algebras verifying Hodge-Riemann bilinear relations.

Let $B^\bul$ be any finite dimensional graded vector space. Let $d$ be an integer such that $\dim(B^k) = \dim(B^{d-k})$ for all $k\in\Z$. We will refer to $d$ as the \emph{fundamental degree} of $B^\bul$. (The terminology is justified by the observation that $B^d$ is the piece which contains the fundamental class of a Poincaré duality space.)

Let $Q$ be a symmetric bilinear form on $B^{\bul\leq d/2}=\bigoplus_{k \leq d/2} B^k$ which decomposes as $Q = \bigoplus_{k \leq d/2} Q^k$ with $Q^k$ a symmetric bilinear form on $B^k$. We denote by $\HR(B^\bul,Q)$ the Hodge-Riemann bilinear relations for the pair defined to be the property that for any non-negative integer $k\leq d/2$, $Q^k$ is a perfect pairing of signature
\[ \sum_{i\leq k} (-1)^i\Bigl(\dim(B^i)-\dim(B^{i-1})\Bigr). \]
In particular, small perturbations of the symmetric bilinear forms $Q^k$ preserve the Hodge-Riemann bilinear relations. The following lemma is straightforward.

\begin{lemma} \label{lem:HR_sum}
Let $(B^\bul, Q)$ and $(C^\bul, Q')$ be two graded vector spaces of the same fundamental degree endowed with symmetric bilinear forms as above. If any two vector spaces among $(B^\bul,Q)$, $(C^\bul,Q')$ and their orthogonal sum $(B^\bul\oplus C^\bul,Q\oplus Q')$ verify Hodge-Riemann bilinear relations, so does the third.
\end{lemma}

For any integer $k$, we define the shift of the pair $(B^\bul, Q)$ by $k$ denoted by $(B^\bul, Q)[k]$ by
\[ (B^\bul, Q)[k] \coloneqq (B^{\bul+k}, (-1)^k Q). \]
This pair is of fundamental degree $d-2k$. The sign $(-1)^k$ in front of $Q$ is chosen so that the Hodge-Riemann bilinear relations are preserved by shifts.

\smallskip

Back to the situation we are interested in, assume now that $B^\bul$ is a finite dimensional graded algebra of fundamental degree $d$ which is zero in negative degrees. Assume moreover there is a degree map $\deg\colon B^d \to \R$. For an element $\ell$ in $B^1$, we denote by $Q_\ell = \bigoplus_{k \leq d/2} Q^k_\ell$ the quadratic form defined as in Equation \eqref{eqn:def_Q_l}, that is, $Q^k_{\ell}(a,b)\coloneqq\deg(\ell^{d-2k}ab)$ for $a,b \in B^k$.

We have the following result.

\begin{prop}[Stability of $\HR$ under products] \label{prop:HR_product}
Let $B^\bul$ and $C^\bul$ be two finite dimensional graded algebras of respective fundamental degrees $d$ and $d'$. Let $\ell \in B^1$ and $\ell' \in C^1$. Then
\[ \HR(B^\bul, Q_\ell)\text{ and\/ }\HR(C^\bul, Q_{\ell'}) \quad \Longrightarrow \quad \HR(B^\bul \otimes C^\bul, Q_{\ell \otimes 1 + 1 \otimes \ell'}). \]
\end{prop}

\begin{proof}
Assume the assumptions on the left hand side of the implication hold. For each piece $B^k$, $k \leq d/2$, we get the primitive decomposition $B^k =P_\ell^k \oplus \ell P_\ell^{k-1} \oplus \dots \oplus \ell^{k}P_\ell^0$. We take an orthogonal basis $(b_{k,i})_{i \in [\dim(B^k)-\dim(B^{k-1})]}$ of the primitive part $P_\ell^k$. We can thus decompose $B^\bul$ as an orthogonal sum
\[ B^\bul = \bigoplus_{k=0}^{d/2} \bigoplus_i \Vect\bigl(b_{k,i}, \ell b_{k,i}, \dots, \ell^{d-2k}b_{k,i}\bigr). \]
The restriction of the bilinear form $Q_\ell$ on the term corresponding to some $b_{k,i}$ is isomorphic to the pair $(\rquot{\R[\x]}{(\x^{d-2k+1})}, Q_\x)[-k]$ for the degree map $\deg(\x^{d-2k}) = |\deg(b_{k,i}^2\ell^{d-2k})|$. Using a similar decomposition $(c_{l,j})_{l,j}$ for $C^\bul$, we get an orthogonal decomposition of the product:
\begin{align*}
  \bigoplus_{k=0}^{d/2} \bigoplus_{l=0}^{d'\!/2} \bigoplus_{i} \bigoplus_{j} b_{k,i}c_{l,j} \bigl(\rquot{\R[\x]}{(\x^{d-2k+1})} \otimes \rquot{\R[\y]}{(\y^{d'-2l+1})}, Q_{\x \otimes 1 + 1 \otimes \y}\bigr)[-k - l].
\end{align*}
Hence, we are reduced to verify the statement in the case the two algebras are $\rquot{\R[\x]}{(\x^{r+1})}$ and $\rquot{\R[\y]}{(\y^{s+1})}$, for two non-negative integers $r$ and $s$, and $\ell$ and $\ell'$ are multiplications by $\x$ and $\y$, respectively. The proof in this case can be obtained either by using Hodge-Riemann property for the complex projective variety $\mathbb C\P^{r} \times \mathbb C\P^s$, or by the direct argument given in~\cite{BBFK02}*{Proposition 5.7}, or still by the combinatorial argument in~\cite{McD11}*{Lemma 2.2} and~\cite{AHK}*{Lemma 7.8} based on the use of Gessel-Viennot-Lindstr\"om lemma~\cites{GV85, Lin73}.
\end{proof}

\smallskip
We use the notations of Section~\ref{sec:Keel}. Let $d$ be the dimension of $\Sigma$. By Keel's lemma, we have $A^d(\Sigma) \simeq A^d(\Sigma')$, and under this isomorphism, the canonical element $\varpi!_\Sigma$ of $\Sigma$ gets identified with the canonical element $\varpi!_{\Sigma'}$ of $\Sigma'$. We denote by $\alpha$ a lifting of $\varpi!_{\Sigma^\sigma}$ in $A^{d-\dims \sigma}(\Sigma)$. By Lemma~\ref{lem:lifting}, the canonical element $\varpi!_{\Sigma'}$ can be identified with $-T^{\dims \sigma} \alpha$ in Keel's decomposition $A^\bul(\Sigma') \simeq \rquot{A^\bul(\Sigma)[T]}{(\J T+P(T))}$, and we have the compatibilities of the degree maps given in the lemma.

The proof of Theorem~\ref{thm:ascent_descent} is based on a deformation argument. The idea is used as a way to derive Grothendieck's standard conjecture of Lefschetz and Hodge type for a blow-up from the result on the base, see~\cite{Ito05}*{Section 3}. A similar argument is used in~\cite{AHK}.

\smallskip
Let $h$ be a strictly convex conewise linear function on $\Sigma$, $h'$ the conewise linear function on $\Sigma'$ induced by $h$, and denote by $\ell \in A^1(\Sigma)$ and $\ell'\in A^1(\Sigma')$ the corresponding elements. Consider the decomposition given by Theorem~\ref{thm:keel}
\begin{equation}\label{eq:keel-appendix}
A^\bul(\Sigma')\simeq A^\bul(\Sigma)\oplus A^{\bul-1}(\Sigma^\sigma)T \oplus \dots \oplus A^{\bul-\dims{\sigma}+1}(\Sigma^\sigma)T^{\dims{\sigma}-1}.
\end{equation}

We define the pair $(D^\bul, Q)$ by
\begin{equation}\label{eq:algebra_D}
\bigl(A^\bul(\Sigma),Q_\ell\bigr) \quad \oplus \quad \bigl(A^\bul(\Sigma^\sigma) \otimes \rquot{\R[T]}{(T^{\dims\sigma-1})},Q_{\ell^\sigma\otimes1 + 1\otimes T}\bigr)[-1].
\end{equation}
The pair $(\rquot{\R[T]}{(T^{\dims\sigma-1})}, Q_T)$ trivially verifies the Hodge-Riemann bilinear relations. Once shifted by one, it can be identified with the graded vector space $T\cdot\rquot{\R[T]}{(T^{\dims\sigma-1})}$ endowed with the suitable bilinear pairing. Via this identification and the decompositions given in~\eqref{eq:keel-appendix} and \eqref{eq:algebra_D}, we identify $D^\bul$ with $A^\bul(\Sigma')$ as graded real vector spaces. We thus get a bilinear form $Q$ on $A^{\bul\leq d/2}(\Sigma')$.

For positive $\epsilon>0$, we define the linear automorphism
\[ S_\epsilon\colon A^\bul(\Sigma')\to A^\bul(\Sigma') \]
of degree $0$ which is the identity map $\id$ on $A^\bul(\Sigma)$, and multiplication by $\epsilon^{k-\frac{\dims{\sigma}}2}$ on each $A^\bul(\Sigma^\sigma)T^k$ in the direct sum decomposition~\eqref{eq:keel-appendix}. We obtain the symmetric bilinear forms
\[Q_{\ell'+\epsilon T}\circ S_\epsilon \colon A^{\bul\leq d/2}(\Sigma') \times A^{\bul\leq d/2}(\Sigma') \to \R.\]
Here $S_\epsilon$ acts by diagonal action.

\begin{lemma}\label{lem:limit}
As $\epsilon$ tends to zero, the bilinear forms $Q_{\ell'+\epsilon T}\circ S_\epsilon$ admit a limit $Q_{\lim}$ defined on $A^{\bul\leq d/2}(\Sigma')$. Via the identification $A^{\bul\leq d/2}(\Sigma') = D^{\bul\leq d/2}$ as graded real vector spaces, we have $Q_{\lim} = Q$.
\end{lemma}

We postpone the proof of the lemma. Using the lemma, we deduce the theorem.

\begin{proof}[Proof of Theorem \ref{thm:ascent_descent}]
Assume that $\HR(A^\bul(\Sigma^\sigma), \ell^\sigma)$ holds. Applying Proposition \ref{prop:HR_product} and Lemma \ref{lem:HR_sum} to Equation \eqref{eq:algebra_D}, we get that $\HR(D^\bul, Q)$ is equivalent to $\HR(A^\bul(\Sigma),\ell)$.

\subsubsection*{Proof of the ascent property} Since $\HR(A^\bul(\Sigma),\ell)$ holds, so does $\HR(D^\bul, Q)$. By Lemma~\ref{lem:limit}, we can identify $Q$ with the limit bilinear form $Q_{\lim}$. Since Hodge-Riemann bilinear relations remains true for small perturbations of $Q$, we deduce that the property $\HR(A^\bul(\Sigma'),Q_{\ell'+\epsilon T}\circ S_\epsilon)$ holds for any small enough value of $\epsilon>0$. Since the automorphism $S_\epsilon$ preserves the signature, we infer that $\HR(A^\bul(\Sigma'),\ell'+\epsilon T)$ holds. Since $T$ corresponds to $-x_\rho$ via the decomposition~\ref{eq:keel-appendix}, we infer that the property $\HR(\Sigma', \ell'-\epsilon x_\rho)$ holds, and the ascent property in Theorem~\ref{thm:ascent_descent} follows.

\subsubsection*{Proof of the descent property}
In order to get $\HR(A^\bul(\Sigma), \ell)$, it suffices to prove $\HR(D^\bul,Q)$. Since we assume that $Q_\ell^k$ is non degenerate for $k\leq d/2$, so is $Q^k$.

By the hypothesis, we have $\HR(A^\bul(\Sigma'), \ell'-\epsilon x_\rho)$ for small enough values of $\epsilon>0$. Since $S_\epsilon$ is an automorphism, we deduce that $\HR(A^\bul(\Sigma'), Q_{\ell'-\epsilon x_\rho}\circ S_\epsilon)$ holds for $\epsilon>0$ small enough. We are in the situation where the limit $Q$ of the family $Q_{\ell'+\epsilon T}\circ S_\epsilon$ is non-degenerate, and we can apply a reasoning similar to Proposition~\ref{prop:HRbis} to deduce that the signature of $Q$ is the same as the signature of $Q_{\ell'+\epsilon T}\circ S_\epsilon$ for any small enough positive $\epsilon$. Hence, $\HR(D^\bul,Q)$ is verified.

This finishes the proof of Theorem~\ref{thm:ascent_descent}.
\end{proof}

We are thus left to prove Lemma~\ref{lem:limit}.

\begin{proof}[Proof of Lemma \ref{lem:limit}] Consider the decomposition \eqref{eq:keel-appendix}. The decomposition for degree $k$ part $A^{k}(\Sigma')$ of the Chow ring has pieces $A^{k}(\Sigma)$ and $A^{k-i}(\Sigma^\sigma) T^{i}$ for $1 \leq i \leq \min\{k, \dims \sigma -1\}$. Consider the bilinear form on the Chow ring of $\Sigma'$ given by the degree map of $A^\bul(\Sigma')$, as in Section~\ref{sec:Keel}. Each piece $A^{k-i}(\Sigma^\sigma)T^i$ is orthogonal to the piece $A^{d-k}(\Sigma)$ as well as to all the piece $A^{d-k-j}(\Sigma^\sigma)T^j$ for $j < \dims \sigma -i$, in Keel's decomposition of $A^{d-k}(\Sigma')$, see the matrix~\eqref{eq:matrix}.

We now work out the form of the matrix of the bilinear forms $Q_{\ell'+\epsilon T}\circ S_\epsilon$ in degree $k$ with respect to~\eqref{eq:keel-appendix}. First, note that for $a,b\in A^k(\Sigma)$, we have
\[ Q_{\ell'+\epsilon T}\circ S_\epsilon(a,b)=\deg_{A^\bul(\Sigma')}(a\ell^{\prime d-2k}b)+\O(\epsilon^{\dims \sigma}) = \deg_{A^\bul(\Sigma)}(a\ell^{d-2k}b)+\O(\epsilon^{\dims\sigma}) . \]

Second, for an element $a\in A^k(\Sigma)$ and an element $b\in A^{k-i}(\Sigma^\sigma)T^i$, we get the existence of an element $c\in A^{d-2k-\dims\sigma+2i}(\Sigma')$ such that
\[ Q_{\ell'+\epsilon T}\circ S_\epsilon(a,b)=\deg_{A^\bul(\Sigma')}(ac(\epsilon T)^{\dims\sigma-i}\epsilon^{-\dims\sigma/2+i}b)=\O(\epsilon^{\dims \sigma/2}). \]

Third, for an element $a\in A^{k-i}(\Sigma^\sigma)T^i$ and an element $b\in A^{k-j}(\Sigma^\sigma)T^j$, with $i+j\leq\dims\sigma$, we get, setting $n\coloneqq\binom{d-2k}{\dims\sigma-i-j}$,
\begin{align*}
  &Q_{\ell'+\epsilon T}\circ S_\epsilon(a,b)\\
  &\qquad= \deg_{A^\bul(\Sigma')}\Big(\epsilon^{-\dims\sigma/2+i}a\big(n(\epsilon T)^{\dims\sigma-i-j}\ell'^{d-2k-\dims\sigma+i+j}+\O(\epsilon^{\dims\sigma-i-j+1})\big)\epsilon^{-\dims\sigma/2+j}b\Big) \\
  &\qquad=-n\deg_{A^\bul(\Sigma^\sigma)}(a(\ell^\sigma)^{d-\dims\sigma-2k+i+j}b)+\O(\epsilon).
\end{align*}
Finally, if $a\in A^{k-i}(\Sigma^\sigma)T^i$ and $b\in A^{k-j}(\Sigma^\sigma)T^j$, with $i+j>\dims\sigma$, we get
\[ Q_{\ell'+\epsilon T}\circ S_\epsilon(a,b)=\deg_{A^\bul(\Sigma)}(\epsilon^{-\dims\sigma/2+i}a\O(1)\epsilon^{-\dims\sigma/2+j}b)=\O(\epsilon^{i+j-\dims\sigma}). \]

To conclude, we observe that the limit $\lim_{\epsilon\to 0}Q_{\ell'+\epsilon T}\circ S_\epsilon$ exists, and moreover, the pair $(A^\bul(\Sigma'), \lim_{\epsilon\to 0}Q_{\ell'+\epsilon T}\circ S_\epsilon)$ is identified with the pair $(D^\bul, Q)$,
as claimed.
\end{proof}

\subsection{Proof of Theorem~\ref{thm:chow-KP_stable}} \label{subsec:proof_chow-KP_stable}

We now proceed to the proof of $\mT$-stability of being Chow-Kähler. We use the $\mT$-stability meta Lemma~\ref{lem:stability_meta_lemma} with $\Csh$ the class of effective quasi-projective unimodular tropical fans and $P$ the predicate that a tropical fan $\Sigma$ verify the Kähler package for the Chow ring, Definition~\ref{def:KP-Chow}. It is easy to see that elements of $\Bsho_1$ are in $\Csh$ and verify the Kähler package for the Chow ring. We show that properties \ref{lem:ml1}-\ref{lem:ml2}-\ref{lem:ml3}-\ref{lem:ml4} stated in the lemma are verified.

\subsubsection{Closedness under products \ref{lem:ml1}} \label{sec:product_KP}
This follows directly from Proposition \ref{prop:HR_product}, using the observation that any strictly convex function $f$ on $\Sigma^1 \times \Sigma^2$ is of the form $\pi_1^*(f_1) + \pi_2^*(f_2)$ where $f_i$ is a strictly convex element of $\Sigma^i$, and $\pi_i\colon \Sigma \to \Sigma^i$ is the projection, for $i \in \{1,2\}$.

\subsubsection{Closedness under tropical modifications \ref{lem:ml2}} \label{sec:trop_modif}

Let $(\Sigma,\omega)$ be a fan verifying the Kähler package for the Chow ring. In particular $\Sigma$ is div-faithful. Let $(\~\Sigma,\omega)$ be a tropical modification of $(\Sigma,\omega)$ along an effective divisor $\Delta=\div(f)$, for $f$ a meromorphic function on $\Sigma$. Applying Theorem~\ref{thm:invariant_chow_tropical_modification}, we obtain the Chow stability of tropical modification, namely that the natural ring morphism $A^\bul(\Sigma) \to A^\bul(\~\Sigma)$ is an isomorphism.

\smallskip
Let $\~\ell=\~\ell(\~h) \in A^1(\~\Sigma)$ be an element associated to a strictly convex function $\~h$ on $\~\Sigma$. Adding a linear map if necessary, we can assume that $\~h$ is zero on the new ray of the tropical modification. There thus exists a function $h$ on $\Sigma$ such that $\~h$ is the pullback of $h$ by the projection $\pi\colon \~\Sigma \to \Sigma$. It is easy to see that $h$ is strictly convex on $\Sigma$. The corresponding class $\ell \in A^1(\Sigma)$ coincides with $\~\ell$ under the isomorphism given by Theorem~\ref{thm:invariant_chow_tropical_modification}. We deduce $\HR(\~\Sigma, \~\ell)$ from $\HR(\Sigma,\ell)$.

\subsubsection{Closedness under blow-ups and blow-downs \ref{lem:ml3}-\ref{lem:ml4}} \label{sec:hr_bd_bu}

Let $(\Sigma,\omega)$ be a unimodular tropical fan of dimension $d$ and let $h$ be a strictly convex conewise linear function on $\Sigma$. Denote by $\ell$ the corresponding ample element of $A^1(\Sigma)$. Let $\Sigma'$ be the fan obtained from $\Sigma$ by unimodular blow-up of a cone $\sigma\in\Sigma$. Denote by $\rho$ the new ray in $\Sigma'$. The function $h$ defines a conewise linear function on $\Sigma'$ that we denote by $h'$, and we get the corresponding element $\ell'$ associated to $h'$ in $A^1(\Sigma')$ given by
\[\ell'=\sum_{\zeta\in \Sigma_1}h(\e_\zeta)x_\zeta+\bigl(\sum_{\zeta \subface \sigma \\ \dims\zeta=1}h(\e_\zeta)\bigr)x_\rho.\]

\begin{prop}
For any small enough $\epsilon>0$, the element $\ell' - \epsilon x_\rho$ of $A^1(\Sigma')$ is ample.
\end{prop}

\begin{proof}
This is the element of $A^1(\Sigma')$ associated to the function $h' - \epsilon {\mathbf 1}_\rho$, with $\mathbf{1}_\rho$ the conewise linear function on $\Sigma'$ which takes value one on $\e_\rho$ and value zero on all the other rays. The conewise linear function $h' - \epsilon \mathbf{1}_\rho$ on $\Sigma'$ is strictly convex.
\end{proof}

Using Lemma~\ref{lem:stability_meta_lemma}, we can assume that $\Sigma^\tau$ verifies $\HR(\ell^\tau)$ for all $\tau\neq \conezero$ in $\Sigma$.

\begin{prop} \label{prop:HR-trans}
Notations as above, the following statements are equivalent.
\begin{enumerate}
\item \label{hr:trans1} We have $\HR(\Sigma, \ell)$.
\item \label{hr:trans2} The property $\HR(\Sigma',\ell'-\epsilon x_\rho)$ holds for any small enough $\epsilon>0$.
\end{enumerate}
\end{prop}

\begin{proof}
By assumption $\HR(\Sigma^\sigma, \ell^\sigma)$ holds. Theorem~\ref{thm:ascent_descent} leads to the implication $\ref{hr:trans1} \Rightarrow \ref{hr:trans2}$.

For the other implication $\ref{hr:trans2} \Rightarrow\ref{hr:trans1}$, by assumption $\HR(\Sigma^\zeta, \ell^\zeta)$ hold for any ray $\zeta$ in $\Sigma$. From Proposition~\ref{prop:local_HR}, we deduce that $\HL(\Sigma, \ell)$ holds. Applying now the descent part of Theorem~\ref{thm:ascent_descent} gives the result.
\end{proof}

As a consequence, we conclude by Proposition \ref{prop:HR-oneall} that $\Sigma$ is Chow-Kähler if and only if $\Sigma'$ is Chow-Kähler.

\subsubsection{Proof of Theorem~\ref{thm:chow-KP_stable}}

At this point, we have verified all the cases of Lemma~\ref{lem:stability_meta_lemma}. We deduce that the property of being Chow-Kähler is $\mT$-stable among the class of effective quasi-projective unimodular fans. \qed

\section{Further discussions and examples} \label{sec:examples}

In this final section, we provide a collection of examples to which we referred in the text, and complement this with remarks and questions in order to clarify the concepts introduced in the paper.

\subsection{An alternate definition of irreducible components} \label{sec:alternate_irreducible}

In Section~\ref{sec:irreducible}, we defined the notion of irreducible components of a tropical fan $\Sigma$. An alternate definition is the following. A subfan $\Delta$ of $\Sigma$ is an irreducible component of $\Sigma$ if it is the support of a nonzero element in $\MW_d(\Sigma)$, and if it is minimal among the subfans with this property. In the case $\Sigma$ is normal, this definition coincides with the one given in Section~\ref{sec:irreducible}. For general tropical fans however, the two definitions are different. Note that the irreducible components in this new definition might not induce a partition of $\Sigma_d$ as the following example shows.

\begin{example}[Irreducible components of non-normal fans] \label{ex:irreducible_components_not_well_defined}
Let $\Sigma = \tikz[scale=.18]\draw (0:-1)--(0:1) (45:-1.2)--(45:1.2) (90:-1)--(90:1);$ be the $1$-skeleton of $\Sigma_{U_{3,3}}$ in $\R^2$. There are five irreducible components in the sense of minimal support of a nonzero element of $\MW_1(\Sigma)$. These are $ \tikz[scale=.18, baseline=-5pt]\draw (0:-1)--(0:1);\,$, $\tikz[scale=.18]\draw (45:-1.2)--(45:1.2);\,$, $\tikz[scale=.18]\draw (90:-1)--(90:1);\,$, $\tikz[scale=.18]\draw (0:0)--(0:1) (45:-1)--(45:0) (90:0)--(90:1);$, and $\tikz[scale=.18]\draw (0:-1)--(0:0) (45:0)--(45:1) (90:-1)--(90:0);$. They do not induce a partition of $\Sigma_1$.
\end{example}

\subsection{The cross}

In this section, we consider the cross $\Delta=\tikz[scale=.18]\draw (0:-1)--(0:1) (90:-1)--(90:1);$ in $\R^2$ with four rays which is a reduced tropical fan. The cross is the simplest \emph{singular} tropical fan: it is neither irreducible nor div-faithful. In fact, it does not verify Poincaré duality. It can be used moreover in various ways in constructing interesting (counter-)examples. We propose some of them here.

\smallskip
In what follows, for $k\in \N$, $\Lambda^k$ is the complete fan in $\R^k$ with facets the $2^k$ orthants given by coordinate axes. The cross $\Delta$ is a subfan of $\Lambda^2$. We endow $\R^2$ with coordinates $x$ and $y$.

\begin{example}[A degenerate tropical modification and desingularization] \label{ex:nontrivial_degenerate_tropical_modification}
Let $f$ be the meromorphic function on $\Delta$ given by
\[f=
\begin{cases} -x & \textrm{ if } x\geq 0 \textrm{ and } y=0\\
\,y & \textrm{ if } x = 0 \textrm{ and } y\geq 0\\
\,0 & \textrm{ otherwise.}
\end{cases}
\]
Then, the divisor $\div(f)$ is trivial. The tropical modification $\~\Delta \coloneqq \tropmod{f}{\Delta}$ is a reduced tropical line. In particular, $\~\Delta$ verifies Poincaré duality, which is not the case for $\Delta$. Moreover $\MW_1(\~\Delta)$ is of dimension 1, though $\MW_1(\Delta)$ is of dimension 2. This shows that the div-faithfulness assumption is needed concerning Minkowski weights in Theorem \ref{thm:invariant_chow_tropical_modification}.

This example is particularly interesting: it seems that tropical modification tends to \emph{desingularize} tropical fans as blow-up does in algebraic geometry.
\end{example}

\begin{example} (A tropical modification which is no more saturated) \label{ex:nontrivial_linear_tropical_modification}
We use the notations of the previous example. Even though $\Delta$ is saturated, $\tropmod{2f}{\Delta}$ is not saturated at $\conezero$. Hence, saturation is not $\mT$-stable.
\end{example}

\begin{example}[Poincaré duality for a tropical fan does not imply Poincaré duality for its star fans]
The cross is the divisor of the meromorphic function $h$ on $\Lambda^2$ defined by $\min(0,x,y,x+y)$. Consider $\Sigma\coloneqq \tropmod{\Delta}{\Lambda^2}$. Since $\Lambda^2$ is div-faithful and saturated, by Theorem \ref{thm:invariant_chow_tropical_modification}, we have the equality of Chow rings $A^\bul(\Sigma) \simeq A^\bul(\Lambda^2)$. This implies that the Chow ring of $\Sigma$ verifies Poincaré duality. On the other hand, the cross appears as the star fan of the ray $\rho$ in $\Sigma$ corresponding to the tropical modification. We infer that Poincaré duality for a fan does not necessary imply Poincaré duality for its star fans.
\end{example}
\begin{example}[Poincaré duality for a non-normal fan] \label{ex:PD_chow_not_smooth}
We consider the previous example $\Sigma\coloneqq \tropmod{\Delta}{\Lambda^2}$. Note that it is not normal but its Chow ring verifies Poincaré duality. Moreover, $\Sigma$ is at the same time irreducible, principal and div-faithful at $\conezero$. However, it is neither locally irreducible nor div-faithful.
\end{example}

\begin{example}[Two non-principal unimodular tropical fans] \label{ex:non-principal_fan}
Consider the two-dimensional fan $\Sigma=\Delta\times\Lambda$ in $\R^3$, \ie, it has the six rays of the axes and support $\{x = 0\} \cup \{y = 0\} \subset \R^3$. The divisor $\{y = z = 0 \}$ in $\Sigma$ is not principal. This means $\Sigma$ is not principal at $\conezero$.

Following the idea of Example \ref{ex:PD_chow_not_smooth}, one can go further and create a fan which is principal at $\conezero$ but not globally principal: it suffices to take $\tropmod{\Sigma}{\Lambda^3}$.
\end{example}

\begin{example}[Non-irreducible unimodular fans]
We can find higher dimensional analogues of the cross in $\R^2$. Consider any unimodular fan $\Delta'$ with support $\{x_1 = x_2 = 0\} \cup \{x_3 = x_4 = 0\}$ in $\R^4$. The fan $\Sigma$ is normal, div-faithful and principal, but it is not irreducible.
\end{example}

\begin{example}[A normal div-faithful unimodular fan which is neither irreducible nor principal] \label{ex:div-faithful_not_irreducible}
Let $\Delta'$ be the fan defined in the previous example. Set $\Sigma = \Delta'\times\Lambda$. We get a tropical fan $\Sigma$ of dimension three which is still normal and div-faithful but it is neither irreducible, nor principal: the divisor $\{x_1 = x_2 = 0\}\times\{0\}$ is not a principal divisor in $\Sigma'$.
\end{example}

\begin{example}[A quasilinear fan which is not principal] \label{ex:q-principal}
Let $(\e_1,\e_2)$ be the standard basis of $\Z^2$. Let $\Sigma=\tikz[scale=.2, baseline=-4pt, transform shape]{\clip(0,0)circle(1);\fill[lightgray](0,0)circle(1);\draw(0,0)edge(1,0)edge(-1,-2)--(0,1);}$ be the complete fan in $\R^2$ with rays $\R_{\geq 0}\e_1, \R_{\geq 0}\e_2,$ and $\R_{\geq 0}(-\e_1-2\e_2)$. Let $(\Delta, \omega!_\Delta)$ be the generalized tropical line in $\R^2$ obtained as the divisor of the holomorphic function $f=\min(2x_1, x_2, 0)$ on $\Sigma$. Note that the orientation $\omega!_\Delta$ takes value 1 on $\R_{\geq 0}\e_1, \R_{\geq 0}(-\e_1-2\e_2),$ and value 2 on $\R_{\geq 0}\e_2$. The tropical modification $\~\Sigma$ of $\Sigma$ along the divisor of $f$ is quasilinear and it is principal. In $\~\Sigma$ we have a 2-dimensional cone $\sigma$ of weight~$2$. The unimodular stellar subdivision $\Sigma'$ of $\~\Sigma$ at cone $\sigma$ is not principal: the star of the new ray $\rho$ is isomorphic to the fan $\Lambda$ in $\R$ with weight 2.

This shows that principality is not $\mT$-stable in general. As stated in Theorem~\ref{thm:stability_principal}, principality is $\mT$-stable in the class of unitary tropical fans.
\end{example}

\subsection{The fan over the one-skeleton of the cube} \label{subsec:cube}

\begin{example}[A non-principal locally irreducible unimodular fan] \label{ex:cube}
Consider the standard cube $\mbox{\mancube}$ with vertices $(\pm1, \pm1, \pm1)$, and let $\Sigma$ be the two-dimensional fan with rays generated by the vertices, and with facets generated by the edges of the cube. The fan $\Sigma$ is locally irreducible and tropical but it is not unimodular. We obtain a unimodular fan after changing the lattice. In what follows, we work with the lattice $N \coloneqq \sum_{\zeta\in\Sigma_1}\Z\e_\zeta$.

Direct computation of the Chow ring proves that the image of $A^1(\Sigma)$ inside $A^1(\Sigma)^\dual$ is a sublattice of full rank and of index two. The map $A^1(\Sigma) \to A^1(\Sigma)^\dual$ is not surjective, hence $\Sigma$ is not principal at $\conezero$ though it is $\Q$-principal. Indeed, the divisor $D=\R\cdot(1,1,1)$ is not principal, but $2D$ is.
\end{example}

\subsection{Saturation, unimodularity, and desingularization} \label{subsec:examples_non_saturated_non_unimodular}
Let $(\e_1,\e_2)$ be the standard basis of $\Z^2$ and let $\e_0=-\e_1-\e_2$. Denote by $\rho_i = \R_{\geq 0}\e_i$, $i\in\{0,1,2\}$ the corresponding rays.

\smallskip
In the following, we consider the lattice $N \coloneqq \Z\e_1 + \frac 13\Z(\e_1 - \e_2)$ in $\R^2$. Note that $\Z^2 \subset N$ is a sublattice of index three in $N$. The dual lattice $M \coloneqq N^\dual$ is of index three in $(\Z^2)^\dual$. The fans $\Delta$ and $\Sigma$ treated in the examples of this section verify $\ssN_{\Delta} = \ssN_\Sigma=N$.

\begin{example} \label{ex:F1_vs_N}
Let $\Delta$ be the one-dimensional fan in $N_\R$ with rays $\rho_0, \rho_1, \rho_2$. The fan $\Delta$ is tropical and unimodular but it is not saturated. Any element $f$ of $(\Z^2)^\dual \setminus M$ induces a meromorphic function on $\Delta$ which is linear but not integral linear. (The meromorphic function $3f$ will be integral linear.) As a consequence, $A^1(\Delta)$ has torsion: the element $x_{\rho_1}-x_{\rho_2}$ is non-zero, but $3\cdot(x_{\rho_1}-x_{\rho_2})$ vanishes in $A^1(\Delta)$.
\end{example}

\begin{example}[Necessity of assumptions in Theorem~\ref{thm:invariant_chow_tropical_modification}] \label{ex:necessity_invariant_chow_tropical_modification}
Consider the fan $\Delta$ defined in the previous example. Let $f$ be the element in $(\Z^2)^\dual$ which takes value one on $\e_1$ and vanishes on $\e_2$. The divisor $\div(f)$ is trivial. Let $\~\Delta$ be the degenerate tropical modification of the fan $\Sigma$ with respect to $f$. Then, $\~\Delta$ is just a tropical line. Therefore, it is saturated. In particular, the Chow ring of $\~\Delta$ is torsion-free. This example shows that the hypotheses made in Theorem~\ref{thm:invariant_chow_tropical_modification} about the saturation or torsion-freeness are needed.

As in Example~\ref{ex:nontrivial_degenerate_tropical_modification}, tropical modification here desingularizes the original tropical fan.
\end{example}

\begin{example} \label{ex:non_unimodular_projective_fan}
Let $\Sigma$ be the complete fan in $N_\R$ with rays $\rho_0$, $\rho_1,$ and $\rho_2$. The fan $\Sigma$ is tropical and saturated but it is not unimodular. As in Example~\ref{ex:F1_vs_N}, $A^1(\Sigma)$ has torsion.
\end{example}

\begin{example} \label{ex:non_saturated_smooth_fan}
We continue with the complete fan $\Sigma$ of the previous example. Let $f$ be the holomorphic function on $\Sigma$ which takes value $-3$ on $\e_0$, and vanishes on $\e_1$ and $\e_2$. The divisor of $f$ is the reduced divisor $\Delta =\div(f)$ of Example~\ref{ex:F1_vs_N} with rays $\rho_0$, $\rho_1$, and $\rho_2$.

Let $\Sigma'$ be a unimodular subdivision of $\Sigma$. For instance, we can add the rays generated by the primitive vectors $\tfrac13\e_i+\tfrac23\e_j$ for any pair of distinct $i,j \in \{0,1,2\}$.

We view $f$ as a holomorphic function on $\Sigma'$, and set $\~\Sigma \coloneqq \tropmod{f}{\Sigma'}$, the tropical modification of $\Sigma'$ along the divisor $\Delta=\div(f)$. Denote by $\rho$ the new ray in $\~\Sigma$. Then, $\~\Sigma$ is a unimodular fan that is saturated at $\conezero$ but not at $\rho$. Even worse, there is no way to modify the lattice $N$ in order to make $\~\Sigma$ saturated without changing the intersection $\~\Sigma \cap N$.

Since $\Sigma'$ is a complete unimodular fan, it is div-faithful and its Chow ring verifies Poincaré duality. Applying Theorem \ref{thm:invariant_chow_tropical_modification}, we infer that $A^\bul(\~\Sigma)$ verifies Poincaré duality. However, the Chow ring of the star fan $\~\Sigma^\rho$ has torsion. In $A^1(\~\Sigma^\rho)$, the element $\i^*_{\conezero\subface\rho}(x_{\rho_1}-x_{\rho_2})$ has order 3. (The image by the Gysin map of this element is trivial in $A^2(\~\Sigma)$.)

Once again, as in Examples~\ref{ex:nontrivial_degenerate_tropical_modification} and~\ref{ex:necessity_invariant_chow_tropical_modification}, a tropical modification can solve these \emph{singularity} issues. Consider the conewise linear function $h$ on $\~\Sigma$ defined by
\[ h = \min(\e_1^\dual,\e_\rho^\dual), \]
where $(\e_1^\dual,\e_2^\dual,\e_\rho^\dual)$ denotes the dual (rational) basis of $(\e_1,\e_2, \e_\rho)$. Since $\e_\rho^\dual$ coincides with $f$ on the graph of $f$, it is straightforward to see that $h$ takes integer values on $\e_0, \e_1, \e_2$, and on any vector of the form $\tfrac13\e_i+\tfrac23\e_j$, for distinct $i,j \in \{0,1,2\}$. This means that $h$ is holomorphic on a suitable subdivision of $\~\Sigma$. Consider now the tropical modification $\widehat \Sigma$ of $\~\Sigma$ with respect to $h$ on $\~\Sigma$. The fan $\widehat \Sigma$ is saturated and verifies $\Ploc\PD_\Z$, that is, all the star fans of $\widehat \Sigma$ verify Poincaré duality.
\end{example}

In view of the preceding examples, we formulate the following questions.

\begin{question}[Saturation and tropical modification]
Is it true that for any tropical fan $\Sigma$, there exists a sequence of tropical modifications that turns $\Sigma$ into a saturated tropical fan?
\end{question}

\begin{question}[Saturation and torsion-freeness]
Is it true in general that the Chow ring $A^\bul(\Sigma, \Z)$ with integer coefficients of a saturated unimodular fan $\Sigma$ is torsion-free?
\end{question}

\subsection{Quasilinear vs generalized Bergman}

\begin{example}[A quasilinear tropical fan which is not generalized Bergman] \label{ex:stable_not_Bergman}
Let $\Sigma$ be the two-dimensional skeleton of the fan of the projective space of dimension three: it has four rays generated by vectors $\e_1, \e_2, \e_3$ and $\e_0= -\e_1-\e_2-\e_3$, where $(\e_1,\e_2,\e_3)$ is a basis of $N$, and the six facets $\R_{\geq 0}\e_i+\R_{\geq 0}\e_j$ for $0\leq i<j\leq 3$ (\cf Figure \ref{fig:Bergman_fan_U34} where we keep only four rays). We have $\supp\Sigma = \supp{\ssSigma_{\Ma}}$ for the uniform matroid $\Ma = U^3_{4}$ of rank three on four elements. Let $\Delta$ be the tropical curve in $\supp\Sigma$ with four rays $\rho_1, \rho_2, \rho_3$ and $\rho_4$ with
\[ \rho_1 =\R_{\geq 0} (\e_1+\e_2), \qquad \rho_2 =\R_{\geq 0} \e_2, \qquad \rho_3 = \R_{\geq 0}(\e_0+ 2 \e_3), \qquad \rho_4 = \R_{\geq 0} (\e_0 + \e_1). \]
Then we have $\Delta \simeq \Sigma_{\Ma'}$ with $\Ma' = U^2_4$ the uniform matroid of rank two on four elements. The tropical modification $\tropmod\Delta\Sigma$ is quasilinear but does not have the same support as the Bergman fan of any matroid. It means it is not generalized Bergman.
\end{example}

\subsection{Collection of examples around the Kähler package for the Chow ring}

In this section, all the Chow rings are with rational coefficients.

\begin{example}
Let $\Sigma$ be the complete Bergman fan associated to the uniform matroid $U_{3,3}$ on the ground set $\{0,1,2\}$ as illustrated in Figure~\ref{fig:Bergman_fan_U33}.
\begin{figure}[ht]
\caption{The Bergman fan of $U_{3,3}$} \label{fig:Bergman_fan_U33}
\begin{tikzpicture}
\clip (0, 0) circle (2);
\fill[color=gray!10] (-2,-2) rectangle (2,2);
\draw[color=gray!25, scale=.6] (-4,-4) grid (4,4);
\draw (0 , 0)
  edge["$\rho_1$", near end]     (2 , 0)
  edge["$\rho_{12}$"{above left=-4pt}]  (2 , 2)
  edge["$\rho_{2}$", near end]   (0 , 2)
  edge["$\rho_{02}$", near end]  (-2, 0)
  edge["$\rho_{0}$"{below right=-4pt}]   (-2,-2)
  edge["$\rho_{01}$", near end]  (0 ,-2);
\end{tikzpicture}
\end{figure}
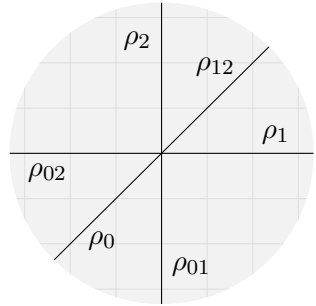
The nonempty proper flats of $U_{3,3}$ are $\prCl(U_{3,3})\coloneqq\{1,12,2,02,0,01\}$ (where $ij$ denotes the set $\{i,j\}$). The rays of $\Sigma$ correspond to the flats of $U_{3,3}$ and are denoted $\rho_F$ with $F\in\prCl(U_{3,3})$. These rays corresponds to elements of $A^1(\Sigma)$ denoted $x_F$, for $F\in\prCl(U_{3,3})$. We have
\[ \dim(A^1(\Sigma))=\card{\prCl(U_{3,3})}-\dim\bigl((\Q^2)^\dual\bigr)=4. \]
We have a natural degree map $\deg\colon A^2(\Sigma)\to\Q$ which verifies the following.
\[ \deg(x_Fx_G)=\begin{cases}
  -1 & \text{if $F=G$,} \\
  1  & \text{if $F\neq G$ and $F \subset G$ or $F \supset G$,} \\
  0  & \text{if $F$ and $G$ are not comparable.}
\end{cases} \]

\smallskip
Let $\ell\coloneqq x_{1}+x_{12}+x_{2}+x_{02}+x_{0}+x_{01}\in A^1(\Sigma)$. This element corresponds to a strictly convex conewise linear function on $\Sigma$, thus, it must verify the Hodge-Riemann bilinear relations $\HR(\Sigma,\ell)$. We check that this is indeed the case. We have to show that $\deg(\ell^2)>0$ and that
\[ \begin{array}{rccc}
Q_1\colon & A^1(\Sigma)  \times A^1(\Sigma) & \to     & \Q,       \\
          & x            ,       y          & \mapsto & \deg(xy),
\end{array} \]
has signature $-2=2\dim(A^0)-\dim(A^1)$. We have
\[\deg(\ell x_1)=\deg(x_1^2+x_1x_{12}+x_1x_{01})=1, \]
and, by symmetry, $\deg(\ell x_F)=1$ for every $F\in\prCl(U_{3,3})$. Thus, $\deg(\ell^2)=6$.

For $Q_1$ we have an orthogonal basis $(x_0, x_1, x_2, x_1+x_{12}+x_2)$. Moreover $\deg(x_0^2)=\deg(x_1^2)=\deg(x_2^2)=-1$ and

\begin{equation}\label{eq:degree_one}
  \deg((x_1+x_{12}+x_2)^2)=\deg(x_1^2+x_{12}^2+x_2^2+2x_1x_{12}+2x_{12}x_2)=1.
\end{equation}

Thus, the signature of $Q_1$ is $-2$ and the Hodge-Riemann bilinear relations are verified.
\end{example}

\begin{example}
We give a second example in the Chow ring treated in the preceding example. Let $\ell'=x_1+x_{12}+x_2$. The conewise linear function associated to $\ell'$ is not strictly convex. Around $\rho_0$, the function is zero which implies that $\i^*_{\conezero\ssubface\rho_0}(\ell')=0$. Nevertheless, $\ell'$ verifies the Hodge-Riemann bilinear relations $\HR(\Sigma, \ell')$. We already checked in \eqref{eq:degree_one} that $\deg(\ell'^2) = 1 >0$. Since $Q_1$ does not depend on $\ell'$, we have already checked that $Q_1$ has signature $-2$. This gives an example of a non-ample element verifying the Hodge-Riemann bilinear relations.
\end{example}

\begin{example}[Non-convex functions verifying $\HL$ and $\HR$] \label{ex:non_convex_HR}
For any tropical fan $\Sigma$, the set of elements of $A^1(\Sigma)$ which verify $\HL$ is either empty or it is the complement of a finite union of hypersurfaces (defined by the vanishing of the determinant of the linear map $A^k(\Sigma) \to A^{d-k}(\Sigma)$, $k\leq d/2$, given by the multiplication by the $(d-2k)$-th power of the element in $A^1(\Sigma)$). The essential part of the present example is thus to describe a non-convex conewise linear function on a tropical fan $\Sigma$ whose associated element in $A^1(\Sigma)$ verifies $\HR$. Note that such an element does not exist in the case of complete fans, \cf Remark \ref{rk:HR_complete_fan}.

\begin{figure}[ht]
  \caption{The Bergman fan of $U_{3,4}$} \label{fig:Bergman_fan_U34}
  \begin{tikzpicture}[scale = 1.4]
    \let\i\relax

    \BFcoordinates

    \foreach \i in {1,2,3,0} {
      \draw (I) -- (I\i);
      \draw ($(I)!1.1!(I\i)$) node {\i};
    }
    \foreach \i/\j/\cc in {1/3/blue, 2/3/green, 1/0/black, 3/0/violet, 1/2/red, 2/0/brown} {
      \draw (I) -- (I\i\j);
      \draw ($(I)!1.1!(I\i\j)$) node {\ifnumequal{\j}{0}{\j\i}{\i\j}};
      \roundface\i\j[\cc, opacity=.1];
    }
  \end{tikzpicture}
\end{figure}
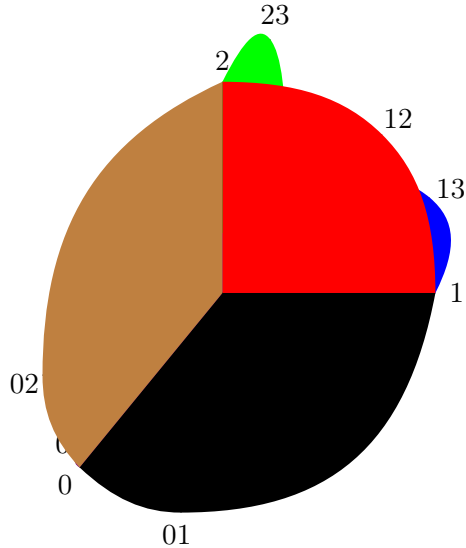
Take the fan $\Sigma$ associated to the matroid $U_{3,4}$ (\cf Figure \ref{fig:Bergman_fan_U34}) where we keep only the rays labeled by the flats $0,1,2,3,01,23$. Let $(\e_1^\dual, \e_2^\dual, \e_3^\dual)$ be the dual basis of $(\e_1, \e_2, \e_3)$. Let $f$ be the restriction to $\Sigma$ of the function
\[ \max(0,-\e_1^\dual,\e_3^\dual-\e_2^\dual,\e_3^\dual-\e_2^\dual-\e_1^\dual). \]
The following table gives the values of $f$ and of $\div(f)$:
\[ \begin{array}{r||c|c|c|c|c|c}
  F                  & 0 & 1 & 2 & 3 & 01 & 23 \\\hline
  f(\e_F)^{\phantom|}& 1 & 0 & 0 & 1 & 0  & 0  \\
  \div(f)(\zeta_F) & -1 & -1 & -1 & -1 & -1  & -1
\end{array} \]
It follows that $-f$ is holomorphic on $\Sigma$. From its definition, we get that $f$ is convex. Moreover, the values of $\div(f)$ taken on rays show that $f$ is strictly convex around each ray. In particular, Proposition \ref{prop:local_HR} ensures that we have $\HL(\Sigma, f)$. Since $f$ is convex, it is a limit of ample elements, and by continuity we also deduce $\HR(\Sigma, f)$. However, $f$ is not strictly convex around~$\conezero$ since $f(\e_{01}+\e_{23}) = f(\e_{01}) + f(\e_{23}) = 0$.

Since $\HL$ and $\HR$ are open conditions, we can find a non-convex function in the neighborhood of $f$ which also verifies $\HL$ and $\HR$.
\end{example}

\begin{remark}[Comparison with the case of a complete fan] \label{rk:HR_complete_fan}
The phenomenon discussed in the previous example cannot occur in the case of a complete fan. Let $\Sigma$ be a complete fan and $f$ be a meromorphic function on $\Sigma$.
\begin{thm}
The following statements are equivalent:
\begin{enumerate}
  \item $f$ is strictly convex on $\Sigma$.
  \item the properties $\HL(\Sigma^\sigma, f^\sigma)$ and $\HR(\Sigma^\sigma, f^\sigma)$ hold for all $\sigma\in\Sigma$.
  \item $\div(-f)$ is effective and has full support $\Sigma_{d-1}$.
\end{enumerate}
\end{thm}
We omit the proof.
\end{remark}

\begin{example}[A convex function which does not come from a convex function on the ambient space]
The difference of behavior is to be compared with the fact that a convex function on a non-complete fan might not be extendable to a convex function on the ambient space. Let $\Sigma$ and $f$ be as in Example \ref{ex:non_convex_HR}. Consider the ray $\rho_{001}$ obtained from the unimodular blow-up of the cone between $\rho_{0}$ and $\rho_{01}$, and denote by $\Sigma'$ the new tropical fan. Let $g$ be a holomorphic function on $\Sigma'$ which coincides with $f$ on all the rays of the original fan and takes a value on $\e_{001}$ which is slightly smaller that $f(\e_{001})$. Then, $g$ is still convex but we have
\[ g(\e_{001}) + g(\e_{23}) < 1 = g(e_0) = g(\e_{001} + \e_{23}). \]
This means that $g$ cannot come from a convex function on the ambient space. It is possible to turn $g$ into a strictly convex function verifying the same non-extendability property, for instance by slightly increasing the value of $g$ on $\e_{01}$.
\end{example}

\begin{example}[An ample element which does not verify $\HR$]
We study here an interesting example discovered by Babaee and Huh \cite{BH17}. We would like to thank Edvard Aksnes and Kris Shaw for drawing our attention to the relevance of this example.

A full description of the fan is given in \cites{BH17, Aks19, Piq-thesis}, and we will study it more thoroughly in our future work which extends the present paper to generalization of fans. We just mention its main properties.

This is a normal unimodular tropical fan $\Sigma$ of dimension two living in $\R^4$. Its Chow ring verifies Poincaré duality. However, $\Sigma$ does not verify the Hodge-Riemann bilinear relations: $\Sigma$ is quasi-projective, but the pairing $A^1(\Sigma) \times A^1(\Sigma) \to \Z$ has more than one positive eigenvalues. In particular, we note that $\Sigma$ is not quasilinear.
\end{example}

\begin{example}[Positivity for tropical fans whose orientation takes negative weights]
There are examples of tropical fans whose orientation takes negative values that verify an analogue of the Kähler package for the Chow ring: it will be enough to take a Kähler tropical fan $(\Sigma, \omega!_\Sigma)$ with positive orientation $\omega!_\Sigma$, and consider the fan $(\Sigma, -\omega!_\Sigma)$. A generalized tropical line with both positive and negative weights verifies as well the Kähler package for the Chow ring. It will be interesting to formulate an appropriate notion of Kähler package for tropical fans whose orientations take both positive and negative values. This requires formulating the right notion of convexity in this setting.
\end{example}

\appendix

\section{Multi-magmoids} \label{sec:multimagmoid}

\newcommand{\ee}{\underline{e}}
\newcommand{\ff}{\underline{f}}
\newcommand{\Set}{\textrm{Set}}
\newcommand{\op}{\mathrm{op}}
\newcommand{\gen}[1]{\langle #1 \rangle}
\newcommand{\M}{\mathrm{C}}
\newcommand{\cN}{\mathrm{D}}
\newcommand{\abs}{\dims}
\newcommand{\rmS}{\mathrm{S}}
\newcommand{\rmB}{\mathrm{B}}
\newcommand{\rmI}{\mathrm{I}}
\tossub\op

In this appendix, we define an algebraic structure called \emph{multi-magmoid} in order to give a more conceptual approach to the notion of $\mT$-stability introduced in Section \ref{sec:operations}. We also prove Proposition \ref{prop:intersection-property}.

A \emph{multi-magmoid} $\M$ is the data of a ground set, also denoted $\M$, and of a finite set of multi-valued binary operators,
\[\op!_i = \op!_{\M, i}\colon \M \times \M \to 2^\M \quad \textrm{ for } i=1, \dots, k\]
for a positive integer $k$. If $k = 1$ and if $\op!_1(c,c')$ is a singleton for all $c$ and $c'$ in $\M$, then $\M$ is a \emph{magma}.

Note that if $\cN$ is a subset of $\M$, then the multi-magmoid structure on $\M$ induces a multi-magmoid structure on $\cN$ with operators
\[ \begin{array}{rccl}
  \op!_{\cN,i}\colon & \cN \times \cN & \to & 2^\cN, \\
  & d, d' & \mapsto & \op!_{\M,i}(d,d') \cap \cN.
\end{array} \]

A subset $\rmS$ of $\M$ is called \emph{$\M$-stable} if for any $i \in [k]$ and any $s,s' \in \rmS$, $\op!_i(s,s') \subseteq \rmS$. If $\rmS$ is included in a subset $\cN$ of $\M$, we say that $\rmS$ is \emph{$\M$-stable in $\cN$} if it is $\cN$-stable for the induced multimagmoid on $\cN$. If $\rmB$ is a subset of $\M$, we denote by $\gst{\rmB}{\M}$ the smallest $\M$-stable subset of $\M$ that contains $B$ and called it the \emph{$\M$-stable subset generated by $\rmB$}. This exists by Proposition~\ref{prop:basic_properties_multi-magmoid}.

A subset $\rmI$ of $\M$ is called an \emph{ideal of $\M$} if for any $i \in [k]$, and for any $c \in \M$ and any $j \in \rmI$, $\op!_i(j,c) \subseteq \rmI$ and $\op!_i(c,j) \subseteq \rmI$.

A subset $\rmS$ of $\M$ is called \emph{strongly $\M$-stable} if it is $\M$-stable and if $\M \setminus \rmS$ is an ideal of $\M$. Equivalently, this means that $\rmS$ is $\M$-stable and, for $c,c' \in \M$ and for $i \in [k]$, $\op!_i(c,c') \cap \rmS \neq \emptyset$ implies that both $c$ and $c'$ belongs to $\rmS$.

\begin{prop} \label{prop:basic_properties_multi-magmoid}
We have the following basic properties.
\begin{enumerate}
\item The intersection of two $\M$-stable subsets of\/ $\M$ is $\M$-stable.
\item The intersection of two ideals of\/ $\M$ is an ideal.
\item Ideals of\/ $\M$ are $\M$-stable.
\item If\/ $\rmI$ is an ideal of\/ $\M$ and\/ $\rmS$ is a $\M$-stable subset of\/ $\M$, then $\rmI \cup \rmS$ is stable.
\item If\/ $\rmS$ is strongly $\M$-stable, then for any $\rmB$, we have $\gst{\rmB \cap \rmS}{\M} = \gst{\rmB}{\M} \cap \rmS$.
\end{enumerate}
\end{prop}

\begin{proof}
The first four points are immediate. We prove the last point. Since $\rmS$ is strongly $\M$-stable, $\M \setminus \rmS$ is an ideal. Point (4) implies that $\rmS' \coloneqq (\M \setminus \rmS) \cup \gst{\rmB \cap \rmS}{\M}$ is stable. Moreover, since $\gst{\rmB \cap \rmS}{\M} \subseteq \rmS$, this is a disjoint union. In addition, $\rmS'$ contains $\rmB$, therefore, it contains $\gst{\rmB}{\M}$. We infer that $\gst{\rmB \cap \rmS}{\M} = \rmS'\cap \rmS \supseteq \gst{\rmB}{\M} \cap \rmS$. The other inclusion is trivial, and we get the equality.
\end{proof}

\smallskip
We now make explicit the link between $\mT$-stability introduced in Section \ref{sec:operations} and multi-magmoids. Let $\mT$ be the set of tropical fans up to isomorphism. We assume that a formal empty fan belongs to $\mT$ in order to allow the empty divisor for the tropical modification. We define the following binary operators on $\mT$.
\begin{itemize}
\item the cartesian product
\[ \op!_1(\Sigma, \Sigma') \coloneqq \{ \Sigma \times \Sigma' \}. \]
\item the tropical modification
\[ \op!_2(\Sigma, \Delta) \coloneqq \left\{ \tropmod{f}{\Sigma} \st f \in \mer(\Sigma), \div(f) = \Delta \right\}. \]
\item the blow-up $\op!_3(\Sigma, \Delta)$ defined as the set of blow-ups $\Sigma_{(\rho)}$ of $\Sigma$ along a ray $\rho$ that lives in the interior of a face $\sigma$ of $\Sigma$ such that we have $\Sigma^\sigma = \Delta$.
\item the blow-down
\[ \op!_4(\Sigma', \Delta) \coloneqq \left\{ \Sigma \st \Sigma' \in \op!_3(\Sigma, \Delta) \right\}. \]
\end{itemize}

\begin{prop}
The set $\mT$ with the four multi-valued binary operators $\op!_1, \op!_2, \op!_3, \op!_4$ defined above is a multi-magmoid.
\end{prop}

If $\Csh$ is a subset of $\mT$, then the $\mT$-stability in $\Csh$ introduced in Section \ref{sec:operations} coincide with the notion introduced in this section. Also note that the class of effective, \resp reduced, \resp unitary, tropical fans are strongly $\mT$-stable. Together with Proposition \ref{prop:basic_properties_multi-magmoid}, this implies Proposition~\ref{prop:intersection-property}.

\let\v\cech
\bibliographystyle{alpha}
\bibliography{$HOME/bibliography/bibliography}
\end{document}